\documentclass[twoside,reqno,11pt]{amsart}
\usepackage{a4wide}
\usepackage{amsmath}
\usepackage{amsfonts}
\usepackage{amssymb}
\usepackage{times}
\usepackage{graphicx}
\usepackage{color}
\usepackage[color,arrow, matrix, curve]{xy}
\usepackage{pgf,tikz}
\usepackage{etex}
\usepackage{tikz}
\usetikzlibrary{arrows,backgrounds}
\usetikzlibrary{decorations.pathreplacing}
\usetikzlibrary{decorations.markings}
\usetikzlibrary{shapes}
\newenvironment{pic}[1][]%
{\begin{aligned}\begin{tikzpicture}[#1]}%
{\end{tikzpicture}\end{aligned}}
\tikzstyle{string}=[line width=1.25pt]
\tikzstyle{thickstring}=[line width=2.5pt]
\tikzstyle{dot}=[circle, draw=black, fill=gray, inner sep=.4ex, line width=1.25pt]
\tikzstyle{whitedot}=[circle, draw=black, fill=white, inner sep=.4ex, line width=1.25pt]
\tikzstyle{blackdot}=[circle, draw=black, fill=black, inner sep=.4ex, line width=1.25pt]
\tikzstyle{cross}=[preaction={draw=white, -, line width=6pt}]
\tikzstyle{coupon}=[rectangle,draw=black,minimum size=10pt,line width=1.25pt]
\tikzstyle{unit}=[triangle,line width=1.25pt,fill=white,draw,inner sep=1pt,minimum width=1.1cm]
\tikzstyle{counit}=[triangle,hflip,line width=1.25pt,fill=white,draw,inner sep=1pt,minimum width=1.1cm]
\tikzstyle{iso2cell}=[iso,line width=1.25pt,fill=white,draw,inner sep=1pt,minimum width=1.1cm]
\tikzstyle{pullback2cell}=[iso,pullback,line width=1.25pt,fill=white,draw,inner sep=1pt,minimum width=1.1cm]
\tikzstyle{2cell}=[rectangle,line width=1.25pt,fill=white,draw,inner sep=1pt,minimum size=10pt,minimum width=1.1cm]
\tikzset{arrow/.style={decoration={
    markings,
    mark=at position #1 with \arrow{>}},
    postaction=decorate}
}
\tikzset{reverse arrow/.style={decoration={
    markings,
    mark=at position #1 with \arrow{<}},
    postaction=decorate}
}
\makeatletter
\newif\ifpullback
\pgfkeys{/tikz/pullback/.is if=pullback}
\pgfkeys{/tikz/pullback=false}
\pgfdeclareshape{iso}{%
    \inheritsavedanchors[from=rectangle]
    \inheritanchorborder[from=rectangle]
    \inheritanchor[from=rectangle]{north}
    \inheritanchor[from=rectangle]{north west}
    \inheritanchor[from=rectangle]{north east}
    \inheritanchor[from=rectangle]{center}
    \inheritanchor[from=rectangle]{west}
    \inheritanchor[from=rectangle]{east}
    \inheritanchor[from=rectangle]{mid}
    \inheritanchor[from=rectangle]{mid west}
    \inheritanchor[from=rectangle]{mid east}
    \inheritanchor[from=rectangle]{base}
    \inheritanchor[from=rectangle]{base west}
    \inheritanchor[from=rectangle]{base east}
    \inheritanchor[from=rectangle]{south}
    \inheritanchor[from=rectangle]{south west}
    \inheritanchor[from=rectangle]{south east}
    \savedanchor\centerpoint
    {
        \pgf@x=0pt
        \pgf@y=0pt
    }
    \anchor{text}{%
         \pgf@x=0pt
         \pgf@y=0pt
         \advance\pgf@y by -0.2\ht\pgfnodeparttextbox
    }
    \backgroundpath{%
      \northeast \pgf@xa=\pgf@x \pgf@ya=\pgf@y
      \southwest \pgf@xb=\pgf@x \pgf@yb=\pgf@y
      \begingroup
         \tikz@mode
         \iftikz@mode@draw
            \pgfpathmoveto\northeast
            \pgfpathlineto{\pgfqpoint{\pgf@xb}{\pgf@ya}}
            \pgfpathmoveto\southwest
            \pgfpathlineto{\pgfqpoint{\pgf@xa}{\pgf@yb}}
            \pgfusepath{stroke}
         \fi
         \iftikz@mode@fill
            \pgfsetlinewidth{\pgflinewidth}
            \advance\pgf@ya by -0.5\pgflinewidth%
            \advance\pgf@yb by  0.5\pgflinewidth%
            \pgfpathrectanglecorners{\pgfqpoint{\pgf@xb}{\pgf@yb}}{\pgfqpoint{\pgf@xa}{\pgf@ya}}
            \pgfusepath{fill}
         \fi
         \ifpullback
            \northeast \pgf@xa=\pgf@x \pgf@ya=\pgf@y
            \southwest \pgf@xb=\pgf@x \pgf@yb=\pgf@y
            \pgf@xc=\pgf@xb
            \advance\pgf@xc by  3\pgflinewidth
            \pgf@yc=\pgf@ya
            \advance\pgf@yc by -2\pgflinewidth
            \pgfpathmoveto{\pgfqpoint{\pgf@xc}{\pgf@ya}}
            \pgfpathlineto{\pgfqpoint{\pgf@xc}{\pgf@yc}}
            \pgfpathlineto{\pgfqpoint{\pgf@xb}{\pgf@yc}}
            \pgfusepath{stroke}
         \fi
    \endgroup
    }
}
\newif\ifhflip
\pgfkeys{/tikz/hflip/.is if=hflip}
\pgfkeys{/tikz/hflip=false}
\pgfdeclareshape{triangle}{%
     \saveddimen{\halfbaselength}{%
         \pgf@x=0.5\wd\pgfnodeparttextbox
         \pgfmathsetlength\pgf@xc{\pgfkeysvalueof{/pgf/inner xsep}}%
         \advance\pgf@x by \pgf@xc%
         \advance\pgf@x by \ht\pgfnodeparttextbox
         \pgfmathsetlength\pgf@xb{\pgfkeysvalueof{/pgf/minimum width}}%
         \divide\pgf@xb by 2
         \ifdim\pgf@x<\pgf@xb%
             \pgf@x=\pgf@xb%
         \fi%
     }
     \saveddimen\triangleheight{%
         \pgfmathsetlength\pgf@xc{\pgfkeysvalueof{/pgf/inner ysep}}%
         \advance\pgf@x by \pgf@xc%
         \pgfmathsetlength\pgf@xb{\pgfkeysvalueof{/pgf/minimum height}}%
         \divide\pgf@xb by 2
         \ifdim\pgf@x<\pgf@xb%
             \pgf@x=\pgf@xb%
         \fi%
     }
     \savedanchor\centerpoint{%  midpoint on base line
         \pgf@x=0pt 
         \pgf@y=0pt
     }
     \anchor{center}{\centerpoint}
     \anchor{text}{%
         \pgf@x=-0.5\wd\pgfnodeparttextbox
         \ifhflip
            \pgf@y=-1.2\ht\pgfnodeparttextbox
            \advance\pgf@y by \dp\pgfnodeparttextbox
            \advance\pgf@y by -3pt
         \else
            \pgf@y=\dp\pgfnodeparttextbox
            \advance\pgf@y by -\dp\pgfnodeparttextbox
            \advance\pgf@y by 4pt
         \fi
     }
     \anchor{left}{%
        \pgf@x=-\halfbaselength
        \pgf@y=0pt
     }
     \anchor{right}{%
        \pgf@x=\halfbaselength
        \pgf@y=0pt
     }
     \anchor{tip}{%
        \pgf@x=0pt
        \ifhflip
           \pgf@y=-\triangleheight
        \else
           \pgf@y=\triangleheight
        \fi 
     }
     \anchor{a}{%
        \pgf@x=-\halfbaselength
        \divide\pgf@x by 2
        \pgf@y=0pt
     }
     \anchor{b}{%
        \pgf@x=\halfbaselength
        \divide\pgf@x by 2
        \pgf@y=0pt
     }
     \backgroundpath{%
 	    \pgf@xa=\halfbaselength
            \pgf@ya=\triangleheight 
 	    \ifhflip
            \pgfpathmoveto{\pgfqpoint{0pt}{-\pgf@ya}}
	        \pgfpathlineto{\pgfqpoint{-\pgf@xa}{0pt}}
	        \pgfpathlineto{\pgfqpoint{\pgf@xa}{0pt}}
	        \pgfpathclose
	    \else 
            \pgfpathmoveto{\pgfqpoint{0pt}{\pgf@ya}}
	        \pgfpathlineto{\pgfqpoint{-\pgf@xa}{0pt}}
	        \pgfpathlineto{\pgfqpoint{\pgf@xa}{0pt}}
	        \pgfpathclose
            \fi
     }
     \anchorborder{%
		\pgfkeysgetvalue{/pgf/shape border rotate}{\rotate}%
		\edef\externalx{\the\pgf@x}%
		\edef\externaly{\the\pgf@y}%
		\centerpoint%
		\pgf@xa\externalx\relax%
		\pgf@ya\externaly\relax%
		\advance\pgf@xa\pgf@x%
		\advance\pgf@ya\pgf@y%
		\edef\externalx{\the\pgf@xa}%
		\edef\externaly{\the\pgf@ya}%
		\pgfmathanglebetweenpoints{\centerpoint}{\pgfqpoint{\externalx}{\externaly}}%
                \pgfmathsubtract@{\pgfmathresult}{\rotate}%
		\ifdim\pgfmathresult pt<0pt\relax%
			\pgfmathadd@{\pgfmathresult}{360}%
		\fi%
		\let\externalangle\pgfmathresult%
        \pgf@xc=-\halfbaselength
        \pgf@yc=0pt
		\pgfmathdivide@{\externalangle}{50}
		\pgfmathparse{\halfbaselength*\pgfmathresult}
		\advance\pgf@xc by \pgfmathresult pt%
		\pgf@y=0pt
		\pgf@x=\pgf@xc
        }
}
\makeatother
\def\GL{\mathsf{GL}}
\def\CGL{\mathsf{CGL}}  % Ron's version is : CharGL
\def\RGL{\mathsf{RGL}}
\def\H{\mathsf{H}}
\def\HB{\mathsf{H}^\bullet}
\def\HBI#1{\mathsf{H}^{\bullet,#1}}
\def\O{\mathsf{O}}
\def\CO{\mathsf{CO}}
\def\SL{\mathsf{SL}}
\def\Sp{\mathsf{Sp}}
\def\CSp{\mathsf{CSp}}  % iff Ron's version of \CGL then CharSp
\def\U{\mathsf{U}}
\def\mul{\mathsf{m}}
\def\comul{\Delta}
\def\antip{\mathsf{S}}
\def\sw{\mathsf{sw}}
\def\conv{\star}
\def\unit{\mathsf{e}}
\def\Symm{\mathsf{Sym}}
\def\sfA{\mathsf{A}}
\def\sfa{\mathsf{a}}
\def\sfoa{\overline{\mathsf{a}}}
\def\sfb{\mathsf{b}}
\def\sfB{\mathsf{B}}
\def\sfc{\mathsf{c}}
\def\sfd{\mathsf{d}}
\def\sfoc{\overline{\mathsf{c}}}
\def\sfe{\mathsf{e}}
\def\sff{\mathsf{f}}
\def\sfof{\overline{\mathsf{f}}}
\def\sfg{\mathsf{g}}
\def\sfu{\mathsf{u}}
\def\sfv{\mathsf{v}}
\def\ox{\overline{x}}
\def\oz{\overline{z}}
\def\End{\textsf{End}}
\def\Sym{\textsf{Sym}}
\def\NSym{\textsf{NSym}}
\def\ot{\otimes}
\def\Id{\mathrm{I\kern-0.2ex d}}
\def\la{\langle}
\def\ra{\rangle}
\def\lla{\langle\kern-0.5ex\langle}
\def\rra{\rangle\kern-0.5ex\rangle}
\def\lsqb{[\kern-0.4ex[}
\def\rsqb{]\kern-0.4ex]}
\def\nn{\nonumber \\}

\def\ed{\end{document}}
\newcounter{mycnt}
\def\themycnt{\thesection.\arabic{mycnt}}
\def\mybenv#1{\refstepcounter{mycnt}%
       \vskip 3pt\noindent{\bf #1~~\themycnt}:~}
\def\myeenv{\mbox{~}\hfill\rule{1ex}{1ex}\vskip 3pt}
\def\qed{\mbox{~}\hfill$\Box$}
\renewcommand{\theequation}{\arabic{section}-\arabic{equation}}
\begin{document}
\title[Hopf algebras, distributive pairings and hash products]{%
       Hopf Algebras, Distributive (Laplace) Pairings\\
       and Hash Products:\\
       {\small A unified approach to tensor product decompositions\\
               of group characters}}
\author{Bertfried Fauser}
\address{%
School of Computer Science\\
The University of Birmingham\\
Edgbaston-Birmingham, W. Midlands, B15 2TT, England}
%\curraddr{}
\email{b.fauser@cs.bham.ac.uk Bertfried.Fauser@rhul.ac.uk}
%\thanks{}
\author{Peter D. Jarvis}
\address{%
School of Mathematics and Physics, University of Tasmania,
Private Bag 37, GPO, Hobart Tas 7001, Australia}
\email{Peter.Jarvis@utas.edu.au}
\author{Ronald C. King}
\address{%
School of Mathematics, University of Southampton,
Southampton SO17 1BJ, England}
\email{R.C.King@soton.ac.uk}

%%16W30 Coalgebras, bialgebras, Hopf algebras [See also 16S40, 57T05];
%%      rings, modules, etc. on which these act
%%05E05 Symmetric functions
%%33D52 Basic orthogonal polynomials and functions associated with
%%      root systems (Macdonald polynomials, etc.)
%%11E57 Classical groups [See also 14Lxx, 20Gxx]
%%43A40 Character groups and dual objects


\date{January 30, 2013}

\begin{abstract}
We show for bicommutative graded connected Hopf algebras that a certain
distributive (Laplace) subgroup of the convolution monoid of 2-cochains
parameterizes certain well behaved Hopf algebra deformations.
% We define higher derived hash products parameterized by this subgroup,
% providing a general theory of such deformations. 
Using the Laplace group, or its Frobenius subgroup, we define higher
derived hash products, and develop a general theory to study their
main properties. Applying our results to the (universal)
bicommutative graded connected Hopf algebra of symmetric functions,
we show that classical tensor product and character decompositions,
such as those for the general linear group, mixed co- and contravariant
or rational characters, orthogonal and symplectic group characters, 
Thibon and reduced symmetric group characters, are special cases of
higher derived hash products. In the Appendix we discuss a relation to
formal group laws.
\end{abstract}
%--TITLE----------------------------------------------------------------
\maketitle
%--TOC--comment out for final version-----------------------------------
{\tiny \tableofcontents}
%--PAPER STARTS HERE----------------------------------------------------
\section{Introduction and main results}

The main problem we want to address in this paper is that of showing
that there is a unique way to produce decompositions of products of
irreducible (or indecomposable) group representations or restricted
group characters for a variety of subgroups of the general linear
group $\GL(N)$. Let $V$ be a vector space of dimension $N$ together
with a $\GL(N)$ action and let $\RGL$ be the category of (finite
dimensional) $\GL(N)$-representations with a tensor product turning it
into a ring. Given a basis of irreducible representations
$V^{\lambda}$ with characters $s_{\lambda}$ one wants to compute the
coefficients $c^{\lambda}_{\mu,\nu}$ of the decomposition
\begin{align}
   V^{\mu}\ot V^{\nu}= \oplus_{\lambda\in\mathcal{P}}\oplus^{c^{\lambda}_{\mu,\nu}} V^{\lambda}
&&\textrm{or for characters}&&
   s_{\mu}\cdot s_{\nu}= \sum_{\lambda}c^{\lambda}_{\mu,\nu} s_{\lambda}
\end{align}
with sums over all partitons $\lambda \in \mathcal{P}$. 
Of course the coefficients in the $\GL$ case are the Littlewood Richardson
coefficients. We study \emph{restricted groups}, that is, subgroups $\H$ of
$\GL$ defined by polynomial identities. Classical examples of such groups
are orthogonal and symplectic groups but also discrete subgroups like
the symmetric group. The task at hand is to compute the decomposition
coefficients for irreducible (or indecomposable) representations and
characters of these subgroups. This is usually done in a large $N$ limit
(inductive limit) to avoid syzygies, the so-called \emph{modification rules}.
The corresponding formal characters are called universal characters, and
we restrict ourselves here to this case. Modification rules, at least
for classical groups, have been worked out on a case-by-case basis.
For characters of restricted groups we have demonstrated
in~\cite{fauser:jarvis:king:wybourne:2005a} that finding the decomposition
coefficients can be achieved by Hopf algebra deformations using a
twisted multiplication (and comultiplication). This process used an
isomorphism on the underlying module which was constructed using
plethystic Schur function series. Here we take another path generalizing,
in the commutative setting, a process developed by Rota and
Stein~\cite{rota:stein:1994a,rota:stein:1994b}. This process directly
deforms the multiplication using certain 2-cocycles of distributive
algebra valued pairings called Laplace pairings.

In this paper we develop a general theory of Hopf algebra deformations
aimed at the tensor product decomposition of universal restricted
group characters. For that we first provide, in 
Section~\ref{sec-HashProd}, a general theory of deformations in
bicommutative graded connected Hopf algebras $\HB$. To that end we
construct in Section~\ref{subsec:MonoidOfPairings} certain subgroups of
the monoid of 2-cocycles which will \emph{parameterize} the twists inducing
the deformations. We examine two special such subgroups. One is the subgroup
of distributive or Laplace pairings. Such pairings enjoy two straightening
or expansion laws, which can be seen as distributive laws. The second,
studied in Section~\ref{subsec:FobeniusLaplace}, is a subgroup of the
Laplace group consisting of pairings endowed with the further property
of being Frobenius, that is carrying a commutative Frobenius algebra
strcuture. We will show that deformations based on Frobenius Laplace
pairings induce Hopf algebra morphisms, while the general Laplace
pairings need a further deformation of the comultiplication to
remain Hopf. Our general development will shed some light on the
deformation process, and identify a certain condition (e), used by Rota
and Stein, as being equivalent to the Frobenius property. The deformation
process does not deform the comultiplication, and in general produces
neither Hopf nor bialgebras (but only $\HB$-comodule algebras). The
Frobenius property characterizes the cases where the deformation remains
Hopf. In cases for which there exists a module isomorphism between the
characters, we give in Section~\ref{sec-HashProdComul} a deformation of
the comultiplication necessary to retain the Hopf algebra property. 

We assume in the following that all Hopf algebras are both 
biassociative and biunital and do not mention this again explicitly. 
Our theoretical framework then starts with an \emph{ambient Hopf algebra}
assumed to be a graded connected bicommutative Hopf algebra. This Hopf
algebra is used to produce convolution monoids of $k$-cochains. Then
we define in Section~\ref{sec-HashProd} \emph{higher derived hash products}
by this process. The higher multiplications are obtained by multiple
convolutions with Laplace pairings. A derived multiplication will be
deformed by a derived pairing, that is a pairing which is composed with
a 1-cocycle. We need this generality for developing all the examples in
Section~\ref{sec-Applications} dealing with applications. In the
theory Sections~\ref{sec-ComHopf} and \ref{sec-HashProd} we show
among other results:
\begin{itemize}
\item Laplace pairings form a Laplace subgroup of the convolution monoid
      of 2-cochains.
\item Frobenius Laplace pairings form a subgroup of the Laplace group.
\item Deformations by Frobenius Laplace pairings induce Hopf algebra
      morphisms, while non Frobenius pairings do not. We give in
      Section~\ref{sec-HashProdComul} a deformation theory of the
      comultiplication which still produces Hopf algebra morphisms. 
\item Higher derived pairings are obtained by iterating the deformation
      process.            
\end{itemize}
The theory can in principle be extended to $n$-ary multiplications using
higher convolution monoids.

In Section~\ref{sec-Applications} we deal with group characters. For
that enterprise we specialize the ambient Hopf algebra to that of
symmetric functions $\Sym$ (or $\Sym\otimes \Sym$ for rational $\GL(N)$
characters). This is the universal (positive self 
adjoint~\cite{zelevinsky:1981b}) graded connected bicommutative Hopf
algebra. By the Cartier--Milnor--Moore theorem this Hopf algebra is
generated by polynomial generators. We show that our deformation
process is capable of producing the product decompositions for a variety
of classical groups.
\begin{itemize}
\item The branching $\GL(N)\downarrow\GL(N-1)$ is described by the
      trivial deformation. The Hopf algebra for this case is related
      to an additive formal group law.
\item The hash product describes the branching $\GL(N+M+NM)\downarrow
      \GL(N)\times \GL(M)$. We call the associated characters Thibon
      characters. These characters are related to stable permutation
      characters and Young polynomials and embody a multiplicative
      formal group law.
\item Changing the ambient Hopf algebra to $\Sym\ot\Sym$ we show that
      a derived hash product governs the product decomposition of
      rational $\GL(N)$ characters of mixed co- and contravariant
      $\GL(N)$-representations.
\item A similar derived hash product, again on $\Sym$, produces the
      Newell-Littlewood product formulae for orthogonal and symplectic
      characters.
\item Deforming the inner multiplication, induced by the symmetric
      group via the Frobenius characteristic map in $\Sym$, we show
      that a higher derived hash product produces the Murnaghan-Littlewood
      formula for reduced symmetric group characters obtained from a
      branching $\GL(N)$ $\downarrow \GL(N-1) \downarrow O(N-1) \downarrow S_n$.
      We investigate how Bernstein vertex operators are involved in
      the deformation process.
\end{itemize}
These character formulae are not new, and some of them have even been
derived using $\lambda$-ring or Hopf algebra techniques (references follow
in the text, but we want to mention
explicitly~\cite{scharf:thibon:wybourne:1993a,scharf:thibon:1994a}).
What we want to emphasize here is our unified approach which produces
\emph{all} of these results from a \emph{single source}.

To accomplish our goals we use graphical calculus throughout the paper,
especially in the development of the theory in Section~\ref{sec-ComHopf}.
The benefit of graphical calculus in manipulations is that it shows the
general basis free Hopf algebraic core of the process, and frees it from
its underlying combinatorial complexity. However, for concrete, say
machine, computations highly optimized combinatorics is indeed needed.
As the underlying modules come with a distinguished basis of irreducibles
(indecomposables) the basis dependence is a necessary part of the
interpretation of the result as character decompositions.

%-----------------------------------------------------------------------
\section*{Acknowledgement}
This work is the result of a collaboration over several years following
on from our paper \cite{fauser:jarvis:king:wybourne:2005a} with the late
Brian G Wybourne. BF gratefully acknowledges the Alexander von Humboldt
Stiftung for \emph{sur place} travel grants to visit the School of
Mathematics and Physics, University of Tasmania, and the University of
Tasmania for an honorary Research Associate appointment. Likewise PDJ
acknowledges longstanding support from the Alexander von Humboldt Foundation,
in particular for visits to the Max Planck Institute for Mathematics in
the Sciences, Leipzig. PDJ also acknowledges the Australian American
Fulbright Foundation for the award of a senior Fulbright scholarship.
RCK acknowledges support for part of this research from the Leverhulme
Foundation in the form of an Emeritus Fellowship. We thank the DAAD for
financial support allowing a visit of RCK and PDJ to Erlangen, and we
are grateful for financial support from the Emmy-Noether Zentrum f\"ur
Algebra, at the University of Erlangen. The authors thank these
institutions for hospitality during our collaborative visits, especially
Friedrich Knop (Erlangen) for discussions about the Murnaghan Littlewood
product formula for stable characters. The main body of this work was done
under a `Research in Pairs' grant from the Mathematisches Forschungsinstitut
Oberwolfach, and the authors also wish to record their appreciation of
this award and the hospitality extended to them in Oberwolfach. Finally,
this work could not have been completed without the generous support
of the Quantum Computing Group, Department of Computer Science,
University of Oxford, for hosting a research visit.
%------------------------------------------------------------ SECTION --
\section{Commutative Hopf algebras and distributive pairings}\label{sec-ComHopf}
%----------------------------------------------------------SUBSECTION --
\subsection{Graded commutative Hopf algebras}

We work with an abstract graded or filtered commutative cocommutative
Hopf algebra $\HB$. Biassociativity and biunitality are always assumed.
As a module $\HB=\oplus_n \H^n$ is $\mathbb{Z}_{\geq0}$-graded
with non negative degrees. A morphism $\sff$ of graded modules decomposes into
a family of maps $\sff^i : \H^i\rightarrow \H^i$, while for example the
multiplication map $\mul : \H^i \ot \H^j \rightarrow \H^{i+j}$ respects
grades additively. In the filtered case, we get for a pairing
$\sfa : \H^i\ot \H^j \rightarrow \oplus_{r=0}^{i+j} \H^r$. A special
case, used in applications below, is the inner multiplication acting
on a single degree only $* : \H^i\ot\H^i\rightarrow\H^i$. A detailed
description of graded Hopf algebras can be found in~\cite{milnor:moore:1965a}.
We make use of Heyneman-Sweedler index
notation~\cite{heyneman:sweedler:1969a,heyneman:sweedler:1970a},
see~\eqref{def-Sweedler-index}, and we distinguish Sweedler indices for
different comultiplications by using different brackets.

\mybenv{Definition}\label{def-gcHopfAlgebra}
Let $\HB=\H^0+\H^+=\oplus_{k\geq0}\H^k$ be a
$\mathbb{Z}_{\geq0}$-graded module over a commutative ring $\Bbbk$. Then
a connected graded bicommutative \emph{Hopf algebra} $\HB$ is given by
the following data:
\begin{align}
\mul &: \HB \ot \HB \rightarrow \HB :: (x,y) \mapsto xy \\
\eta &: \Bbbk \rightarrow \HB :: 1 \mapsto \eta(1)=1_{\HB}\\
\label{def-Sweedler-index}
\comul &: \HB \rightarrow \HB \ot \HB :: x \mapsto
    \comul(x) = x_{(1)} \ot x_{(2)} = {\sum}_i x_{i,1}\ot x_{i,2}\\
\epsilon &: \HB \rightarrow \Bbbk\cong\H^0 :: x \mapsto \epsilon(x){}
   = \epsilon(x^0+x^+) = x^{0} \\
\antip &: \HB \rightarrow \H^{\bullet,op} :: x \mapsto \antip(x)
\end{align}
such that 
\begin{align*}
 \mul^{(3)} &:= \mul\circ(\mul\ot\Id) = \mul\circ(\Id\ot\mul)	
&&
\textrm{associativity~\eqref{grph-assoc-com-mul}}
&&&
\\
 \mul\circ(\Id\ot\eta) &= \Id = \mul(\eta\ot\Id){}
&&
\textrm{unital~\eqref{grph-unit-counit}}
&&&
\\
 \comul^{(3)} &:= (\comul\ot\Id)\circ\comul = (\Id\ot\comul)\circ\comul	
&&
\textrm{coassociativity~\eqref{grph-assoc-com-comul}}
&&&
\\
 (\Id\ot\epsilon)\circ\comul &= \Id = (\epsilon\ot\Id)\circ\comul
&&
\textrm{counital~\eqref{grph-unit-counit}}
&&&
\\
 \mul &= \mul\circ \sw,\hskip2ex
\comul = \sw\circ\comul
&&
\textrm{bicommutativity~(\ref{grph-assoc-com-mul},\ref{grph-assoc-com-comul})}
&&&
\\
\comul(xy) &= \comul(x)\comul(y)
&&
\textrm{algebra homomorphism~\eqref{grph-bialgebra-and-antipode}}
&&&
\\
\ker\epsilon &= \H^+\hskip1ex \textrm{and}\hskip1ex \H^0\cong\Bbbk
&&
\textrm{connectedness}
&&&
\\
\antip(x_{(1)}) x_{(2)}
 &= \eta\circ\epsilon(x) 
  = \unit(x){}
&&
\textrm{antipode~\eqref{grph-bialgebra-and-antipode}}
&&&
\end{align*}
are satisfied (equation numbers point to their graphical representations).
\myeenv
In a graded connected bialgebra the antipode $\antip$ exists
automatically due to the inversion formula to be discussed below. We
will use graphical calculus, and use the example of a Hopf algebra to
introduce this notion. We read diagrams downwards along the
pessimistic arrow of time. Multiplication and comultiplication
are depicted by undecorated nodes, while the switch map
$\sw : x\ot y \rightarrow y\ot x$ is represented as a
crossing with no under or over information.
\begin{align}\label{grph-assoc-com-mul}
    \begin{pic}
      \node (i1) at (0,0.5) {};
      \node (i2) at (0.5,0.5) {};
      \node (i3) at (1,0.5) {};
      \node (m1) at (0.25,0.25) {};
      \node (m2) at (0.625,-0.25) {};
      \node (o) at (0.625,-0.5) {};
      \draw[thick,out= 270,in=180] (i1.center) to (m1.center);
      \draw[thick,out= 270,in=  0] (i2.center) to (m1.center);
      \draw[thick,out= 270,in=180] (m1.center) to (m2.center);
      \draw[thick,out= 270,in=  0] (i3.center) to (m2.center);
      \draw[thick] (m2.center) to (o.center);
   \end{pic} 
&\cong
    \begin{pic}
      \node (i1) at (0,0.5) {};
      \node (i2) at (0.5,0.5) {};
      \node (i3) at (1,0.5) {};
      \node (m1) at (0.75,0.25) {};
      \node (m2) at (0.375,-0.25) {};
      \node (o) at (0.375,-0.5) {};
      \draw[thick,out= 270,in=180] (i1.center) to (m2.center);
      \draw[thick,out= 270,in=180] (i2.center) to (m1.center);
      \draw[thick,out= 270,in=  0] (m1.center) to (m2.center);
      \draw[thick,out= 270,in=  0] (i3.center) to (m1.center);
      \draw[thick] (m2.center) to (o.center);
   \end{pic} 
;&&&
   \begin{pic}
      \node (i1) at (0,0.5) {};
      \node (i2) at (0.5,0.5) {};
      \node (d1) at (0,0.125) {};
      \node (d2) at (0.5,0.125) {};
      \node (m1) at (0.25,-0.125) {};
      \node (o) at (0.25,-0.5) {};
      \draw[thick,out=270,in=90] (i1.center) to (d2.center);
      \draw[thick,out=270,in=90] (i2.center) to (d1.center);
      \draw[thick,out=270,in=180] (d1.center) to (m1.center);
      \draw[thick,out=270,in=  0] (d2.center) to (m1.center);
      \draw[thick,out=270,in=90] (m1.center) to (o.center);
   \end{pic}
&\cong
   \begin{pic}
      \node (i1) at (0,0.5) {};
      \node (i2) at (0.5,0.5) {};
      \node (m1) at (0.25,-0.125) {};
      \node (o) at (0.25,-0.5) {};
      \draw[thick,out=270,in=180] (i1.center) to (m1.center);
      \draw[thick,out=270,in=  0] (i2.center) to (m1.center);
      \draw[thick,out=270,in=90] (m1.center) to (o.center);
   \end{pic}
;&&&
\parbox[c]{0.4\textwidth}{associativity and commutativity of the
           multiplication $\mul$.}
\end{align}
%%%%%%%%%%%%%%% comultiplication %%%%%%%%%%%%%%%%%%%%%%%%%%%%%%%%%%%%%%%
\begin{align}\label{grph-assoc-com-comul}
   \begin{pic}
      \node (i1) at (0,-0.5) {};
      \node (i2) at (0.5,-0.5) {};
      \node (i3) at (1,-0.5) {};
      \node (m1) at (0.25,-0.25) {};
      \node (m2) at (0.625,0.25) {};
      \node (o) at (0.625,0.5) {};
      \draw[thick,out= 90,in=180] (i1.center) to (m1.center);
      \draw[thick,out= 90,in=  0] (i2.center) to (m1.center);
      \draw[thick,out= 90,in=180] (m1.center) to (m2.center);
      \draw[thick,out= 90,in=  0] (i3.center) to (m2.center);
      \draw[thick] (m2.center) to (o.center);
   \end{pic}
&\cong
   \begin{pic}
      \node (i1) at (0,-0.5) {};
      \node (i2) at (0.5,-0.5) {};
      \node (i3) at (1,-0.5) {};
      \node (m1) at (0.75,-0.25) {};
      \node (m2) at (0.375,0.25) {};
      \node (o) at (0.375,0.5) {};
      \draw[thick,out= 90,in=180] (i1.center) to (m2.center);
      \draw[thick,out= 90,in=180] (i2.center) to (m1.center);
      \draw[thick,out= 90,in=  0] (m1.center) to (m2.center);
      \draw[thick,out= 90,in=  0] (i3.center) to (m1.center);
      \draw[thick] (m2.center) to (o.center);
   \end{pic}
;&&&
   \begin{pic}
      \node (i1) at (0,-0.5) {};
      \node (i2) at (0.5,-0.5) {};
      \node (d1) at (0,-0.125) {};
      \node (d2) at (0.5,-0.125) {};
      \node (m1) at (0.25,0.125) {};
      \node (o) at (0.25,0.5) {};
      \draw[thick,out=90,in=270] (i1.center) to (d2.center);
      \draw[thick,out=90,in=270] (i2.center) to (d1.center);
      \draw[thick,out=90,in=180] (d1.center) to (m1.center);
      \draw[thick,out=90,in=  0] (d2.center) to (m1.center);
      \draw[thick,out=90,in=90] (m1.center) to (o.center);
\end{pic}
&\cong
   \begin{pic}
      \node (i1) at (0,-0.5) {};
      \node (i2) at (0.5,-0.5) {};
      \node (m1) at (0.25,0.125) {};
      \node (o) at (0.25,0.5) {};
      \draw[thick,out=270,in=180] (i1.center) to (m1.center);
      \draw[thick,out=270,in=  0] (i2.center) to (m1.center);
      \draw[thick,out=270,in=90] (m1.center) to (o.center);
   \end{pic}
;&&&
\parbox[c]{0.4\textwidth}{coassociativity and cocommutativity of the
          comultiplication $\comul$.}
\end{align}
We refer to commutativity and cocommutativity collectively as bicommutativity,
and to associativity and coassociativity as biassociativity. Moreover, we
use associativity to define graphically iterated co/multiplications like
$\mul^{(3)}=\mul\ot(\Id\ot \mul)$
as nodes with many inputs and one output, similarly for or
$\comul^{(3)}=(\comul\ot \Id)\circ\comul$ with a reflected diagram, one input
many outputs. Tangle diagrams are referred to as $n$-$m$-tangles having
$n$ inputs and $m$ outputs. Equality of tangles up to allowed homotopies,
which do not produce crossings neither create nor delete extrema, is
denoted by $\cong$.
%%%%%unit and counit %%%%%%%%%%%%%%%%%%%%%%%%%%%%%%%%%%%%%%%%%%%%%%%%%%%
\begin{align}\label{grph-unit-counit}
%%%%%%%%%%%%%%%%%%%%%%%%%%%%%%%%%%%%%%%%%%%%%%%%%%%%%%%%%%%%%%%%%%%%%%%%
%% unit
   \begin{pic}
      \node[circle,inner sep=2pt,draw,thick] (i1) at (0,0.2) {};
      \node[xshift=-3pt] at (i1.west) {$\eta $};
      \node (i2) at (0.5,0.5) {};
      \node (m1) at (0.25,-0.125) {};
      \node (o) at (0.25,-0.5) {};
      \draw[thick,out=270,in=180] (i1) to (m1.center);
      \draw[thick,out=270,in=  0] (i2.center) to (m1.center);
      \draw[thick,out=270,in=90] (m1.center) to (o.center);
   \end{pic}
&\cong
   \begin{pic}
      \node (i) at (0,0.5) {};
      \node (o) at (0,-0.5) {};
      \draw[thick] (i.center) to (o.center);
   \end{pic}
\cong
   \begin{pic}
      \node[circle,inner sep=2pt,draw,thick] (i1) at (0.5,0.2) {};
      \node[xshift=3pt] at (i1.east) {$\eta $};
      \node (i2) at (0,0.5) {};
      \node (m1) at (0.25,-0.125) {};
      \node (o) at (0.25,-0.5) {};
      \draw[thick,out=270,in=  0] (i1) to (m1.center);
      \draw[thick,out=270,in=180] (i2.center) to (m1.center);
      \draw[thick,out=270,in=90] (m1.center) to (o.center);
   \end{pic}
;&&&
%%%%%%%%%%%%%%%%%%%%%%%%%%%%%%%%%%%%%%%%%%%%%%%%%%%%%%%%%%%%%%%%%%%%%%%%
%% counit
   \begin{pic}
      \node[circle,inner sep=2pt,draw,thick] (i1) at (0,-0.2) {};
      \node[xshift=-3pt] at (i1.west) {$\epsilon $};
      \node (i2) at (0.5,-0.5) {};
      \node (m1) at (0.25,0.125) {};
      \node (o) at (0.25,0.5) {};
      \draw[thick,out=90,in=180] (i1) to (m1.center);
      \draw[thick,out=270,in=  0] (i2.center) to (m1.center);
      \draw[thick,out=270,in=90] (m1.center) to (o.center);
   \end{pic}
&\cong
   \begin{pic}
      \node (i) at (0,0.5) {};
      \node (o) at (0,-0.5) {};
      \draw[thick] (i.center) to (o.center);
   \end{pic}
\cong
   \begin{pic}
      \node[circle,inner sep=2pt,draw,thick] (i1) at (0.5,-0.2) {};
      \node[xshift=3pt] at (i1.east) {$\epsilon $};
      \node (i2) at (0,-0.5) {};
      \node (m1) at (0.25,0.125) {};
      \node (o) at (0.25,0.5) {};
      \draw[thick,out=90,in=0] (i1) to (m1.center);
      \draw[thick,out=270,in=180] (i2.center) to (m1.center);
      \draw[thick,out=270,in=90] (m1.center) to (o.center);
   \end{pic}
;&&&
\parbox[c]{0.3\textwidth}{unit $\eta$ and counit $\epsilon$ for
          multiplication $\mul$ and comultiplication $\comul$.}
\end{align}
To form a bialgebra, multiplication and comultiplication have to
satisfy a compatibility law. This demands that the comultiplication is
an algebra homomorphism and \emph{vice versa} the multiplication is
a coalgebra homomorphism, shown by the horizontal symmetry of the diagram.
The antipode is defined as the convolutive inverse of the identity morphism
$\Id : \HB \rightarrow \HB$, not to be confused with the linear inverse.
%% BIALGEBRA COMPATIBILITY %%%%%%%%%%%%%%%%%%%%%%%%%%%%%%%%%%%%%%%%%%%%%
\begin{align}\label{grph-bialgebra-and-antipode}
\begin{pic}
      \node (i1) at (0,0.5) {};
      \node (i2) at (0.5,0.5) {};
      \node (u) at (0.25,0.2) {};
      \node (d) at (0.25,-0.2) {};
      \node (o1) at (0,-0.5) {};
      \node (o2) at (0.5,-0.5) {};
      \draw[thick,out=270,in=180] (i1.center) to (u.center);
      \draw[thick,out=270,in=  0] (i2.center) to (u.center);
      \draw[thick] (u.center) to (d.center);
      \draw[thick,out=180,in=90] (d.center) to (o1.center);
      \draw[thick,out=  0,in=90] (d.center) to (o2.center);
   \end{pic}
&\cong
   \begin{pic}
      \node (i1) at (0.25,0.5) {};
      \node (i2) at (0.75,0.5) {};
      \node (u1) at (0.25,0.3) {};
      \node (u2) at (0.75,0.3) {};
      \node (m1) at (0,0) {};
      \node (m4) at (1,0) {};
      \node (d1) at (0.25,-0.3) {};
      \node (d2) at (0.75,-0.3) {};
      \node (o1) at (0.25,-0.5) {};
      \node (o2) at (0.75,-0.5) {};
      \draw[thick] (i1.center) to (u1.center);
      \draw[thick] (i2.center) to (u2.center);
      \draw[thick,out=180,in=90] (u1.center) to (m1.center);
      \draw[thick,out=  0,in=180] (u1.center) to (d2.center);
      \draw[thick,out=180,in=  0] (u2.center) to (d1.center);
      \draw[thick,out=  0,in=90] (u2.center) to (m4.center);
      \draw[thick,out=270,in=180] (m1.center) to (d1.center);
      \draw[thick,out=270,in=  0] (m4.center) to (d2.center);
      \draw[thick] (d1.center) to (o1.center);
      \draw[thick] (d2.center) to (o2.center);
   \end{pic}
;&&&
%% ANTIPODE %%%%%%%%%%%%%%%%%%%%%%%%%%%%%%%%%%%%%%%%%%%%%%%%%%%%%%%%%%%%
   \begin{pic}
      \node (i) at (0.25,0.5) {};
      \node (u) at (0.25,0.35) {};
      \node[circle,inner sep=1pt,draw,thick,fill=white] (a) at (0,0) {$\antip$};
      \node (m) at (0.5,0) {};
      \node (d) at (0.25,-0.35) {};
      \node (o) at (0.25,-0.5) {};
\begin{pgfonlayer}{background}      
      \draw[thick] (i.center) to (u.center);
      \draw[thick,out=180,in=90] (u.center) to (a.center);
      \draw[thick,out=  0,in=90] (u.center) to (m.center);
      \draw[thick,out=270,in=180] (a.center) to (d.center);
      \draw[thick,out=270,in=  0] (m.center) to (d.center);
      \draw[thick] (d.center) to (o.center);
\end{pgfonlayer}
   \end{pic}
&\cong
   \begin{pic}
      \node (i) at (0,0.5) {};
      \node (o) at (0,-0.5) {};
      \node[circle,inner sep=2pt,draw,thick] (eps) at (0,0.15) {};
      \node[xshift=3pt] at (eps.east) {$\epsilon$};
      \node[circle,inner sep=2pt,draw,thick] (eta) at (0,-0.15) {};
      \node[xshift=3pt] at (eta.east) {$\eta$};
      \draw[thick] (i.center) to (eps);
      \draw[thick] (eta) to (o.center);
   \end{pic}
;&&&
\parbox[c]{0.4\textwidth}{bialgebra compatibility law and the antipode
          $\antip$ as convolutive inverse of identity $\Id$. Identities
          are not drawn.}
\end{align}
The antipode is a left and right inverse as we are working in a
bicommutative setting. This completes the diagrams for the Hopf algebra
$\HB$

A further important notion available for graded connected comodules
(Hopf algebras) is that of a cut comultiplication
$\comul^\prime : \HB \rightarrow \H^+\ot \H^+$. It is defined as that
part of a comultiplication which splits its input in a nontrivial way.
On $\H^0$ set $\Delta^{\prime}=0$ and on $\H^+$ define
$\comul^\prime := \comul - (\eta\ot \Id) - (\Id \ot \eta)$. Graphically
we denote the cut by double lines, only elements in $\ker\epsilon=\H^+$
can pass through this gate.
\begin{align}\label{grph-cutCoproduct}
   \begin{pic}
      \node (i) at (0,0.75) {};
      \node[iso2cell,minimum width=5pt] (a) at (0,0) {};
      \node (o) at (0,-0.75) {};
      \draw[thick] (i.center) to (a);
      \draw[thick] (a) to (o.center);
   \end{pic}
:\cong   
   \begin{pic}
	   \node (i1) at (0,0.75) {};
	   \node (o1) at (0,-0.75) {};
	   \draw[thick] (i1.center) to (o1.center);   
   \end{pic}
-
   \begin{pic}
      \node (i) at (0,0.75) {};
      \node (o) at (0,-0.75) {};
      \node[circle,inner sep=2pt,draw,thick] (eps) at (0,0.2) {};
      \node[xshift=3pt] at (eps.east) {$\epsilon$};
      \node[circle,inner sep=2pt,draw,thick] (eta) at (0,-0.2) {};
      \node[xshift=3pt] at (eta.east) {$\eta$};
      \draw[thick] (i.center) to (eps);
      \draw[thick] (eta) to (o.center);
   \end{pic}
;&&&
   \begin{pic}
      \node (i) at (0.25,0.75) {};
      \node (m) at (0.25,0.25) {};
      \node[iso2cell,minimum width=5pt] (c0) at (0.25,0.5) {};
      \node[iso2cell,minimum width=5pt] (c1) at (0,0) {};
      \node[iso2cell,minimum width=5pt] (c2) at (0.5,0) {};
      \node (o1) at (0,-0.75) {};
      \node (o2) at (0.5,-0.75) {};
      \draw[thick] (i.center) to (c0.north);
      \draw[thick] (c0.south) to (m.center);
      \draw[thick,out=180,in=90] (m.center) to (c1.north);
      \draw[thick,out=  0,in=90] (m.center) to (c2.north);
      \draw[thick] (c1.south) to (o1.center);
      \draw[thick] (c2.south) to (o2.center);
   \end{pic}
&\hskip-2ex:\cong
   \begin{pic}
      \node (i) at (0.25,0.75) {};
      \node[iso2cell,minimum width=5pt] (c0) at (0.25,0.5) {};
      \node (m) at (0.25,0.125) {};
      \node (o1) at (0,-0.75) {};
      \node (o2) at (0.5,-0.75) {};
      \draw[thick] (i.center) to (c0.north);
      \draw[thick] (c0.south) to (m.center);
      \draw[thick,out=180,in=90] (m.center) to (o1.center);
      \draw[thick,out=  0,in=90] (m.center) to (o2.center);
   \end{pic}
-
   \begin{pic}
      \node (i) at (0,0.75) {};
      \node[iso2cell,minimum width=5pt] (c0) at (0,0.5) {};
      \node[circle,inner sep=2pt,draw,thick] (eta) at (0.25,0.125) {};
      \node[xshift=3pt] at (eta.east) {$\eta$};
      \node (o1) at (0,-0.75) {};
      \node (o2) at (0.25,-0.75) {};
      \draw[thick] (i.center) to (c0.north);
      \draw[thick] (c0.south) to (o1.center);
      \draw[thick] (eta) to (o2.center);
   \end{pic}
-
\hskip-1ex
   \begin{pic}
      \node (i) at (0.25,0.75) {};
      \node[iso2cell,minimum width=5pt] (c0) at (0.25,0.5) {};
      \node[circle,inner sep=2pt,draw,thick] (eta) at (0,0.125) {};
      \node[xshift=-3pt] at (eta.west) {$\eta$};
      \node (o1) at (0,-0.75) {};
      \node (o2) at (0.25,-0.75) {};
      \draw[thick] (i.center) to (c0.north);
      \draw[thick] (c0.south) to (o2.center);
      \draw[thick] (eta) to (o1.center);
   \end{pic}
;&&&
\parbox[c]{0.3\textwidth}{the cut comultiplication $\comul^\prime$.}
\end{align}
Connectedness (Definition~\ref{def-gcHopfAlgebra}) is needed to ensure
that $\Bbbk\ot\eta(1)\cong \Bbbk\cong \H^0$. This ensures that the two
terms $\eta\ot \Id$ and $\Id\ot \eta$ eliminate \emph{all} trivial parts
of the coproduct and allow recursive expansions.

The cut coproduct can be utilized to produce a recursive formula for
the antipode, using an inclusion-exclusion principle, which may however
result in a good deal of overcounting and a lot of cancellation.
\begin{align}\label{grph-recursive-antipode}
   \begin{pic}
      \node (i) at (0,0.75) {};
      \node[iso2cell,minimum width=5pt] (c0) at (0,0.5) {};      
      \node[circle,inner sep=1pt,draw,thick] (a) at (0,0) {$\antip$};
      \node (o) at (0,-0.75) {};
      \draw[thick] (i.center) to (c0.north);
      \draw[thick] (c0.south) to (a);
      \draw[thick] (a) to (o.center);
   \end{pic}
&\hskip1ex\cong
   %\begin{pic}
      %\node (i) at (0,0.5) {};
      %\node (o) at (0,-0.5) {};
      %\node[circle,inner sep=2pt,draw,thick] (eps) at (0,0.15) {};
      %\node[xshift=3pt] at (eps.east) {$\epsilon$};
      %\node[circle,inner sep=2pt,draw,thick] (eta) at (0,-0.15) {};
      %\node[xshift=3pt] at (eta.east) {$\eta$};
      %\draw[thick] (i.center) to (eps);
      %\draw[thick] (eta) to (o.center);
   %\end{pic}
-
   \begin{pic}
      \node (i) at (0,0.75) {};
      \node[iso2cell,minimum width=5pt] (c0) at (0,0) {};      
      \node (o) at (0,-0.75) {};
      \draw[thick] (i.center) to (c0.north) (c0.south) to (o.center);
   \end{pic}
-
   \begin{pic}
      \node (i) at (0.25,1) {};
      \node[iso2cell,minimum width=5pt] (c0) at (0.25,0.75) {};
      \node (m1) at (0.25,0.5) {};
      \node[iso2cell,minimum width=5pt] (c1) at (0,0.25) {};
      \node[iso2cell,minimum width=5pt] (c2) at (0.5,0.25) {};
      \node[circle,inner sep=1pt,draw,thick] (a) at (0,-0.15) {$\antip$};
      \node (d) at (0.5,-0.35) {};
      \node (m2) at (0.25,-0.6) {};
      \node (o) at (0.25,-0.75) {};
      \draw[thick] (i.center) to (c0.north);
      \draw[thick] (c0.south) to (m1.center);
      \draw[thick,out=180,in=90] (m1.center) to (c1.north);
      \draw[thick,out=  0,in=90] (m1.center) to (c2.north);
      \draw[thick] (c1.south) to (a.north);
      \draw[thick] (c2.south) to (d.center);
      \draw[thick,out=270,in=180] (a.south) to (m2.center);
      \draw[thick,out=270,in=0] (d.center) to (m2.center);
      \draw[thick] (m2.center) to (o.center);
   \end{pic}
%%%%%%%%%%%%%%%%%%%%%%%%%%%%%%%%%%%%%%%%%%%%%%%%%%%%%%%%%%%%%%%%%%%%%%%%
%% antip recursion initial cond
&&&\textrm{with}\hskip2ex
   \begin{pic}
      \node[circle,inner sep=2pt,draw,thick] (eta) at (0,0.6) {};
      \node[xshift=3pt] at (eta.east) {$\eta $};
      \node[circle,inner sep=1pt,draw,thick] (a) at (0,0) {$\antip $};
      \node (o) at (0,-0.6) {};
      \draw[thick] (eta) to (a.north);
      \draw[thick] (a.south) to (o.center);
      \end{pic}
   &\cong
      \begin{pic}
      \node[circle,inner sep=2pt,draw,thick] (eta) at (0,0.6) {};
      \node[xshift=3pt] at (eta.east) {$\eta $};
      \node (o) at (0,-0.6) {};
      \draw[thick] (eta) to (o.center);
   \end{pic}
;&&&
\parbox[c]{0.35\textwidth}{recursive antipode formula, with initial
          condition $\antip(1)=1$.}
\end{align}
This formula
%\footnote{Some times called Connes-Kreimer antipode formula.}
is a special case of the inversion formula by
Milnor--Moore~\cite{milnor:moore:1965a} as we will see below.
%-----------------------------------------------------------------------
\subsection{Convolution algebras over a Hopf algebra}
%-----------------------------------------------------------------------
From now on we fix a graded connected bicommutative (biassociative)
Hopf algebra $\HB$ over a commutative ring $\Bbbk$. This will be our
ambient Hopf algebra used to define further structure.

One encounters several ways to compose maps in a symmetric monoidal
category.\footnote{We are interested here in the category of finite dimensional
$\Bbbk$-modules, and later in categories of group representations.} 
We have \emph{sequential composition} of maps
\begin{align}
   \begin{pic}
      \node (i) at (0,1) {};
      \node[rectangle,draw,thick] (f) at (0,0.4) {$\sff$};
      \node[rectangle,draw,thick] (g) at (0,-0.4) {$\sfg$};
      \node (o) at (0,-1) {};
      \draw[thick] (i.center) to (f.north);
      \draw[thick] (f.south) to (g.north);
      \draw[thick] (g.south) to (o.center);
   \end{pic}
&\cong
   \begin{pic}
      \node (i) at (0,1) {};
      \node[rectangle,draw,thick] (fg) at (0,0) {$\sfg\circ \sff$};
      \node (o) at (0,-1) {};
      \draw[thick] (i.center) to (fg.north);
      \draw[thick] (fg.south) to (o.center);
   \end{pic}
\end{align}
and \emph{parallel composition} of maps acting on different parts of
a tensor product space
\begin{align}
\begin{pic}
\node (i1) at (0,1) {};
\node (i2) at (0.75,1) {};
\node[rectangle,draw,thick] (f) at (0,0) {$\sff $};
\node[rectangle,draw,thick] (g) at (0.75,0) {$\sfg $};
\node (o1) at (0,-1) {};
\node (o2) at (0.75,-1) {};
\draw[thick] (i1.center) to (f.north);
\draw[thick] (i2.center) to (g.north);
\draw[thick] (f.south) to (o1.center);
\draw[thick] (g.south) to (o2.center);
\end{pic}
&\cong
   \begin{pic}
      \node (i1) at (-0.25,1) {};
      \node (i2) at (0.25,1) {};
      \node[rectangle,draw,thick] (fg) at (0,0) {$\sff\ot \sfg$};
      \node (o1) at (-0.25,-1) {};
      \node (o2) at (0.25,-1) {};
      \draw[thick] (i1.center) to (i1 |- fg.north);
      \draw[thick] (i1 |- fg.south) to (o1.center);
      \draw[thick] (i2.center) to (i2 |- fg.north);
      \draw[thick] (i2 |- fg.south) to (o2.center);
   \end{pic}
\end{align}
It is important to note here, that in the graded setting in use here,
if $\sff^i$ is of degree $i$ and $\sfg^j$ is of degree $j$, then
$\sff^i\ot \sfg^j$ acts on a space $\H^{i}\ot \H^{j}\subseteq\H^{i+j}$
of degree $i+j$. In this paper we do \emph{not} distinguish in the
graphics between these spaces, viewing them as subspaces of $\HB$. In
that sense our lines (or strings) for $\HB$ are \emph{cables} and graded
pieces like $\H^i$ would be \emph{wires} in a language used in spin
networks. This allows us to merge lines in the exchange law
below~\eqref{grph-exchangelaw}. However, using multiple lines usually
indicates we are working in $\HBI{k} = \HB\otimes\ldots\ot\HB$ for
some $k$.

Parallel and sequential or vertical composition enjoy an exchange law,
which is trivial if presented graphically
\begin{align}\label{grph-exchangelaw}
%% EXCHANGE LAW %%%%%%%%%%%%%%%%%%%%%%%%%%%%%%%%%%%%%%%%%%%%%%%%%%%%%%%%
   \begin{pic}
      \node (i1) at (-0.75,1) {};
      \node (i2) at (0.75,1) {};      
      \node[rectangle,draw,thick] (f) at (0,0) {$(\sfg\ot \sfv)\circ(\sff\ot \sfu)$};
      \node (o1) at (-0.75,-1) {};
      \node (o2) at (0.75,-1) {};
      \draw[thick] (i1.center) to (i1 |- f.north);
      \draw[thick] (o1|- f.south) to (o1.center);
      \draw[thick] (i2.center) to (i2 |- f.north);
      \draw[thick] (o2|- f.south) to (o2.center);
   \end{pic}
&\cong
   \begin{pic}
      \node (i1) at (-0.25,1) {};
      \node (i2) at (0.25,1) {};      
      \node[rectangle,draw,thick] (f) at (0,0.4) {$\sff\ot \sfu$};
      \node[rectangle,draw,thick] (g) at (0,-0.4) {$\sfg\ot \sfv$};
      \node (o1) at (-0.25,-1) {};
      \node (o2) at (0.25,-1) {};
      \draw[thick] (i1.center) to (i1 |- f.north);
      \draw[thick] (i1 |- f.south) to (i1 |- g.north);
      \draw[thick] (i1 |- g.south) to (o1.center);
      \draw[thick] (i2.center) to (i2 |- f.north);
      \draw[thick] (i2 |- f.south) to (i2 |- g.north);
      \draw[thick] (i2 |- g.south) to (o2.center);
   \end{pic}
\cong
   \begin{pic}
      \node (i1) at (0,1) {};
      \node (i2) at (1,1) {};
      \node[rectangle,draw,thick] (f) at (0,0.4) {$\sff$};
      \node[rectangle,draw,thick] (g) at (0,-0.4) {$\sfg$};
      \node[rectangle,draw,thick] (u) at (1,0.4) {$\sfu$};
      \node[rectangle,draw,thick] (v) at (1,-0.4) {$\sfv$};
      \node (o1) at (0,-1) {};
      \node (o2) at (1,-1) {};
      \draw[thick] (i1.center) to (f.north);
      \draw[thick] (f.south) to (g.north);
      \draw[thick] (g.south) to (o1.center);
      \draw[thick] (i2.center) to (u.north);
      \draw[thick] (u.south) to (v.north);
      \draw[thick] (v.south) to (o2.center);
   \end{pic}
\cong
   \begin{pic}
      \node (i1) at (0,1) {};
      \node (i2) at (1.2,1) {};
      \node[rectangle,draw,thick,inner sep=2pt] (f) at (0,0) {$\sfg\circ\sff$};
      \node[rectangle,draw,thick,minimum height=14pt] (g) at (1.2,0) {$\sfv\circ\sfu$};
      \node (o1) at (0,-1) {};
      \node (o2) at (1.2,-1) {};
      \draw[thick] (i1.center) to (f.north);
      \draw[thick] (f.south) to (o1.center);
      \draw[thick] (i2.center) to (g.north);
      \draw[thick] (g.south) to (o2.center);
   \end{pic}
\cong
   \begin{pic}
      \node (i1) at (-0.75,1) {};
      \node (i2) at (0.75,1) {};      
      \node[rectangle,draw,thick] (f) at (0,0) {$(\sfg\circ \sff)\ot(\sfv\circ \sfu)$};
      \node (o1) at (-0.75,-1) {};
      \node (o2) at (0.75,-1) {};      
      \draw[thick] (i1.center) to (i1 |- f.north);
      \draw[thick] (i1 |- f.south) to (o1.center);
      \draw[thick] (i2.center) to (i2 |- f.north);
      \draw[thick] (i2 |- f.south) to (o2.center);
   \end{pic}
\end{align}
Some of the graphical rearrangements which we use below will make free
use of such manipulations.

A further task is to define several convolution monoids, subgroups of
which will play a central role in parameterizing Hopf algebra
deformations. These convolution monoids will be defined on $n$-$1$-maps.

\mybenv{Definition}\label{def-convolution-monoid}
Let $\sff,\sfg \in \hom(\HB,\HB)$. We define the \emph{convolution monoid}
of $1$-$1$-maps with convolution multiplication
$\conv : \hom(\HB,\HB) \ot \hom(\HB,\HB) \rightarrow \hom(\HB,\HB) ::
(\sff,\sfg) \mapsto \sff\conv \sfg := \mul\circ(\sff\ot\sfg)\circ\comul$
and unit $\sfe=\eta\circ\epsilon$. Graphically $\conv$ is given as:
\begin{align}
   \begin{pic}
      \node (i) at (0.4,1) {};
      \node (u) at (0.4,0.75) {};
      \node[rectangle,draw,thick] (f) at (0,0) {$\sff$};
      \node[rectangle,draw,thick] (g) at (0.8,0) {$\sfg$};
      \node (d) at (0.4,-0.75) {};
      \node (o) at (0.4,-1) {};
      \draw[thick] (i.center) to (u.center);
      \draw[thick,out=180,in=90] (u.center) to (f.north);
      \draw[thick,out=  0,in=90] (u.center) to (g.north);
      \draw[thick,out=270,in=180] (f.south) to (d.center);
      \draw[thick,out=270,in=  0] (g.south) to (d.center);
      \draw[thick] (d.center) to (o.center);
   \end{pic}
&\cong
   \begin{pic}
      \node (i) at (0,1) {};
      \node[rectangle,draw,thick] (fg) at (0,0) {$\sff\conv \sfg$};
      \node (o) at (0,-1) {};
      \draw[thick] (i.center) to (fg.north);
      \draw[thick] (fg.south) to (o.center);
   \end{pic}
;&&&
   \begin{pic}
      \node (i) at (0,1) {};
      \node[rectangle,draw,thick] (fg) at (0,0) {$\sfe$};
      \node (o) at (0,-1) {};
      \draw[thick] (i.center) to (fg.north);
      \draw[thick] (fg.south) to (o.center);
   \end{pic}
&\cong
   \begin{pic}
      \node (i) at (0,1) {};
      \node (o) at (0,-1) {};
      \node[circle,inner sep=2pt,draw,thick] (eps) at (0,0.15) {};
      \node[xshift=3pt] at (eps.east) {$\epsilon$};
      \node[circle,inner sep=2pt,draw,thick] (eta) at (0,-0.15) {};
      \node[xshift=3pt] at (eta.east) {$\eta$};
      \draw[thick] (i.center) to (eps);
      \draw[thick] (eta) to (o.center);
   \end{pic}
&&&
\end{align}
The \emph{unit} fulfils, $\sfe\conv\sff=\sff$ by unitality of $\mul$
and counitality of $\comul$. The convolution $\conv$ is commutative
and associative if and only if $\mul,\comul$ are bicommutative and
biassociative.
\myeenv

As it is easily inferred which convolution is at hand we do not introduce
further notation such as $\conv^{n}$.
%-----------------------------------------------------------------------
We shall need the inverse of a map $\sff \in \hom(\HB,\HB)$, which was
given by Milnor and Moore~\cite{milnor:moore:1965a}.

We use the following terminology for maps $\sff\in\hom(\HB,\HB)$: A map
$\sff$ is \emph{normalized} if $\sff\circ\eta = \eta$, in which case $\sff(1)=1$.
For example $\epsilon$ is normalized as we find $\epsilon(\eta(1))=1$.
Normalization applies to duals too and $\sff$ is called \emph{conormalized}
if $\epsilon\circ\sff =\epsilon$. Composing maps in $\hom(\HBI{k},\HB)${}
with a linear form (such as $\epsilon$), shows that linear functionals
$\hom(\HB,\Bbbk)$ can also be equipped with a convolution product.
Indeed, if we let $\comul$ be the diagonal map $\comul(x)=x\ot x$ then
the convolution in $\hom(\HB,\Bbbk)$ given by
$\sff\conv\sfg(x)=\sff(x)\sfg(x)$ with the multiplication from $\Bbbk$ is 
point-wise multiplication of functionals in $\hom(\HB,\Bbbk)$.

\mybenv{Lemma}\label{lem-2convolution}
$\sff\in\hom(\HB,\HB)$ has a \emph{convolutive inverse} $\sfof$
if and only if $\sff(1)\not=0$. This can be computed recursively using the cut
coproduct $\comul'$ as (with normalization $f(1)=1$ and since
$\sff(1)\sfof(1)=1$ we have initial condition $\sfof(1)=1$ ):
\begin{align}
   \begin{pic}
      \node (i) at (0,1) {};
      \node[iso2cell,minimum width=5pt] (c0) at (0,0.75) {};
      \node[rectangle,draw,thick,minimum height=18pt] (a) at (0,0) {$\sfof$};
      \node (o) at (0,-1) {};
      \draw[thick] (i.center) to (c0.north) (c0.south) to (a);
      \draw[thick] (a) to (o.center);
   \end{pic}
=
   %\begin{pic}
      %\node (i) at (0,1) {};
      %\node (o) at (0,-1) {};
      %\node[circle,inner sep=2pt,draw,thick] (eps) at (0,0.15) {};
      %\node[xshift=3pt] at (eps.east) {$\epsilon$};
      %\node[circle,inner sep=2pt,draw,thick] (eta) at (0,-0.15) {};
      %\node[xshift=3pt] at (eta.east) {$\eta$};
      %\draw[thick] (i.center) to (eps);
      %\draw[thick] (eta) to (o.center);
   %\end{pic}
-
   \begin{pic}
      \node (i) at (0,1) {};
      \node[iso2cell,minimum width=5pt] (c0) at (0,0.75) {};
      \node[rectangle,draw,thick,minimum height=18pt] (a) at (0,0) {$\sff$};
      \node (o) at (0,-1) {};
      \draw[thick] (i.center) to (c0.north) (c0.south) to (a);
      \draw[thick] (a) to (o.center);
\end{pic}
-
   \begin{pic}
      \node (i) at (0.25,1.25) {};
      \node (m1) at (0.25,0.75) {};
      \node[iso2cell,minimum width=5pt] (c0) at (0.25,1) {};
      \node[iso2cell,minimum width=5pt] (c1) at (0,0.4) {};
      \node[iso2cell,minimum width=5pt] (c2) at (0.5,0.4) {};
      \node[rectangle,draw,thick,minimum height=18pt] (a1) at (0,-0.15) {$\sfof$};
      \node[rectangle,draw,thick,minimum height=18pt] (a2) at (0.5,-0.15) {$\sff$};
      \node (m2) at (0.25,-0.75) {};
      \node (o) at (0.25,-1) {};
      \draw[thick] (i.center) to (c0.north) (c0.south) to (m1.center);
      \draw[thick,out=180,in=90] (m1.center) to (c1.north);
      \draw[thick,out=  0,in=90] (m1.center) to (c2.north);
      \draw[thick] (c1.south) to (a1.north);
      \draw[thick] (c2.south) to (a2.north);
      \draw[thick,out=270,in=180] (a1.south) to (m2.center);
      \draw[thick,out=270,in=0] (a2.south) to (m2.center);
      \draw[thick] (m2.center) to (o.center);
   \end{pic}
;\quad\hskip1cm
   \begin{pic}
      \node at (0,1) {};
      \node[circle,inner sep=2pt,draw,thick] (eps) at (0,0.6) {};
      \node[xshift=3pt] at (eps.east) {$\eta$};
      \node[rectangle,draw,thick,minimum height=18pt] (a1) at (0,-0.15) {$\sfof$};
      \node (o) at (0,-1) {};
      \draw[thick] (eps.south) to (a1.north);
      \draw[thick] (a1.south) to (o.center);
   \end{pic}
\cong
\begin{pic}
      \node at (0,1) {};
      \node[circle,inner sep=2pt,draw,thick] (eps) at (0,0.6) {};
      \node[xshift=3pt] at (eps.east) {$\eta$};
%      \node at (0.4,0) {$\ot~\Bbbk$};
      \node (o) at (0,-1) {};
      \draw[thick] (eps.south) to (o.center);
   \end{pic}
\end{align}
\vskip-3ex
\myeenv
Let $\sff(1)$ be a unit in $\Bbbk$ we can define a normalized
$\sff^{\prime}$ for any such $1$-$1$-map $\sff\in\hom(\HB,\HB)$. In
the graphics above we have hence assumed that $\sff(1)=1$. 

The convolution of 1-1-maps $\sff \in \hom(\HB,\HB)$ with one input and
one output is a special case of more general $k$-$1$-maps with $k$ inputs
and one output. As these maps and their convolutions will be related to
Sweedler cohomology below, we refer to \emph{invertible} $k$-$1$-maps as
$k$-cochains, see Definition~\ref{def-cochain}, and we consider the
subgroup of such $k$-cochains in $\hom(\HBI{k},\H)$.

Let $\sff,\sfg \in \hom(\HBI{k},\HB)$ be $k$-1-maps. We extend
the above notion of convolution in the obvious manner. First we need a
coproduct on the space $\HBI{k}$, using the symmetric tensor product we
iterate $\comul_\HB$ to $\comul_{\HBI{2}}
= \Sigma_2\circ(\comul\ot\comul)$, with $\Sigma_2=\Id\ot\sw\ot\Id$
acting on $4=2k$ strings with $k=2$. For the general case on $k$ strings we obtain the
comultiplication on $\HBI{k}$ as $\comul_{\HBI{k}} =
\Sigma_k\circ\ot^{k}\comul$ with $\ot^k\comul=\comul\ot\ldots\ot\comul$
and the following definition of $\Sigma_{k}$. 
The comultiplication on $k$ strings produces $2k$ strings,
which we may enumerate as $(1,2,3,\ldots,2k)$. We define $\Sigma_{k}$
to be the (inverse shuffle) permutation on $2k$ strings reordering them
into $(1,3,5,\ldots,2k -1, 2,4,6,\ldots,2k)$. 
%The first and last
%string are fixed. We use the elementary transpositions
%$\sw_{i}=\sw^{(k)}_i = \Id^{i-1}\ot\sw\ot\Id^{2k-i-1}$ exchanging the strings
%$(i,i+1)$ for $i\in\{1,\ldots,2k-1\}$. We define the base case
%$\Sigma_{1}^{(k)} = \Id^{(k)}=\Id^{k}\otimes \Id^{k}$ acting on $2k$ strings
%as identity, and further
%$\Sigma_{l}=\Sigma_{l}^{(k)} = \circ_{i=k-l+2}^{k+l-2}\sw^{(k)}_{i} \circ \Sigma^{(k)}_{l-1}$
%with step-size $2$ for $i$ and iterating up to $l=k$. For example:
%$\Sigma_{2}^{(2)}= (\sw_{2})\circ\Id^{(2)}$,
%$\Sigma_{3}^{(3)}= (\sw_{2}\circ\sw_{4})\circ\Sigma_{2}^{(3)}
                 %=(\sw_{2}\circ\sw_{4})\circ(\sw_{3})\circ\Sigma_{1}^{(3)}
                 %= (\sw_{2}\circ\sw_{4})\circ(\sw_{3})\circ\Id^{(3)}$, 
%$\Sigma_{4}^{(4)}= (\sw_{2}\circ\sw_{4}\circ\sw_{6})\circ(\sw_{3}\circ\sw_{5})\circ(\sw_{4})\circ\Id^{(4)}$
%etc. We drop usually the upper index counting input strings. 
Extending Lemma~\ref{lem-2convolution} we get for general $k$-1-maps:

\mybenv{Lemma}
The spaces $\hom(\HBI{k},\HB)$, together with the convolution product 
$\conv$ based on $\mul,\comul_{\HBI{k}}$ with unit $\sfe=\sfe^k$,{}
form a \emph{$k$-convolution monoid}. Explicitly:
\begin{align}\label{grph:inverse}
  \sff\conv\sfg
     &=  \mul\circ(\sff\ot \sfg)\circ\comul_{\HBI{k}};
     \qquad
       \sfe
      = \eta\circ(\ot^{k}\epsilon)
\\
  \sff\conv\sfg
     &(x_1,\ldots,x_k)
      = \sff(x_{1(1)},\ldots,x_{k(1)}) \sfg(x_{1(2)},\ldots,x_{k(2)})
\end{align}
The monoid is commutative or associative if and only if $\mul,\comul_{\HBI{k}}$
are bicommutative or are biassociative.
\myeenv
In graphical terms the convolution multiplication and unit are illustrated by
\begin{align}
   \begin{pic}
      \node (i1) at (0,1) {};
      \node (i2) at (1,1) {};
      \node at (0.5,0.8) {$\ldots$ };
      \node[rectangle,draw,minimum width=1.3cm,minimum height=18pt] (fg) at (0.5,0) {$\sff\conv\sfg$ };
      \node (o) at (0.5,-1) {};
      \draw[thick] (i1.center) to (i1.center |- fg.north west);
      \draw[thick] (i2.center) to (i2.center |- fg.north east);
      \draw[thick] (fg.south) to (o.center);
   \end{pic}
\cong
   \begin{pic}
      \node (i1) at (0.5,1) {};
      \node at (1,0.8) {$\ldots$ };
      \node (i2) at (1.5,1) {};
      \node (u1) at (0.5,0.75) {};
      \node (u2) at (1.5,0.75) {};
      \node[rectangle,draw,minimum width=1cm] (f) at (0.25,-0.25) {$\sff $};
      \node[rectangle,draw,minimum width=1cm] (g) at (1.75,-0.25) {$\sfg $};
      \node (d) at (1,-0.75) {};
      \node (o) at (1,-1) {};
      \draw[thick] (i1.center) to (u1.center);
      \draw[thick] (i2.center) to (u2.center);
      \draw[thick] (d.center) to (o.center);
      \draw[thick,out=180,in=90] (u1.center) to (f.145);
      \draw[thick,out=  0,in=90] (u1.center) to (g.145);
      \draw[thick,out=180,in=90] (u2.center) to (f.35);
      \draw[thick,out=  0,in=90] (u2.center) to (g.35);
      \draw[thick,out=315,in=180] (f.south) to (d.center);
      \draw[thick,out=225,in=  0] (g.south) to (d.center);
      \node[yshift=0.2cm] at (f.north) {$\ldots$};
      \node[yshift=0.2cm] at (g.north) {$\ldots$};
   \end{pic}
%%%%%%%%%%%%%%%%%%%%%%%%%%%%%%%%%%%%%%%%%%%%%%%%%%%%%%%%%%%%%%%%%%%%%%%%
;\hskip1cm
   \begin{pic}
      \node (i1) at (0,1) {};
      \node (i2) at (1,1) {};
      \node at (0.5,0.8) {$\ldots$ };
      \node[rectangle,draw,minimum width=1.3cm] (fg) at (0.5,0) {$\sfe$ };
      \node (o) at (0.5,-1) {};
      \draw[thick] (i1.center) to (i1.center |- fg.north west);
      \draw[thick] (i2.center) to (i2.center |- fg.north east);
      \draw[thick] (fg.south) to (o.center);
   \end{pic}
\cong
   \begin{pic}
      \node (i1) at (0,1) {};
      \node (i2) at (1,1) {};
      \node at (0.5,0.75) {$\ldots$ };
      \node[circle,inner sep=2pt,draw,thick] (eps1) at (0,0.3) {};
      \node[xshift=3pt] at (eps1.east) {$\epsilon$};
      \node[circle,inner sep=2pt,draw,thick] (eps2) at (1,0.3) {};
      \node[xshift=3pt] at (eps2.east) {$\epsilon$};
      \node[circle,inner sep=2pt,draw,thick] (eta) at (0.5,-0.3) {};
      \node[xshift=3pt] at (eta.east) {$\eta$};
      \node (o) at (0.5,-1) {};
      \draw[thick] (i1.center) to (eps1);
      \draw[thick] (i2.center) to (eps2);
      \draw[thick] (eta) to (o.center);
   \end{pic}
\end{align}
We could also allow for more than one output string, but for what
follows the present setup is general enough.

As said above, we are interested in $k$-cochains, which are invertible
$k$-$1$-maps. Using unitality and a normalization condition
the Milnor--Moore recursive inverse is still available via cut
coproducts as shown in~\eqref{grph:inverse}.
\mybenv{Lemma}\label{lem-inverseCochain}
Let $\sff$ be unital, $\sff(1,\ldots,1)=1$, and normalized,
$f(x_1,\ldots,x_n)=0$ for at least one $x_{i}\in\H^{0}$ and
one other $x_{j}\in\H^{+}$. The \emph{convolutive} inverse
$\overline{\sff}\in \hom(\HBI{k},\HB)$ of a $k$-cochain $\sff$ is given
by the Milnor--Moore recursive formula, which in graphical terms is
given by:
\begin{align}
   \begin{pic}
      \node (i1) at (0,1) {};
      \node[iso2cell,minimum width=5pt] (c1) at (0,0.75) {};
      \node (i2) at (1,1) {};
      \node[iso2cell,minimum width=5pt] (c2) at (1,0.75) {};
      \node at (0.5,0.8) {$\ldots$ };
      \node[rectangle,draw,minimum width=1.1cm,minimum height=17pt] (of) at (0.5,0) {$\overline{\sff}$ };
      \node (o) at (0.5,-1) {};
      \draw[thick] (i1.center) to (c1.north) (c1.south) to (i1 |- of.north);
      \draw[thick] (i2.center) to (c2.north) (c2.south) to (i2 |- of.north);
      \draw[thick] (of.south) to (o.center);
   \end{pic}
\cong
   %\begin{pic}
      %\node (i1) at (0,1) {};
      %\node (i2) at (1,1) {};
      %\node at (0.5,0.75) {$\ldots$ };
      %\node[circle,inner sep=2pt,draw,thick] (eps1) at (0,0.3) {};
      %\node[xshift=3pt] at (eps1.east) {$\epsilon$};
      %\node[circle,inner sep=2pt,draw,thick] (eps2) at (1,0.3) {};
      %\node[xshift=3pt] at (eps2.east) {$\epsilon$};
      %\node[circle,inner sep=2pt,draw,thick] (eta) at (0.5,-0.3) {};
      %\node[xshift=3pt] at (eta.east) {$\eta$};
      %\node (o) at (0.5,-1) {};
      %\draw[thick] (i1.center) to (eps1);
      %\draw[thick] (i2.center) to (eps2);
      %\draw[thick] (eta) to (o.center);
   %\end{pic}
-
   \begin{pic}
      \node (i1) at (0,1) {};
      \node[iso2cell,minimum width=5pt] (c1) at (0,0.75) {};      
      \node (i2) at (1,1) {};
      \node[iso2cell,minimum width=5pt] (c2) at (1,0.75) {};      
      \node at (0.5,0.8) {$\ldots$ };
      \node[rectangle,draw,minimum width=1.1cm,minimum height=17pt] (of) at (0.5,0) {$\sff$ };
      \node (o) at (0.5,-1) {};
      \draw[thick] (i1.center) to (c1.north) (c1.south) to (i1 |- of.north);
      \draw[thick] (i2.center) to (c2.north) (c2.south) to (i2 |- of.north);
      \draw[thick] (of.south) to (o.center);
   \end{pic}
-
   \begin{pic}
      \node (i1) at (0.5,1.25) {};
      \node[iso2cell,minimum width=5pt] (c1) at (0.5,1) {};
      \node at (1,0.9) {$\ldots$ };
      \node (i2) at (1.5,1.25) {};
      \node[iso2cell,minimum width=5pt] (c2) at (1.5,1) {};      
      \node (u1) at (0.5,0.75) {};
      \node (u2) at (1.5,0.75) {};
      \node[iso2cell,minimum width=5pt] (iso1a) at (-0.125,0.25) {};
      \node[iso2cell,minimum width=5pt] (iso1b) at (0.625,0.25) {};
      \node[iso2cell,minimum width=5pt] (iso2a) at (1.325,0.25) {};
      \node[iso2cell,minimum width=5pt] (iso2b) at (2.125,0.25) {};      
      \node[rectangle,draw,minimum width=1cm,minimum height=17pt] (f) at (0.25,-0.25) {$\overline{\sff} $};
      \node[rectangle,draw,minimum width=1cm,minimum height=17pt] (g) at (1.75,-0.25) {$\sff $};
      \node (d) at (1,-0.75) {};
      \node (o) at (1,-1) {};
      \draw[thick] (i1.center) to (c1.north) (c1.south) to (u1.center);
      \draw[thick] (i2.center) to (c2.north) (c2.south) to (u2.center);
      \draw[thick] (d.center) to (o.center);
      \draw[thick,out=180,in=90] (u1.center) to (iso1a.north);
      \draw[thick,out=  0,in=90] (u1.center) to (iso2a.north);
      \draw[thick,out=180,in=90] (u2.center) to (iso1b.north);
      \draw[thick,out=  0,in=90] (u2.center) to (iso2b.north);
      \draw[thick] (iso1a.south) to (iso1a.south |- f.north);
      \draw[thick] (iso2a.south) to (iso2a.south |- g.north);
      \draw[thick] (iso1b.south) to (iso1b.south |- f.north);
      \draw[thick] (iso2b.south) to (iso2b.south |- g.north);
      \draw[thick,out=315,in=180] (f.south) to (d.center);
      \draw[thick,out=225,in=  0] (g.south) to (d.center);
      \node[yshift=0.2cm] at (f.north) {$\ldots$};
      \node[yshift=0.2cm] at (g.north) {$\ldots$};
   \end{pic}
;&&&&
   \begin{pic}
      \node[circle,inner sep=2pt,draw,thick] (i1) at (0,0.75) {};
      \node[xshift=3pt] at (i1.east) {$\eta$};
      \node[circle,inner sep=2pt,draw,thick] (i2) at (1,0.75) {};
      \node[xshift=3pt] at (i2.east) {$\eta$};
      \node at (0.5,0.5) {$\ldots$ };
      \node[rectangle,draw,minimum width=1.1cm] (of) at (0.5,0) {$\overline{\sff}$ };
      \node (o) at (0.5,-1) {};
      \draw[thick] (i1) to (i1 |- of.north);
      \draw[thick] (i2) to (i2 |- of.north);
      \draw[thick] (of.south) to (o.center);
   \end{pic}
\cong
   \begin{pic}
      \node[circle,inner sep=2pt,draw,thick] (i1) at (0,0.75) {};
      \node[xshift=3pt] at (i1.east) {$\eta$};
      \node (o) at (0,-1) {};
      \draw[thick] (i1) to (o.center);
   \end{pic}
;&&&&
   \begin{pic}
      \node (i1) at (0,1) {};
      \node[iso2cell,minimum width=5pt] (c1) at (0,0.75) {};
      \node (i2) at (1.5,1) {};
      \node[iso2cell,minimum width=5pt] (c2) at (1.5,0.75) {};
      \node at (0.4,0.5) {$\ldots$ };
      \node at (1.1,0.5) {$\ldots$ };
      \node[circle,inner sep=2pt,draw,thick] (i3) at (0.75,0.8) {};
      \node[xshift=3pt] at (i3.east) {$\eta$};
      \node[rectangle,draw,minimum width=1.7cm] (of) at (0.75,0) {$\overline{\sff}$ };
      \node (o) at (0.75,-1) {};
      \draw[thick] (i1.center) to (c1.north) (c1.south) to (i1 |- of.north);
      \draw[thick] (i2.center) to (c2.north) (c2.south) to (i2 |- of.north);
      \draw[thick] (i3.south) to (i3 |- of.north);
      \draw[thick] (of.south) to (o.center);
   \end{pic}
=0      
\end{align}
Again, it would be sufficient to assume that $\sff$ is normalized and
invertible in $\Bbbk$ at the cost of a more complicated inversion
formula.
\myeenv

For $k$-$1$-cochains with $k>1$ there is no obvious candidate for an
\emph{antipode} as there is no identity morphism $\Id$ for $k$-$l$
tangles with $k\not=l$. However, we can invert for example the associative
multiplication $\mul^{(k)}$ of our ambient Hopf algebra. In the $k=2$
case we find that $\overline{\mul} = \mul\circ(\antip\ot\antip)$, as
will be shown in Proposition~\ref{prop-twoInverses}.
%-----------------------------------------------------------------------
\subsection{Hopf algebra cohomology}
%-----------------------------------------------------------------------
The main application that we have in mind for convolution monoids is that
of Hopf algebra deformations, involving deformed binary multiplications
which remain associative. The appropriate formulation of conditions for
this emerges out of Sweedler cohomology~\cite{sweedler:1968a} which we
now briefly sketch. In this language the binary multiplications are
classed as $2$-$1$ maps; clearly an extended construction exists for
$n$-ary multiplications, which we shall not consider further here. We
begin with the
\mybenv{Definition}\label{def-cochain}
A \emph{$k$-cochain} is an (unital) normalized invertible $k$-$1$-module
map $\sfc : \HBI{k} \longrightarrow \HB$, not necessarily an algebra
morphism. The space of function $k$-cochains is given by the
linear dual $\HBI{k,*} =\hom(\HBI{k},\Bbbk)$ restricting the maps to
the codomain $\Bbbk\cong\H^{0}$.
\myeenv
Note that any $2$-$0$ map $\sff$ can be promoted to a $2$-$1$ map
by composing it with the unit $\eta\circ\sff$. A cohomology theory of
algebras over a Hopf algebra was devised by Sweedler, from which we
just need the coboundary operator (or differential). One defines the
face operators $m_{i}$ and the degeneracy operators $s_{i}$
\begin{align}
  m_{i} &: \HBI{k+1}\rightarrow \HBI{k}
    :: (x_{1},\ldots,x_{k+1})\mapsto (x_{1},\ldots,x_{i}x_{i+1},\ldots,x_{k+1})
  &&&
  1\leq i \leq k
  \nn
  s_{i} &: \HBI{k+1}\rightarrow\HBI{k+2}
    :: (x_{1},\ldots,x_{k+1}) \mapsto (x_{1},\ldots,x_{i},1,x_{i+1},\ldots,x_{k+2})  
  &&&
  1\leq i\leq k+2
\end{align}
which are $\HB$-module morphisms and coalgebra morphisms.
\mybenv{Definition}\label{def-Sweedler-cohomology}
The coboundary operator (differential)
$\partial_{k-1} : \hom(\HBI{k},\HB) \rightarrow \hom(\HBI{k+1},\HB)$ is
given by the coface maps $\partial_{k-1}^{i}$ as follows
\begin{align}\label{eq-SweddlerCohomology}
\partial_{k-1}^i \sfc_{k}(x_1,\ldots,x_{k+1})
  &:=
    \left\{
     \begin{array}{ll}
       \epsilon( x_1)\sfc_{k}(x_2,\ldots,x_{k+1})  & i=0 \\
       \sfc_{k}( x_1,\ldots,x_i x_{i+1},\ldots,x_{k+1})  & i\not= 0,k+1 \\
       \sfc_{k}( x_1,\ldots,x_{k})\epsilon(x_{k+1})  & i=k+1
     \end{array}
    \right.
  \nonumber \\
(\partial_{k-1} \sfc_{k})_{k+1} \equiv \partial_{k-1} \sfc_{k}
  &:= \partial^0_{k-1} \sfc_{k} \conv 
      \partial^1_{k-1} \sfc^{-1}_{k} \conv 
      \partial^2_{k-1} \sfc_{k} \conv{}
      \ldots\conv{}
      \partial^{k+1}_{k-1} \sfc^{\pm 1}_{k}
   \\
\partial_{k} \partial_{k-1}
  &= \sfe^{k+2} = \eta(\epsilon\otimes \ldots \otimes \epsilon)
  \nonumber
\end{align}
\myeenv
Sweedler cohomology is written multiplicatively as opposed to an ordinary
cohomology with $\mathrm{d}^2=0$ and coface maps joined additively using
alternating signs. Here we have alternating maps and inverses and the
connecting element is the Abelian convolution multiplication. We will
drop degree indices whenever they are clear from the context.
\mybenv{Definition}\label{def-cocycle}
A $k$-cochain $\sfc$ is a $k$-cocycle if $\partial\sfc=\sfe$.
\myeenv

If one wants to avoid convolutive inverses, then the $k$-cocycle property
for a $k$-cochain $\sfc$ can be rewritten as an identity and then
extended to all $k$-cochains. Using the face maps $m_{i}$ a
translation of the cocycle condition $\partial\sfc=\sfe$ implied
by~\eqref{eq-SweddlerCohomology} is given by
\begin{align}\label{eq-algCocycleId}
k &\textrm{~even} &:&&
(\prod_{\conv,{i:\textrm{odd}}} (\sfc\circ m_{i}))
  \conv (\sfc\ot\epsilon)
  &=
  (\epsilon\ot\sfc)\conv
  (\prod_{\conv,{i:\textrm{even}}} (\sfc\circ m_{i}))
\nn
k &\textrm{~odd} &:&&
(\prod_{\conv,{i:\textrm{odd}}} (\sfc\circ m_{i}))
  &=
  (\epsilon\ot\sfc)\conv
  (\prod_{\conv,{i:\textrm{even}}} (\sfc\circ m_{i}))
  \conv (\sfc\ot\epsilon)\,.
\end{align}

In graphical terms, a $k$-cochain $\sfc$ is a map with $k$ input lines
and one output line, its convolutive inverse is the $k$-cochain denoted
by $\sfoc$. Looking at the case of a 1-cochain $\sff$, the above face
maps $\partial^i_{0}\sff$ are given by the 2-1-maps $\epsilon\ot \sff$,
$\sfof\circ \mul_{1}$ and $\sff\ot\epsilon$.
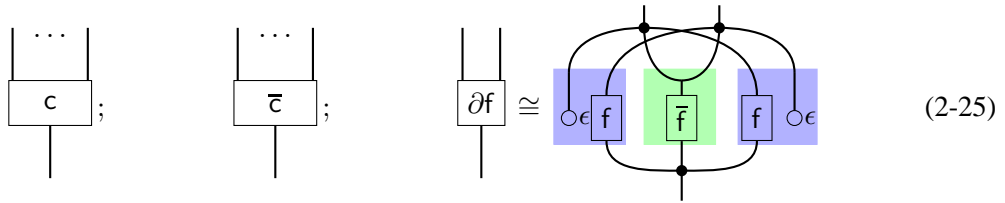
\begin{align}
   \begin{pic}
      \node (i1) at (0,1) {};
      \node at (0.5,0.9) {$\ldots$ };
      \node (i2) at (1,1) {};
      \node[rectangle,draw,minimum width=1.1cm,minimum height=17pt] (c) at (0.5,0) {$\sfc$ };
      \node (o) at (0.5,-1) {};
      \draw[thick] (i1.center) to (i1 |- c.north);
      \draw[thick] (i2.center) to (i2 |- c.north);
      \draw[thick] (o |-c.south) to (o.center);
   \end{pic}
;\hskip1.5cm
   \begin{pic}
      \node (i1) at (0,1) {};
      \node at (0.5,0.9) {$\ldots$ };
      \node (i2) at (1,1) {};
      \node[rectangle,draw,minimum width=1.1cm,minimum height=17pt] (c) at (0.5,0) {$\sfoc$ };
      \node (o) at (0.5,-1) {};
      \draw[thick] (i1.center) to (i1 |- c.north);
      \draw[thick] (i2.center) to (i2 |- c.north);
      \draw[thick] (o |-c.south) to (o.center);
   \end{pic}
;\hskip1.5cm
%%%%%%%%%%%%%%%%%%%%%%%%%%%%%%%%%%%%%%%%%%%%%%%%%%%%%%%%%%%%%%%%%%%%%%%%
%%  2-cocycle
   \begin{pic}
      \node (i1) at (0,1) {};
      \node (i2) at (0.5,1) {};
      \node[rectangle,draw,minimum height=17pt] (c) at (0.25,0) {$\partial\sff$ };
      \node (o) at (0.25,-1) {};
      \draw[thick] (i1.center) to (i1 |- c.north);
      \draw[thick] (i2.center) to (i2 |- c.north);
      \draw[thick] (o |- c.south) to (o.center);
   \end{pic}
\cong
   \begin{pic}
      \node (i1) at (1,1.3) {};
      \node (i2) at (2,1.3) {};
      \node[circle,fill=black,inner sep=1.5pt] (u1) at (1,1) {};
      \node[circle,fill=black,inner sep=1.5pt] (u2) at (2,1) {};
      \node (m1) at (0,0.4) {};
      \node (m2) at (1.5,0.3) {};
      \node (m3) at (3,0.4) {};
      \node[circle,inner sep=2pt,draw] (eps1) at (0,-0.2) {};
      \node[xshift=3pt] at (eps1.east) {$\epsilon$};
      \node[rectangle,draw,minimum height=17pt] (c1) at (0.5,-0.2) {$\sff$ };
      \node[rectangle,draw,minimum height=17pt] (c2) at (1.5,-0.2) {$\sfof$ };
      \node[rectangle,draw,minimum height=17pt] (c3) at (2.5,-0.2) {$\sff$ };
      \node[circle,inner sep=2pt,draw] (eps2) at (3,-0.2) {};
      \node[xshift=3pt] at (eps2.east) {$\epsilon$};
      \node[circle,fill=black,inner sep=1.5pt] (d) at (1.5,-0.9) {};
      \node (o) at (1.5,-1.3) {};
      \draw[thick] (i1.center) to (u1.center);
      \draw[thick] (i2.center) to (u2.center);
      \draw[thick,out=180,in=90] (u1.center) to (m1.center);
      \draw[thick] (m1.center) to (eps1);
      \draw[thick,out=270,in=180] (u1.center) to (m2.center);
      \draw[thick] (m2.center) to (c2.north);
      \draw[thick,out=0,in=90] (u1.center) to (c3.north);
      \draw[thick,out=180,in=90] (u2.center) to (c1.north);
      \draw[thick,out=270,in= 0] (u2.center) to (m2.center);
      \draw[thick,out=  0,in=90] (u2.center) to (m3.center);
      \draw[thick] (m3.center) to (eps2.north);
      \draw[thick,out=270,in=180] (c1.south) to (d.center);
      \draw[thick] (c2.south) to (d.center);
      \draw[thick,out=270,in=  0] (c3.south) to (d.center);
      \draw[thick] (d.center) to (o.center);
      \begin{pgfonlayer}{background}
         \draw[draw,fill,color=blue!30] (-0.2,-0.55) rectangle (0.75cm,0.45cm);
         \draw[draw,fill,color=green!30,xshift=1.2cm] (-0.2,-0.55) rectangle (0.75cm,0.45cm);
         \draw[draw,fill,color=blue!30,xshift=2.45cm] (-0.2,-0.55) rectangle (0.85cm,0.45cm);
      \end{pgfonlayer}
   \end{pic}
\end{align}
In the rightmost tangle, we have for better readability marked the two
comultiplications and the multiplication with black dots, the other
incidence points are just crossings. The tangle reads as follows: the
two input lines are comultiplied by
$\Delta^{(3)} =(\Id\ot\Delta)\circ\Delta$ to provide 3 strands\footnote{%
   We use the convention $\comul^{(0)}=\epsilon$, $\comul^{(1)}=\Id$,
   $\comul^{(2)}=\comul$, $\comul^{(3)}=(\comul\ot\Id)\circ\comul = 
   (\Id\ot\comul)\circ\comul$ etc. which is compliant with notation used
   for `loop operators' in Appendix~\ref{FGLloop}. As a mnemonic we count
   the output lines.}
and then shuffled by $\Sigma_{2}\circ(1\ot\Sigma_{1}\ot 1)$ so that
the outputs of the left comultiplication are fed into the left inputs of
the face maps $\partial^i \sff^{\pm}$ while the outputs of the right
comultiplication are fed into the right inputs of these face maps. 

Note that a $2$-coboundary $\sfg=\partial\sff$ is
automatically a $2$-cocycle, due to $\partial^2=\sfe$. The condition for
a $1$-cochain to be a $1$-cocyle translates by~\eqref{eq-algCocycleId}
into the statement that $\sff$ is an algebra homomorphism.
\begin{align}\label{eq-oneCocycleIsAlgHom}
   \begin{pic}
      \node (i1) at (0,1) {};
      \node (i2) at (0.5,1) {};
      \node[rectangle,draw] (c) at (0.25,0) {$\partial\sff$ };
      \node (o) at (0.25,-1) {};
      \draw[thick] (i1.center) to (i1 |- c.north);
      \draw[thick] (i2.center) to (i2 |- c.north);
      \draw[thick] (o |- c.south) to (o.center);
   \end{pic}
&\cong
   \begin{pic}
      \node (i1) at (0,1) {};
      \node (i2) at (0.5,1) {};
      \node[circle,inner sep=2pt,draw,thick] (eps1) at (0,0.2) {};
      \node[xshift=3pt] at (eps1.east) {$\epsilon$};
      \node[circle,inner sep=2pt,draw,thick] (eps2) at (0.5,0.2) {};
      \node[xshift=3pt] at (eps2.east) {$\epsilon$};
      \node[circle,inner sep=2pt,draw,thick] (eta) at (0.25,-0.2) {};
      \node[xshift=3pt] at (eta.east) {$\eta$};
      \node (o) at (0.25,-1) {};
      \draw[thick] (i1.center) to (eps1);
      \draw[thick] (i2.center) to (eps2);
      \draw[thick] (eta) to (o.center);
   \end{pic}
&&\Leftrightarrow &&
   \begin{pic}
      \node (i1) at (0,1) {};
      \node (i2) at (1,1) {};
      \node (a) at (0.5,0.5) {};
      \node[rectangle,thick,draw] (p) at (0.5,-0.25) {$\sff$ };
      \node (o) at (0.5,-1) {};
      \draw[thick,out=270,in=180] (i1.center) to (a.center);
      \draw[thick,out=270,in=  0] (i2.center) to (a.center);
      \draw[thick] (a.center) to (p.north);
      \draw[thick] (p.south) to (o.center);
   \end{pic}
\cong\quad
   \begin{pic}
      \node (i1) at (0,1) {};
      \node (i2) at (1,1) {};
      \node[rectangle,thick,draw] (p1) at (0,0.4) {$\sff$ };
      \node[rectangle,thick,draw] (p2) at (1,0.4) {$\sff$ };
      \node (a) at (0.5,-0.25) {};
      \node (o) at (0.5,-1) {};
      \draw[thick] (i1.center) to (p1.north);
      \draw[thick] (i2.center) to (p2.north);
      \draw[thick,out=270,in=180] (p1.south) to (a.center);
      \draw[thick,out=270,in=  0] (p2.south) to (a.center);
      \draw[thick] (a.center) to (o.center);
   \end{pic}     
\end{align}
 
Now let $\sff$ be a generic 2-cochain. We want to know when it is a
2-cocycle. To simplify the graphics we box the face maps and color code
the inverse as done above, but do not show the internal lines. We
let $\mul_{1} = \mul\ot 1$ and $\mul_{2}=1\ot \mul$ to indicate where
the multiplication takes place. The 2-cochain $\sff$ is a \emph{2-cocycle}
if and only if $\partial\sff=\sfe = \eta\circ\epsilon\ot\epsilon\ot\epsilon$.
In terms of diagrams
\begin{align}
   \begin{pic}
      \node (i1) at (0,1.5) {};
      \node (i2) at (0.5,1.5) {};
      \node (i3) at (1,1.5) {};
      \node[rectangle,draw,minimum width=1.1cm,minimum height=19pt] (c) at (0.5,-0.1) {$\partial\sff$ };
      \node (o) at (0.5,-1.5) {};
      \draw[thick] (i1.center) to (i1 |- c.north);
      \draw[thick] (i2.center) to (i2 |- c.north);
      \draw[thick] (i3.center) to (i3 |- c.north);
      \draw[thick] (o |- c.south) to (o.center);
   \end{pic}
\cong
   \begin{pic}
      \node (i1) at (1.25,1.5) {};
      \node (i2) at (2.25,1.5) {};
      \node (i3) at (3.25,1.5) {};
      \node (u1) at (1.25,1.3) {};
      \node (u2) at (2.25,1.3) {};
      \node (u3) at (3.25,1.3) {};
      \node[rectangle,draw,fill=blue!30,minimum height=19pt] (c1) at (0,-0.1) {$\epsilon\ot\sff$ };
      \node[rectangle,draw,fill=green!30,minimum height=19pt] (c2) at (1.4,-0.1) {$\sfof\circ m_{1}$ };
      \node[rectangle,draw,fill=blue!30,minimum height=19pt] (c3) at (3.05,-0.1) {$\sff\circ m_{2}$ };
      \node[rectangle,draw,fill=green!30,minimum height=19pt] (c4) at (4.5,-0.1) {$\sfof\ot\epsilon$ };
      \node (d) at (2.25,-1.3) {};
      \node (o) at (2.25,-1.5) {};
      \draw[thick] (i1.center) to (u1.center);
      \draw[thick] (i2.center) to (u2.center);
      \draw[thick] (i3.center) to (u3.center);
      \draw[thick,out=180,in=90] (u1.center) to (c1.north west);
      \draw[thick,out=210,in=90] (u1.center) to (c2.north west);
      \draw[thick,out=300,in=115] (u1.center) to (c3.north west);
      \draw[thick,out=340,in=135] (u1.center) to (c4.north west);
      \draw[thick,out=210,in=60] (u2.center) to (c1.north);
      \draw[thick,out=235,in=90] (u2.center) to (c2.north);
      \draw[thick,out=305,in=90] (u2.center) to (c3.north);
      \draw[thick,out=330,in=120] (u2.center) to (c4.north);
      \draw[thick,out=200,in=45] (u3.center) to (c1.north east);
      \draw[thick,out=240,in=65] (u3.center) to (c2.north east);
      \draw[thick,out=330,in=90] (u3.center) to (c3.north east);
      \draw[thick,out=0,in=90] (u3.center) to (c4.north east);
      \draw[thick,out=270,in=180] (c1.south) to (d.center);
      \draw[thick,out=270,in=155] (c2.south) to (d.center);
      \draw[thick,out=270,in= 25] (c3.south) to (d.center);
      \draw[thick,out=270,in=  0] (c4.south) to (d.center);
      \draw[thick] (d.center) to (o.center);
   \end{pic}
\cong
   \begin{pic}
      \node (i1) at (0,1.5) {};
      \node (i2) at (0.5,1.5) {};
      \node (i3) at (1,1.5) {};
      \node[circle,inner sep=2pt,draw,thick] (eps1) at (0,0.3) {};
      \node[xshift=3pt] at (eps1.east) {$\epsilon$};
      \node[circle,inner sep=2pt,draw,thick] (eps2) at (0.5,0.3) {};
      \node[xshift=3pt] at (eps2.east) {$\epsilon$};
      \node[circle,inner sep=2pt,draw,thick] (eps3) at (1,0.3) {};
      \node[xshift=3pt] at (eps3.east) {$\epsilon$};
      \node[circle,inner sep=2pt,draw,thick] (eta) at (0.5,-0.3) {};
      \node[xshift=3pt] at (eta.east) {$\eta$};
      \node (o) at (0.5,-1.5) {};
      \draw[thick] (i1.center) to (eps1);
      \draw[thick] (i2.center) to (eps2);
      \draw[thick] (i3.center) to (eps3);
      \draw[thick] (eta) to (o.center);
   \end{pic}
\end{align}
As we work in a bicommutative setting one can invert the terms which
contain $\sfof$ to get the identity form~\eqref{eq-algCocycleId} of
the 2-cocycle condition $(\sff\circ\mul_{1})\conv(\sff\ot\epsilon)
=(\epsilon\ot\sff)\conv(f\circ\mul_{2})$. This condition is also
called the multiplicativity constraint, as it ensures the associativity
of a deformed product~\cite{sweedler:1968a}. Dealing with groups and
group characters, we need to stay associative, so this condition has
to be employed below.

In analogy to ordinary cohomology we have the following notation and
simple facts:
\begin{itemize}
\item A $k$-cochain $\sfc$ is a \emph{$k$-cocycle} if $\partial\sfc=\sfe$.
      Such cochains may also be called closed.
\item A $k$-cochain $\sfc$ is a \emph{$k$-coboundary} if there exists
      a $k-1$-cochain $\sfd$ such that $\sfc = \partial\sfd$. Such a
      cochain may also be called exact.
\item A 1-cochain $\sff \in \hom(\HB,\HB)$ is a 2-cocycle (closed
      $\partial\sff=\sfe$) if and only if $\sff$ is an $\mul$-algebra
      homomorphism ($\sff(xy) = \sff(x)\sff(y)$, see~\eqref{eq-oneCocycleIsAlgHom}).
\item The 2-cochain \emph{deformation} $\mul^{\prime}=\sfc\conv\mul$ of
      the multiplication $\mul$ is again an \emph{associative multiplication}
      if and only if $\sfc$ is a 2-cocycle $\partial\sfc=\sfe$.
      A proof for the supersymmetric case can be found
      in~\cite{brouder:fauser:frabetti:oeckl:2002a}. A similar but less
      obvious statement is found in~\cite{rota:stein:1994a}, see
      also~\cite{sweedler:1968a}. 
\item In Sections~\ref{sec-CircleProducts} and~\ref{sec-HashProd} we
      will discuss deformations in more detail. Here we note that a deformation
      $\mul^{\prime}=\sfc\conv\mul$ where $\sfc$ is a coboundary
      $\sfc=\partial\sfd$, belongs to the same cohomology class as
      $\mul$. For example, starting with a Grassmann algebra and
      deforming it with an antisymmetric bilinear form, which is
      automatically a 2-coboundary, results in an isomorphic Grassmann
      algebra, however the grading is changed~\cite{fauser:2001b,fauser:2002c}.
      In contrast, deforming with a symmetric bilinear form, a proper
      2-cocycle, deforms the Grassmann algebra into a Clifford algebra
      inducing a non trivial \emph{quantization}, \cite{fauser:2001b,fauser:2002c}.
\item Composing a $k$-cochain with a linear form enables the analogous
      cohomology to be formulated for ring-valued cochains in
      $\hom(\HBI{k},\Bbbk)$. \emph{Vice versa}, any such $\Bbbk$-valued
      cochain can be promoted to an algebra valued cochain by composing
      it with the unit map $\eta$.
\end{itemize}
%-----------------------------------------------------------------------
\subsection{The monoid of pairings and its Laplace subgroup}\label{subsec:MonoidOfPairings}
%-----------------------------------------------------------------------
As we are interested in deforming binary multiplications, we focus on
pairings. These are instances of general maps $\sfa : A\ot B \rightarrow C$,
see~\cite{schmitt:1999a}. However, in the present context we restrict
ourselves to $\hom(\HBI{2},\HB)$. Using these $2$-cochains or pairings for
deformations, we consider the monoid structure associated with
$2$-cocycles, and we will also need to furnish these pairings with
additional properties, most prominently a distributive property. Such
distributive pairings are also called Laplace pairings.

We denote pairings by sans serif letters from the beginning of the alphabet
$\sfa,\sfb,\ldots \in \hom(\HBI{2},\HB)$, and Hopf algebra elements
$x,y,z\in \H$ from the end of the alphabet. As usual, a pairing is called
unital if $\sfa\circ(\eta\ot\eta)=\eta$, and normalized if
$\sfa\circ(\Id\ot\eta) = \eta\circ\epsilon$ and
$\sfa\circ(\eta\ot\Id) = \eta\circ\epsilon$.

\mybenv{Definition}\label{def-Laplace-pairing}
A \emph{Laplace pairing} (or distributive pairing) for a convolution
algebra over a Hopf algebra $\HB$ fulfils the right and left
straightening law
\begin{align}\label{eq-LaplaceRight}
\sfa(x,yz) &= \sfa(x_{(1)},y)\sfa(x_{(2)},z)
&&&
   \begin{pic}
      \node (i1) at (0,0.75) {};
      \node (i2) at (0.5,0.75) {};
      \node (i3) at (1,0.75) {};
      \node (m) at (0.75,0.3) {};
      \node[circle,draw,inner sep=1pt] (a) at (0.375,-0.3) {$\sfa$ };
      \node (o) at (0.375,-0.75) {};
      \draw[thick,out=270,in=180] (i1.center) .. controls +(down:8mm) and +(left:3mm) .. (a.west);
      \draw[thick,out=270,in=180] (i2.center) to (m.center);
      \draw[thick,out=270,in=  0] (i3.center) to (m.center);
      \draw[thick,out=270,in=  0] (m.center) to (a.east);
      \draw[thick] (a.south) to (o.center);
   \end{pic}
\cong
   \begin{pic}
      \node (i1) at (0.25,0.75) {};
      \node (i2) at (0.75,0.75) {};
      \node (i3) at (1.25,0.75) {};
      \node (m) at (0.25,0.5) {};
      \node[circle,draw,inner sep=1pt] (a1) at (0.25,0) {$\sfa$ };
      \node[circle,draw,inner sep=1pt] (a2) at (1,0) {$\sfa$ };
      \node (d) at (0.6125,-0.5) {};
      \node (o) at (0.6125,-0.75) {};
      \draw[thick] (i1.center) to (m.center);
      \draw[thick,out=180,in= 90] (m.center) .. controls +(-0.3,0) and +(-0.3,0.2) .. (a1.west);
      \draw[thick,out=  0,in=180] (m.center) to (a2.west);
      \draw[thick,out=270,in=0] (i2.center) to (a1.east);
      \draw[thick,out=270,in=0] (i3.center) .. controls +(0,-0.4) and +(0.1,0) .. (a2.east);
      \draw[thick,out=270,in=180] (a1.south) to (d.center);
      \draw[thick,out=270,in=  0] (a2.south) to (d.center);
      \draw[thick] (d.center) to (o.center);
   \end{pic}
\\
%%%%%%%%%%%%%%%%%%%%%%%%%%%%%%%%%%%%%%%%%%%%%%%%%%%%%%%%%%%%%%%%%%%%%%%%
\label{eq-LaplaceLeft}
\sfa(xy,z) &= \sfa(x,z_{(1)})\sfa(y,z_{(2)})
&&&
   \begin{pic}
      \node (i1) at (0,0.75) {};
      \node (i2) at (0.5,0.75) {};
      \node (i3) at (1,0.75) {};
      \node (m) at (0.25,0.3) {};
      \node[circle,draw,inner sep=1pt] (a) at (0.6125,-0.3) {$\sfa$ };
      \node (o) at (0.6125,-0.75) {};
      \draw[thick,out=270,in=180] (i1.center) to (m.center);
      \draw[thick,out=270,in=  0] (i2.center) to (m.center);
      \draw[thick,out=270,in=  0] (i3.center) .. controls +(0,-0.5) and +(0.3,0) .. (a.east);
      \draw[thick,out=270,in=180] (m.center) to (a.west);
      \draw[thick] (a.south) to (o.center);
   \end{pic}
\cong
   \begin{pic}
      \node (i1) at (0,0.75) {};
      \node (i2) at (0.5,0.75) {};
      \node (i3) at (1,0.75) {};
      \node (m) at (1,0.5) {};
      \node[circle,draw,inner sep=1pt] (a1) at (0.25,0) {$\sfa$ };
      \node[circle,draw,inner sep=1pt] (a2) at (1,0) {$\sfa$ };
      \node (d) at (0.6125,-0.5) {};
      \node (o) at (0.6125,-0.75) {};
      \draw[thick,out=270,in=180] (i1.center) .. controls +(0,-0.4) and +(-0.1,0).. (a1.west);
      \draw[thick,out=270,in=180] (i2.center) to (a2.west);
      \draw[thick] (i3.center) to (m.center);
      \draw[thick,out=180,in=0] (m.center) to (a1.east);
      \draw[thick,out=  0,in=180] (m.center) .. controls +(0.3,0) and +(0.3,0.2) .. (a2.east);
      \draw[thick,out=270,in=180] (a1.south) to (d.center);
      \draw[thick,out=270,in=  0] (a2.south) to (d.center);
      \draw[thick] (d.center) to (o.center);
   \end{pic}
\end{align}
We shall assume that such pairings are normalized and unital. 
\myeenv
The above properties are synonymously called distributive laws or
straightening laws and play the role of expansion formulae as far as the
pairing structure is concerned. The name Laplace pairing commemorates
the fact that in the case of a Grassmann Hopf algebra, the above
expansions are the Laplace expansions of determinants in multiple rows
or columns (see~\cite{rota:stein:1994a} and references therein).

In the case of ring valued pairings, distributivity is part of the
definition of a pairing of bi- or Hopf algebras. Let $\sfA,\sfB$
be (here not necessarily commutative and/or cocommutative) bialgebras
or Hopf algebras.
\mybenv{Definition}\label{def-BialgebraPairing}
A bilinear pairing $\varphi : \sfA \times \sfB \rightarrow \Bbbk$ of
bialgebras (Hopf algebras) is called a \emph{bialgebra pairing} if i)
and ii) below are true. In the Hopf algebra case, $\varphi$ is called
a Hopf algebra pairing if in addition iii) also holds.
\begin{itemize}
\item[i)] $\varphi$ is normalized : $\varphi(a,1)=\epsilon_A(a)$ and
      $\varphi(1,b)=\epsilon_B(b)$.
\item[ii)] Distributivity or Laplace property :
     $\varphi(aa',b) = \sum_{(b)} \varphi(a,b_{(2)})\varphi(a',b_{(1)})$
     and
     $\varphi(a,bb')= \sum_{(a)} \varphi(a_{(1)},b)\varphi(a_{(2)},b')$
     where the order of Sweedler indices is important in the
     noncommutative case.
\item[iii)] Relation of antipodes $\varphi(\antip_A(a),b) =
     \varphi(a,\antip_B^{-1}(b))$.
\end{itemize}
\myeenv

Bialgebras or Hopf algebras $\sfA,\sfB$ fulfilling the above conditions
are called matched pairs. Such pairs play a central role in constructing
$\varphi$-Drinfeld doubles in Hopf algebra theory 
(see for example~\cite{kassel:1995a}). Laplace pairings can also be
understood in terms of $\HB$-module algebras and $\HB$-module{}
coalgebras, see~\cite{sweedler:1968a}.
%-----------------------------------------------------------------------

Having introduced the distributive
laws~(\ref{eq-LaplaceRight},\ref{eq-LaplaceLeft}) (or Laplace expansions),
we return to the bi-commutative case. However, in contrast to the above
Definition~\ref{def-BialgebraPairing} of a ring valued Hopf algebra
pairing, we revert to the general context of \emph{algebra valued}
pairings $\sfa\in\hom(\HBI{2},\HB)$, which is in the spirit of Sweedler
cohomology and (as we shall see below) the Rota--Stein deformation process.

\mybenv{Proposition}\label{prop-laplaceImpliesCoycle}
If a general 2-cochain $\sfc$ is a Laplace pairing, then it is a
2-cocycle $\partial\sfc=\sfe$.
\myeenv

\noindent
\textbf{Proof:}
We use the 2-cocycle condition in the form without
inverses~\eqref{eq-algCocycleId}. The right hand side reads
\begin{align}
  (\epsilon\ot\sfc )\conv (\sfc\circ\mul_{2})(x,y,z)
    &=\epsilon(x_{(1)})\sfc (y_{(1)},z_{(1)})\sfc (x_{(2)},y_{(2)}z_{(2)}) \nn
    &=\sfc (y_{(1)},z_{(1)})\sfc (x_{(1)},y_{(2)}) \sfc (x_{(2)},z_{(2)}) \,,
\end{align}
the left hand side reads
\begin{align}
  (\sfc\circ\mul_{1})\conv (\sfc\ot\epsilon )(x,y,z)
    &=\sfc (x_{(1)}y_{(1)},z_{(1)})\sfc (x_{(2)},y_{(2)})\epsilon (z_{(2)}) \nn
    &=\sfc (x_{(1)},z_{(1)})\sfc (y_{(1)},z_{(2)}) \sfc (x_{(2)},y_{(2)}) \,.
\end{align}
Using commutativity of the multiplication and the comultiplication shows
that these expressions are identical.
\qed

\mybenv{Definition} A derived Laplace pairing $\sfa_\phi$ is defined to be
the composition of a Laplace pairing $\sfa$ with a normalized 1-cocycle
$\phi$, $\sfa_\phi=\phi\circ\sfa$.
\begin{align}\label{grph-derivedLaplacePairing}
   \begin{pic}
      \node (i1) at (0,0.75) {};
      \node (i2) at (1,0.75) {};
      \node[circle,draw,inner sep=1pt] (a) at (0.5,0) {$\sfa_\phi$ };
      \node (o) at (0.5,-0.75) {};
      \draw[thick,out=270,in=180] (i1.center) .. controls +(down:8mm) and +(left:2mm) .. (a.west);
      \draw[thick,out=270,in=180] (i2.center) .. controls +(down:8mm) and +(left:-2mm) .. (a.east);
      \draw[thick] (a.south) to (o.center);
   \end{pic}
&\cong
   \begin{pic}
      \node (i1) at (0,0.75) {};
      \node (i2) at (1,0.75) {};
      \node[circle,draw,inner sep=1pt] (a) at (0.5,0.25) {$\sfa$ };
      \node[rectangle,thick,draw,inner sep=1pt] (p) at (0.5,-0.25) {$\phi$ };
      \node (o) at (0.5,-0.75) {};
      \draw[thick,out=270,in=180] (i1.center) .. controls +(down:4mm) and +(left:2mm) .. (a.west);
      \draw[thick,out=270,in=180] (i2.center) .. controls +(down:4mm) and +(left:-2mm) .. (a.east);
      \draw[thick] (a.south) to (p.north);
      \draw[thick] (p.south) to (o.center);
   \end{pic}
;\quad
   %\begin{pic}
      %\node (i1) at (0,0.75) {};
      %\node (i2) at (1,0.75) {};
      %\node (a) at (0.5,0.25) {};
      %\node[rectangle,thick,draw,inner sep=1pt] (p) at (0.5,-0.25) {$\phi$ };
      %\node (o) at (0.5,-0.75) {};
      %\draw[thick,out=270,in=180] (i1.center) to (a.center);
      %\draw[thick,out=270,in=  0] (i2.center) to (a.center);
      %\draw[thick] (a.center) to (p.north);
      %\draw[thick] (p.south) to (o.center);
   %\end{pic}
%&\cong\quad
   %\begin{pic}
      %\node (i1) at (0,0.75) {};
      %\node (i2) at (1,0.75) {};
      %\node[rectangle,thick,draw,inner sep=1pt] (p1) at (0,0.25) {$\phi$ };
      %\node[rectangle,thick,draw,inner sep=1pt] (p2) at (1,0.25) {$\phi$ };
      %\node (a) at (0.5,-0.25) {};
      %\node (o) at (0.5,-0.75) {};
      %\draw[thick] (i1.center) to (p1.north);
      %\draw[thick] (i2.center) to (p2.north);
      %\draw[thick,out=270,in=180] (p1.south) to (a.center);
      %\draw[thick,out=270,in=  0] (p2.south) to (a.center);
      %\draw[thick] (a.center) to (o.center);
   %\end{pic}
%;\quad
%\begin{pic}
      %\node[circle,inner sep=2pt,draw,thick] (eta) at (0,0.5) {};
      %\node[xshift=3pt] at (eta.east) {$\eta$};
      %\node[rectangle,thick,draw,inner sep=1pt] (phi) at (0,-0.2) {$\phi$};
      %\node (o) at (0,-0.75) {};
      %\draw[thick] (eta.south) to (phi.north);
      %\draw[thick] (phi.south) to (o.center);
%\end{pic}
%\cong
%\begin{pic}
      %\node[circle,inner sep=2pt,draw,thick] (eta) at (0,0.5) {};
      %\node[xshift=3pt] at (eta.east) {$\eta$};
      %\node (o) at (0,-0.75) {};
      %\draw[thick] (eta.south) to (o.center);
%\end{pic}
&&&
\parbox[c]{0.35\textwidth}{derived Laplace pairing $\sfa_\phi$.}
\end{align}
\vskip-3ex
\myeenv

We have seen that $k$-cochains form a monoid under the convolution
multiplication with unit $\sfe$. We show now that the (unital normalized)
derived Laplace pairings form a subgroup.

\mybenv{Proposition}\textbf{[derived pairings stay Laplace]}
Let $\sfa$ be Laplace pairing and $\phi$ a 1-cocycle, then the
\emph{derived pairing} $\sfa_\phi:=\phi\circ\sfa$ is a Laplace pairing.
\myeenv

\noindent\textbf{Proof:} The 1-cocycle property for $\phi$ states
that $\phi$ is an $\mul$-algebra homomorphism for the algebra structure
of the ambient Hopf algebra~\eqref{eq-oneCocycleIsAlgHom}. Starting with
the left hand side of~\eqref{eq-LaplaceRight} using the Laplace property of $\sfa${}
and the algebra homomorphism property of $\phi$ shows that $\sfa_\phi${}
is right Laplace. The other case is identical.
\qed

\mybenv{Proposition}\textbf{[inverse]}
A unital normalized derived Laplace pairing $\sfa_\phi$ is invertible
with inverse $\sfoa_\phi$, that is $\sfa_{\phi}\conv\sfoa_{\phi}=\sfe^2$.
The inverse is given by the Milnor--Moore recursion of
Lemma~\ref{lem-inverseCochain}
$\sfoa_\phi = -(\sfa_\phi -\sfe ) - \mul\circ(\sfa_\phi\ot\sfoa_\phi)
\circ\comul'_{\HBI{2}}$, where $\comul'_{\HBI{2}}$ is the proper cut
comultiplication on $\HBI{2}$. The inverse pairing $\sfoa_\phi$ is Laplace.
\myeenv

\noindent
\textbf{Proof:}
Again we need in general only that the pairing is invertible in $\Bbbk$
for $\sfa(1,1)$, but then the recursion is more complicated. 
The first part of the proposition is a corollary of the invertibility of
general unital normalized 2-cochains. Composition with a normalized
1-cocycle $\phi$ does not change the Milnor--Moore recursive inverse,
as $\phi\circ\eta=\eta$ and in the convolution $\phi$ acts just as a composition
$\sfa_\phi\conv\sfoa_\phi= \phi\circ(\sfa\conv\sfoa) = (\sfa\conv\sfoa)_\phi$.
To show that the inverse is Laplace, one first shows that 
$\sfa\circ \mul_{i}$ is invertible with inverse 
$\sfoa\circ\mul_{i}$, inverting one side of the Laplace condition.
For the other side one shows directly that
$\sfa_{\phi}(x_{(1)},y)\sfa_{\phi}(x_{(2)},z)$ has the inverse
$\sfoa_{\phi}(x_{(1)},y)\sfoa_{\phi}(x_{(2)},z)$,
and similarly for the other expansion law, hence $\sfoa$ is Laplace.
\qed

\mybenv{Proposition}\textbf{[multiplicative closure]}\label{prop-Closure}
Let $\sfa,\sfb$ be two Laplace pairings. The convolution product
$\sfc=\sfa\conv\sfb$ is again a Laplace pairing. If $\sfa_\phi$
and $\sfb_\psi$ are derived Laplace pairings such that
$\sfc=\sfa_\phi\conv\sfb_\psi$ we get the Laplace
expansion $\sfc(x,yz) = \sfc(x_{(1)},y)\sfc(x_{(2)},z)$
and its left counterpart.
\myeenv

\noindent
\textbf{Proof:}
Insert the definition $c=\sfa\conv\sfb$ into one of the Laplace
expansions~(\ref{eq-LaplaceRight},\ref{eq-LaplaceLeft}), use the
bialgebra law~\eqref{grph-bialgebra-and-antipode}, reorganize
the diagram using associativity and commutativity and use the
definition of $\sfc$ again. In the case of a derived Laplace pairing
$\sfc=\sfa_\phi\conv\sfb_\psi$ we get in the same way the
expansion law $\sfc(x,yz) = \sfc(x_{(1)},y)\sfc(x_{(2)},z)$
and similarly for the left expansion.
\qed

The previous results are summarized by the
\mybenv{Theorem}\textbf{[Laplace subgroup]}
Unital normalized derived Laplace pairings form the \emph{Laplace subgroup}
of the convolution monoid of Laplace pairings, and of the convolution
monoid of 2-cochains.
\myeenv
%-----------------------------------------------------------------------
\subsection{Frobenius Laplace pairings}\label{subsec:FobeniusLaplace}
%-----------------------------------------------------------------------
We can impose further a \emph{Frobenius} condition on our pairings, which
will be of central importance in our applications to group characters in
Section~\ref{sec-Applications}. However it is interesting in its own
right in the abstract setting studied in this section. 
We equip a Laplace pairing $\sfa$ with the \emph{additional} condition
that it carries a commutative Frobenius algebra structure.

Frobenius algebras can be characterized in different ways, in graphical
terms it is convenient to use a coalgebra
structure~\cite{kock:2003a,fauser:2012a}.
\mybenv{Definition}\textbf{[Commutative Frobenius algebra]}
A commutative Frobenius algebra is given by a commutative unital
multiplication $\sfa$ with unit $\eta^\sfa$ (later also called $M$), and
a cocommutative coalgebra comultiplication $\delta_\sfa$ ($\delta_\sfa(x){}
=x_{[1]}\ot x_{[2]}$) with counit $\epsilon^\sfa$ (later also called $\epsilon^{1}$),
such that the Frobenius law holds
$\sfa(x,y_{[1]})\ot y_{[2]} = \sfa(x,y)_{[1]}\ot
\sfa(x,y)_{[2]} = x_{[1]}\ot \sfa(x_{[2]},y)$ or graphically
\begin{align}
   \begin{pic}
      \node (i1) at (0,1) {};
      \node (i2) at (1,1) {};
      \node[circle,fill=white,draw,inner sep=1pt] (a1) at (0.5,-0.25) {$\sfa$ };
      \node[circle,fill=white,draw,inner sep=1pt] (a2) at (1,0.25) {$\delta_\sfa$ };
      \node (o1) at (0.5,-1) {};
      \node (o2) at (1.5,-1) {};
\begin{pgfonlayer}{background}
      \draw[thick,out=270,in=180] (i1.center) to (a1.center);
      \draw[thick] (i2.center) to (a2.center);
      \draw[thick] (a1.center) to (o1.center);
      \draw[thick,out=200,in=20] (a2.center) to (a1.center);
      \draw[thick,out=0,in=90] (a2.center) to (o2.center);
\end{pgfonlayer}
\end{pic}
&\cong
   \begin{pic}
      \node (i1) at (0,1) {};
      \node (i2) at (1,1) {};
      \node[circle,fill=white,draw,inner sep=1pt] (a1) at (0.5,0.35) {$\sfa$ };
      \node[circle,fill=white,draw,inner sep=1pt] (a2) at (0.5,-0.35) {$\delta_\sfa$ };
      \node (o1) at (0,-1) {};
      \node (o2) at (1,-1) {};
\begin{pgfonlayer}{background}
      \draw[thick,out=270,in=180] (i1.center) to (a1.center);
      \draw[thick,out=270,in=0] (i2.center) to (a1.center);
      \draw[thick] (a1.center) to (a2.center);
      \draw[thick,out=180,in=90] (a2.center) to (o1.center);
      \draw[thick,out=0,in=90] (a2.center) to (o2.center);
\end{pgfonlayer}
\end{pic}
\cong
   \begin{pic}
      \node (i1) at (0.5,1) {};
      \node (i2) at (1.5,1) {};
      \node[circle,fill=white,draw,inner sep=1pt] (a1) at (0.5,0.25) {$\delta_\sfa$ };
      \node[circle,fill=white,draw,inner sep=1pt] (a2) at (1,-0.25) {$\sfa$ };
      \node (o1) at (0,-1) {};
      \node (o2) at (1,-1) {};
\begin{pgfonlayer}{background}
      \draw[thick] (i1.center) to (a1.center);
      \draw[thick,out=270,in=0] (i2.center) to (a2.center);
      \draw[thick,out=180,in=90] (a1.center) to (o1.center);
      \draw[thick,out=20,in=200] (a1.center) to (a2.center);
      \draw[thick] (a2.center) to (o2.center);
\end{pgfonlayer}
\end{pic}
\end{align}
In the graded case, we demand the Frobenius structure to be grade wise
defined, $\sfa : \H^i \ot \H^i \rightarrow \H^i$ and zero for mixed
grades. The unit $\eta^\sfa$ is then a disjoint union of units for each
grade, and hence constitutes a formal series $\eta^\sfa:=\sum_i \eta^{\sfa,i}$,
similarly the counit $\epsilon^\sfa=\sum_i \epsilon^{\sfa,i}$ is defined
grade wise too.
\myeenv

Note that the unit $\eta$ and counit $\epsilon$ of the ambient Hopf
algebra are quite different maps from the unit $\eta^\sfa$ and counit
$\epsilon^\sfa$ of the Frobenius multiplication and comultiplication.
A (finitely generated) Frobenius algebra allows a closed structure to
be defined, depicted as cup and cap tangles via
\begin{align}
   \begin{pic}
      \node (o1) at (0,0) {};
      \node (o2) at (1,0) {};
      \draw[thick,out=90,in=90] (o1.center) .. controls +(0,0.75) and +(0,0.75) .. (o2.center);
   \end{pic}
&\cong
   \begin{pic}
      \node[circle,fill=black,inner sep=2pt] (i1) at (0.5,0.75) {};
      \node[xshift=7pt] at (i1.east) {$M $};
      \node[circle,fill=white,draw,inner sep=1pt] (a2) at (0.5,0) {$\delta_\sfa$ };
      \node (o1) at (0,-0.5) {};
      \node (o2) at (1,-0.5) {};
\begin{pgfonlayer}{background}
      \draw[thick] (i1.center) to (a2.center);
      \draw[thick,out=180,in=90] (a2.center) to (o1.center);
      \draw[thick,out=0,in=90] (a2.center) to (o2.center);
\end{pgfonlayer}
\end{pic}
;&&&
   \begin{pic}
      \node (o1) at (0,0) {};
      \node (o2) at (1,0) {};
      \draw[thick,out=90,in=90] (o1.center) .. controls +(0,-0.75) and +(0,-0.75) .. (o2.center);
   \end{pic}
&\cong
   \begin{pic}
      \node[circle,fill=black,inner sep=2pt] (i1) at (0.5,-0.75) {};
      \node[xshift=7pt] at (i1.east) {$\epsilon^1$ };
      \node[circle,fill=white,draw,inner sep=1pt] (a2) at (0.5,0) {$\sfa$ };
      \node (o1) at (0,0.5) {};
      \node (o2) at (1,0.5) {};
\begin{pgfonlayer}{background}
      \draw[thick] (i1.center) to (a2.center);
      \draw[thick,out=180,in=270] (a2.center) to (o1.center);
      \draw[thick,out=0,in=270] (a2.center) to (o2.center);
\end{pgfonlayer}
\end{pic}
.&&&
\end{align}
Equivalently one can use the cup and cap closed structure to define the
Frobenius comultiplication $\delta_\sfa$~\cite{kock:2003a,fauser:2012a}
via
\begin{align}\label{grph-FrobeniusComultiplication}
\begin{pic}
      \node (i) at (0.6,1) {};
      \node[circle,inner sep=1pt,draw] (a) at (0.6,0) {$\delta_\sfa$ };
      \node (o1) at (0,-1) {};
      \node (o2) at (1.2,-1) {};
      \draw[thick,out=180,in=90] (a.west) to (o1.center);
      \draw[thick,out=  0,in=90] (a.east) to (o2.center);
      \draw[thick] (i.center) to (a.north);
\end{pic}
&:\cong
   \begin{pic}
      \node (u1) at (0,0.5) {};
      \node (u2) at (0.75,0.5) {};
      \node (i) at (1.75,1) {};
      \node[circle,inner sep=1pt,draw] (a) at (1.25,-0.2) {$\sfa$ };
      \node (o1) at (0,-1) {};
      \node (o2) at (1.25,-1) {};
      \draw[thick,out=270,in=180] (u2.center) to (a.west);
      \draw[thick,out=270,in=0] (i.center) to (a.east);
      \draw[thick,out=90,in=90] (u1.center) to (u2.center);
      \draw[thick] (u1.center) to (o1.center);
      \draw[thick] (a.south) to (o2.center);
   \end{pic}
\cong
   \begin{pic}
      \node (u1) at (1,0.5) {};
      \node (u2) at (1.75,0.5) {};
      \node (i) at (0,1) {};
      \node[circle,inner sep=1pt,draw] (a) at (0.5,-0.2) {$\sfa$ };
      \node (o1) at (1.75,-1) {};
      \node (o2) at (0.5,-1) {};
      \draw[thick,out=270,in=  0] (u1.center) to (a.east);
      \draw[thick,out=270,in=180] (i.center) to (a.west);
      \draw[thick,out=90,in=90] (u1.center) to (u2.center);
      \draw[thick] (u2.center) to (o1.center);
      \draw[thick] (a.south) to (o2.center);
   \end{pic}
&&&
\end{align}
A mixed bialgebra structure between the Hopf algebra multiplication $\mul$
and the Frobenius comultiplication $\delta_\sfa$ would read as
follows
\begin{align}\label{grph-mixedBialgebra}
   \begin{pic}
      \node (i1) at (0,0.75) {};
      \node (i2) at (1,0.75) {};
      \node (u) at (0.5,0.4) {};
      \node[circle,inner sep=1pt,draw,fill=white] (d) at (0.5,-0.4) {$\delta_\sfa$ };
      \node (o1) at (0,-0.75) {};
      \node (o2) at (1,-0.75) {};
\begin{pgfonlayer}{background}
      \draw[thick,out=270,in=180] (i1.center) to (u.center);
      \draw[thick,out=270,in=  0] (i2.center) to (u.center);
      \draw[thick] (u.center) to (d.center);
      \draw[thick,out=180,in=90] (d.center) to (o1.center);
      \draw[thick,out=  0,in=90] (d.center) to (o2.center);
\end{pgfonlayer}
   \end{pic}
\cong
   \begin{pic}
      \node (i1) at (0.25,0.75) {};
      \node (i2) at (1.25,0.75) {};
      \node (d1) at (0.25,-0.4) {};
      \node (d2) at (1.25,-0.4) {};
      \node (m1) at (-0.25,-0.05) {};
      \node (m4) at (1.75,-0.05) {};
      \node[circle,draw,inner sep=1pt,fill=white] (u1) at (0.25,0.3) {$\delta_\sfa$ };
      \node[circle,draw,inner sep=1pt,fill=white] (u2) at (1.25,0.3) {$\delta_\sfa$ };
      \node (o1) at (0.25,-0.75) {};
      \node (o2) at (1.25,-0.75) {};
\begin{pgfonlayer}{background}
      \draw[thick] (i1.center) to (u1.center);
      \draw[thick] (i2.center) to (u2.center);
      \draw[thick,out=180,in=90] (u1.center) to (m1.center);
      \draw[thick,out=  0,in=180] (u1.center) to (d2.center);
      \draw[thick,out=180,in=  0] (u2.center) to (d1.center);
      \draw[thick,out=  0,in=90] (u2.center) to (m4.center);
      \draw[thick,out=270,in=180] (m1.center) to (d1.center);
      \draw[thick,out=270,in=  0] (m4.center) to (d2.center);
      \draw[thick] (d1.center) to (o1.center);
      \draw[thick] (d2.center) to (o2.center);
\end{pgfonlayer}
   \end{pic}
&&\textrm{dually}&&
   \begin{pic}
      \node (i1) at (0,0.75) {};
      \node (i2) at (1,0.75) {};
      \node[circle,inner sep=1pt,draw,fill=white] (u) at (0.5,0.4) {$\sfa$ };
      \node (d) at (0.5,-0.4) {};
      \node (o1) at (0,-0.75) {};
      \node (o2) at (1,-0.75) {};
\begin{pgfonlayer}{background}
      \draw[thick,out=270,in=180] (i1.center) to (u.center);
      \draw[thick,out=270,in=  0] (i2.center) to (u.center);
      \draw[thick] (u.center) to (d.center);
      \draw[thick,out=180,in=90] (d.center) to (o1.center);
      \draw[thick,out=  0,in=90] (d.center) to (o2.center);
\end{pgfonlayer}
   \end{pic}
\cong
   \begin{pic}
      \node (i1) at (0.25,0.75) {};
      \node (i2) at (1.25,0.75) {};
      \node[circle,draw,inner sep=1pt,fill=white] (d1) at (0.25,-0.4) {$\sfa$ };
      \node[circle,draw,inner sep=1pt,fill=white] (d2) at (1.25,-0.4) {$\sfa$ };
      \node (m1) at (-0.25,0) {};
      \node (m4) at (1.75,0) {};
      \node (u1) at (0.25,0.4) {};
      \node (u2) at (1.25,0.4) {};
      \node (o1) at (0.25,-0.75) {};
      \node (o2) at (1.25,-0.75) {};
\begin{pgfonlayer}{background}
      \draw[thick] (i1.center) to (u1.center);
      \draw[thick] (i2.center) to (u2.center);
      \draw[thick,out=180,in=90] (u1.center) to (m1.center);
      \draw[thick,out=  0,in=180] (u1.center) to (d2.center);
      \draw[thick,out=180,in=  0] (u2.center) to (d1.center);
      \draw[thick,out=  0,in=90] (u2.center) to (m4.center);
      \draw[thick,out=270,in=180] (m1.center) to (d1.center);
      \draw[thick,out=270,in=  0] (m4.center) to (d2.center);
      \draw[thick] (d1.center) to (o1.center);
      \draw[thick] (d2.center) to (o2.center);
\end{pgfonlayer}
   \end{pic}
;
\end{align}
This allows us to make the following
\mybenv{Definition}\textbf{[Frobenius Laplace pairing]}
A Laplace pairing $\sfa$ is Frobenius, if $\sfa$ is a Laplace pairing
with respect to the ambient Hopf algebra $\HB$, and
$(\HB,\sfa,\delta_\sfa,M,\epsilon^1)$ is grade by grade a commutative
Frobenius algebra.
\myeenv
\mybenv{Proposition}\label{prop-FrobeniusBialgebra}
The Laplace pairing property for a Frobenius Laplace pairing is equivalent
to a \emph{mixed bialgebra structure} $(\HB,\mul,\delta_\sfa)$,
see~\eqref{grph-mixedBialgebra}, where $\delta_\sfa$ is the
Frobenius comultiplication~\eqref{grph-FrobeniusComultiplication}.
\myeenv

\noindent
\textbf{Proof:}
We provide a graphical proof as follows (using $\sfa=\delta_{\sfa}$ due
to the symmetry~\eqref{grph-FrobeniusComultiplication} and for
typographical convenience)
\begin{align}
   \begin{pic}
      \node (i1) at (0,0.5) {};
      \node (i2) at (0.5,0.5) {};
      \node (u) at (0.25,0.2) {};
      \node[circle,inner sep=1pt,draw] (d) at (0.25,-0.2) {$\sfa$ };
      \node (o1) at (0,-0.5) {};
      \node (o2) at (0.5,-0.5) {};
      \draw[thick,out=270,in=180] (i1.center) to (u.center);
      \draw[thick,out=270,in=  0] (i2.center) to (u.center);
      \draw[thick] (u.center) to (d);
      \draw[thick,out=180,in=90] (d) to (o1.center);
      \draw[thick,out=  0,in=90] (d) to (o2.center);
   \end{pic}
&\stackrel{Frob.}{\cong}
   \begin{pic}
      \node (i0) at (-0.5,0.5) {};
      \node (i1) at (0,0.5) {};
      \node (i2) at (0.5,0.75) {};
      \node (i3) at (1,0.75) {};
      \node (m) at (0.75,0.3) {};
      \node[circle,draw,inner sep=1pt] (a) at (0.375,-0.3) {$\sfa$ };
      \node (o0) at (-0.5,-0.75) {};
      \node (o) at (0.375,-0.75) {};
      \draw[thick,out=270,in=180] (i1.center) .. controls +(down:8mm) and +(left:3mm) .. (a.west);
      \draw[thick,out=270,in=180] (i2.center) to (m.center);
      \draw[thick,out=270,in=  0] (i3.center) to (m.center);
      \draw[thick,out=270,in=  0] (m.center) to (a.east);
      \draw[thick] (a.south) to (o.center);
      \draw[thick,out=90,in=90] (i0.center) to (i1.center);
      \draw[thick] (i0.center) to (o0.center);
   \end{pic}
\stackrel{Lap.}{\cong}
   \begin{pic}
      \node (i0) at (-0.5,0.75) {};
      \node (i1) at (0.25,0.75) {};
      \node (i2) at (0.75,0.75) {};
      \node (i3) at (1.25,0.75) {};
      \node (m) at (0.25,0.5) {};
      \node[circle,draw,inner sep=1pt] (a1) at (0.25,0) {$\sfa$ };
      \node[circle,draw,inner sep=1pt] (a2) at (1,0) {$\sfa$ };
      \node (d) at (0.6125,-0.5) {};
      \node (o) at (0.6125,-0.75) {};
      \draw[thick] (i1.center) to (m.center);
      \draw[thick,out=180,in= 90] (m.center) .. controls +(-0.3,0) and +(-0.3,0.2) .. (a1.west);
      \draw[thick,out=  0,in=180] (m.center) to (a2.west);
      \draw[thick,out=270,in=0] (i2.center) to (a1.east);
      \draw[thick,out=270,in=0] (i3.center) .. controls +(0,-0.4) and +(0.1,0) .. (a2.east);
      \draw[thick,out=270,in=180] (a1.south) to (d.center);
      \draw[thick,out=270,in=  0] (a2.south) to (d.center);
      \draw[thick] (d.center) to (o.center);
      \draw[thick,out=90,in=90] (i0.center) to (i1.center);
      \draw[thick] (i0.center) to (o0.center);
   \end{pic}
\stackrel{dual.}{\cong}
   \begin{pic}
      \node (u1) at (-0.25,0.75) {};
      \node (u2) at (-0.25,0.5) {};
      \node (i2) at (0.75,0.75) {};
      \node (i3) at (1.25,0.75) {};
      \node (m) at (-0.5,0) {};
      \node[circle,draw,inner sep=1pt,fill=white] (a1) at (0.25,0) {$\sfa$ };
      \node[circle,draw,inner sep=1pt,fill=white] (a2) at (1,0) {$\sfa$ };
      \node (d) at (0.6125,-0.5) {};
      \node (o1) at (-0.5,-0.75) {};
      \node (o2) at (0.6125,-0.75) {};
\begin{pgfonlayer}{background}
      \draw[thick,out=270,in=0] (i2.center) to (a1.center);
      \draw[thick,out=270,in=0] (i3.center) .. controls +(0,-0.4) and +(0.1,0) .. (a2.east);
      \draw[thick,out=0,in=180] (u1.center) to (a2.center);
      \draw[thick,out=0,in=180] (u2.center) to (a1.center);
      \draw[thick,out=180,in=180] (u1.center) .. controls +(-0.5,0) and +(-0.4,0) .. (m.center);
      \draw[thick,out=180,in=0] (u2.center) to (m.center);
      \draw[thick] (m.center) to (o1.center);
      \draw[thick,out=270,in=180] (a1.south) to (d.center);
      \draw[thick,out=270,in=  0] (a2.south) to (d.center);
      \draw[thick] (d.center) to (o2.center);
\end{pgfonlayer}
   \end{pic}
%%%%%%%%%%%%%%%%%%%%%%%%%%%%%%%%%%%%%%%%%%%%%%%%%%%%%%%%%%%%%%%%%%%%%%%%
\\
&\stackrel{Frob.}{\cong}
   \begin{pic}
      \node (i1) at (1.5,0.75) {};
      \node (i2) at (2.5,0.75) {};
      \node (u) at (0,0) {};
      \node[circle,draw,inner sep=1pt,fill=white] (a1) at (1,0) {$\sfa$ };
      \node[circle,draw,inner sep=1pt,fill=white] (a2) at (2,0) {$\sfa$ };
      \node (m1) at (0.5,-0.5) {};
      \node (m2) at (1.5,-0.5) {};
      \node (o1) at (0.5,-0.75) {};
      \node (o2) at (1.5,-0.75) {};
      \node (spacing-dummy) at (0,-1.5) {};
\begin{pgfonlayer}{background}
      \draw[thick,out=270,in=90] (i1.center) to (a1.center);
      \draw[thick,out=270,in=0] (i2.center) to (a2.center);
      \draw[thick,out=90,in=180] (u.center) .. controls +(0,1.75) and +(-0.3,-0.2) .. (a2.center);
      \draw[thick,out=270,in=180] (u.center) to (m1.center);
      \draw[thick,out=180,in=0] (a1.center) to (m1.center);
      \draw[thick] (m1.center) to (o1.center);
      \draw[thick,out=0,in=180] (a1.center) to (m2.center);
      \draw[thick,out=270,in=0] (a2.center) to (m2.center);
      \draw[thick] (m2.center) to (o2.center);
\end{pgfonlayer}
   \end{pic}
%%%
\stackrel{Frob.+com.}{\cong}
   \begin{pic}
      \node (i1) at (0.25,0.75) {};
      \node (i2) at (1.25,0.75) {};
      \node (d1) at (0.25,-0.4) {};
      \node (d2) at (1.25,-0.4) {};
      \node (m1) at (-0.25,0) {};
      \node (m4) at (1.75,0) {};
      \node[circle,draw,inner sep=1pt,fill=white] (u1) at (0.25,0.3) {$\sfa$ };
      \node[circle,draw,inner sep=1pt,fill=white] (u2) at (1.25,0.3) {$\sfa$ };
\begin{pgfonlayer}{background}
      \node (o1) at (0.25,-0.75) {};
      \node (o2) at (1.25,-0.75) {};
      \draw[thick] (i1.center) to (u1);
      \draw[thick] (i2.center) to (u2);
      \draw[thick,out=180,in=90] (u1.center) to (m1.center);
      \draw[thick,out=  0,in=180] (u1.center) to (d2.center);
      \draw[thick,out=180,in=  0] (u2.center) to (d1.center);
      \draw[thick,out=  0,in=90] (u2.center) to (m4.center);
      \draw[thick,out=270,in=180] (m1.center) to (d1.center);
      \draw[thick,out=270,in=  0] (m4.center) to (d2.center);
      \draw[thick] (d1.center) to (o1.center);
      \draw[thick] (d2.center) to (o2.center);
\end{pgfonlayer}
   \end{pic}
\end{align}
\vskip-6ex
\qed

\noindent
This result relates a Frobenius Laplace pairing with the closed
structure used to dualize the Hopf algebra $\HB$, hence with a duality
pairing $\la -\mid-\ra : \HBI{*}\times \HB \rightarrow \Bbbk$.
By general Hopf algebra theory the multiplication and comultiplication
on $\HBI{*}$ are then given by $\comul^*$ and $\mul^*$. For the relation
between Hopf and Frobenius algebras see~\cite{fauser:2012a} and
references therein.
%-----------------------------------------------------------------------

The right hand side of~\eqref{grph-mixedBialgebra} depicts condition (e) of
Rota and Stein~\cite{rota:stein:1994a}, further discussed in
Sections~\ref{sec-CircleProducts} and~\ref{sec-HashProd}. 
Proposition~\ref{prop-FrobeniusBialgebra} shows then that the condition
(e) of Rota and Stein demands that the Laplace pairing in use is Frobenius.
%-----------------------------------------------------------------------

From Proposition~\ref{prop-FrobeniusBialgebra} it follows that
\mybenv{Corollary}
Let $\sfa,\sfb$ be Frobenius Laplace pairings, then $\sfc=\sfa\conv\sfb$
is also Frobenius Laplace.
\myeenv

\mybenv{Corollary}
Let $(\sfa,\Delta)$ be a mixed bialgebra with $\sfa$ Frobenius and let
$\phi\in\hom(\HB,\HB)$ be such that $\comul\circ\phi = (\phi\ot\phi)\circ\comul${}
is a \emph{coalgebra} morphism, then $(\sfa_{\phi},\comul)$ forms a mixed
bialgebra~\eqref{grph-mixedBialgebra}.
\myeenv
The derived pairing $\sfa_{\phi}$ yields via dualization no longer a
\emph{commutative} Frobenius comultiplication. We obtain a left
dual $\delta_{\sfa_{\phi}}^{l}=(1\otimes\phi)\circ\delta_{\sfa}$ or a
right dual $\delta_{\sfa_{\phi}}^{r}=(\phi\otimes 1)\circ\delta_{\sfa}$
and either of ($\sfa_{\phi},\delta_{\sfa_{\phi}}^{l/r})$ is in general
not commutative Frobenius.

We show finally two convolutive inverses for the associative Hopf algebra
multiplication $\mul$ and a Laplace Frobenius pairing $\sfa$. As these
maps are unital but not necessarily normalized we could use a generalized
Milnor--Moore inverse. However the recursive inverse suffers from cancellation.
We can use the Hopf algebra antipode in these two cases to give a direct
inverse. Recall that we work in a bicommutative setting and it is immaterial
that the antipode maps $\mul$ to $\mul^{op}$, the opposite product, as
$\mul^{op}= \mul\circ\sw =\mul$.
\mybenv{Proposition}\label{prop-twoInverses}
Let $\overline{\mul} = \antip\circ\mul$ and $\sfoa = \antip\circ\sfa$ 
where $\mul$ is a Hopf algebra mutiplication and $\sfa$ is a Frobenius Laplace pairing.
Then $\overline{\mul}$ and $\sfoa$ are the convolutive inverses of $\mul$ and $\sfa$,
respectively. Moreover $\overline{\mul}=\mul\circ(\antip\ot\antip)$.
\myeenv

\noindent
\textbf{Proof:} We give the graphical proof for $\sfoa\conv\sfa=\sfe$,
the other case $\overline{\mul}\conv\mul=\sfe$ is obtained merely by replacing 
$\sfa$ by $\mul$ throughout.
\begin{align}
   \begin{pic}
      \node (i1) at (0.25,0.75) {};
      \node (i2) at (1,0.75) {};
      \node (m1) at (0.25,0.5) {};
      \node (m2) at (1,0.5) {};
      \node[circle,draw,inner sep=1pt,fill=white] (a1) at (0.25,0) {$\sfoa$ };
      \node[circle,draw,inner sep=1pt,fill=white] (a2) at (1,0) {$\sfa$ };
      \node (d) at (0.6125,-0.5) {};
      \node (o) at (0.6125,-0.75) {};
\begin{pgfonlayer}{background}
      \draw[thick] (i1.center) to (m1.center);
      \draw[thick] (i2.center) to (m2.center);
      \draw[thick,out=180,in=180] (m1.center) to (a1.center);
      \draw[thick,out=  0,in=180] (m1.center) to (a2.center);
      \draw[thick,out=180,in=0] (m2.center) to (a1.center);
      \draw[thick,out=  0,in=0] (m2.center) to (a2.center);
      \draw[thick,out=270,in=180] (a1.center) to (d.center);
      \draw[thick,out=270,in=  0] (a2.center) to (d.center);
      \draw[thick] (d.center) to (o.center);
\end{pgfonlayer}
   \end{pic}
\stackrel{Ansatz}{\cong}
   \begin{pic}
      \node (i1) at (0.25,1) {};
      \node (i2) at (1,1) {};
      \node (m1) at (0.25,0.75) {};
      \node (m2) at (1,0.75) {};
      \node[circle,draw,inner sep=1pt,fill=white] (a1) at (0.25,0.25) {$\sfa$ };
      \node[circle,draw,inner sep=1pt,fill=white] (a2) at (1,0.25) {$\sfa$ };
      \node[circle,draw,inner sep=0pt,fill=white] (s) at (0.25,-0.25) {$\antip$ };
      \node (d) at (0.6125,-0.75) {};
      \node (o) at (0.6125,-1) {};
\begin{pgfonlayer}{background}
      \draw[thick] (i1.center) to (m1.center);
      \draw[thick] (i2.center) to (m2.center);
      \draw[thick,out=180,in=180] (m1.center) to (a1.center);
      \draw[thick,out=  0,in=180] (m1.center) to (a2.center);
      \draw[thick,out=180,in=0] (m2.center) to (a1.center);
      \draw[thick,out=  0,in=0] (m2.center) to (a2.center);
      \draw[thick] (a1.center) to (s.center);
      \draw[thick,out=270,in=180] (s.center) to (d.center);
      \draw[thick,out=270,in=  0] (a2.center) to (d.center);
      \draw[thick] (d.center) to (o.center);
\end{pgfonlayer}
   \end{pic}
\stackrel{bialg.}{\cong}
   \begin{pic}
      \node (i1) at (0,1) {};
      \node (i2) at (1,1) {};
      \node[circle,draw,inner sep=1pt,fill=white] (a) at (0.5,0.6) {$\sfa$ };
      \node (u) at (0.5,0.2) {};
      \node[circle,draw,inner sep=0pt,fill=white] (s) at (0.25,-0.2) {$\antip$ };
      \node (d) at (0.5,-0.6) {};
      \node (o) at (0.5,-1) {};
\begin{pgfonlayer}{background}
      \draw[thick,out=270,in=180] (i1.center) to (a.center);
      \draw[thick,out=270,in=0] (i2.center) to (a.center);
      \draw[thick] (a.center) to (u.center);
      \draw[thick,out=180,in=90] (u.center) to (s.center);
      \draw[thick,out=270,in=180] (s.center) to (d.center);
      \draw[thick,out=0,in=0] (u.center) to (d.center);
      \draw[thick] (d.center) to (o.center);
\end{pgfonlayer}
\end{pic}
\stackrel{antip.}{\cong}
   \begin{pic}
      \node (i1) at (0,1) {};
      \node (i2) at (1,1) {};
      \node[circle,draw,inner sep=1pt,fill=white] (a) at (0.5,0.6) {$\sfa$ };
      \node[circle,inner sep=2pt,draw,thick] (u) at (0.5,0) {};
      \node[xshift=3pt] at (u.east) {$\epsilon$};
      \node[circle,inner sep=2pt,draw,thick] (d) at (0.5,-0.3) {};
      \node[xshift=3pt] at (d.east) {$\eta$};
      \node (o) at (0.5,-1) {};
\begin{pgfonlayer}{background}
      \draw[thick,out=270,in=180] (i1.center) to (a.center);
      \draw[thick,out=270,in=0] (i2.center) to (a.center);
      \draw[thick] (a.center) to (u.north);
      \draw[thick] (d.south) to (o.center);
\end{pgfonlayer}
\end{pic}
\stackrel{counitality}{\cong}
   \begin{pic}
      \node (i1) at (0,1) {};
      \node (i2) at (0.6,1) {};
      \node[circle,inner sep=2pt,draw,thick] (eps1) at (0,0.3) {};
      \node[xshift=3pt] at (eps1.east) {$\epsilon$};
      \node[circle,inner sep=2pt,draw,thick] (eps2) at (0.6,0.3) {};
      \node[xshift=3pt] at (eps2.east) {$\epsilon$};
      \node[circle,inner sep=2pt,draw,thick] (eta) at (0.3,-0.3) {};
      \node[xshift=3pt] at (eta.east) {$\eta$};
      \node (o) at (0.3,-1) {};
      \draw[thick] (i1.center) to (eps1);
      \draw[thick] (i2.center) to (eps2);
      \draw[thick] (eta) to (o.center);
   \end{pic}
\end{align}
%\vskip-2ex
\noindent
To obtain the final result, already noted at the end of section 2.2, 
it is only necessary to note that $\antip\circ\mul = \mul\circ(\antip\ot\antip)$.
\qed

Hence for Frobenius Laplace pairings we find a direct way to produce
the inverse using the antipode, and do not need the generalized
recursive formula.

We close this section showing that the subgroup of Frobenius Laplace
pairings is nontrivial. An element $p$ of $\HB$ is (1-1) primitive if the
comultiplication reads $\comul(p)=p\ot 1 + 1\ot p$. For $\textrm{char}(\Bbbk)=0$
using the Cartier-Milnor-Moore theorem we can construct a
Poincar\'{e}--Birkhoff--Witt basis in terms of primitive elements. If we
define $P^{1}:=\{x\in\HB \mid \epsilon(x)=0,~\comul(x)=x\ot 1 + 1\ot x\}$
we can decompose $\HB=\Bbbk \oplus (\oplus_{k\geq 1} P^{k})$ with
$P^{k}=\ot^{k}P^{1}$. Elements in $P^{1}$ are polynomial generators.
Rota and Stein showed that a Laplace pairing can be defined by its action on
primitive elements.
\mybenv{Theorem}
Let $\HB$ be a graded connected bicommutative Hopf algebra with a
Poincar\'{e}--Birkhoff--Witt basis of primitive elements
$\H^{+}=\oplus_{k\geq 1} P^{k}$, then a Laplace pairing $\sfa$ is given by the
following data ($x,y,z$ generators of $P^{1}$, $w=\prod_{i} x_{i},
w^{\prime}=\prod_{j} x^{\prime}_{j}$):
\begin{itemize}
\item[i)] $\sfa : \H^{0} \ot\H^{0} \rightarrow \H^{0} 
   :: \sfa(1,1)=\phi_{1,1}^{1}1=1$ normalization.
\item[ii)] $\sfa : P^{1} \ot P^{1} \rightarrow P^{1} 
   :: \sfa(x,y) = \sum_{z}\phi_{x,y}^{z} z$ on primitive elements $x,y,z\in P^{1}$
   ($\phi$ can be zero, producing the trivial pairing).
\item[iii)] $\sfa : P^{i} \ot P^{j} \rightarrow 0 
   :: \sfa(w,w^{\prime}) = 0$ for $i\not= j$, grading.
\item[iv)] $\sfa : P^{n} \ot P^{n} \rightarrow P^{n} 
   :: \sfa(w,w^{\prime}) = \sfa(\prod_{i} x_{i}, \prod_{j} x^{\prime}_{j}) 
                = \prod_{i,j} \sfa(x_{i(j)}, x^{\prime}_{j(i)})$ for $n\ge 1$ by the
   Laplace expansions~(\ref{eq-LaplaceRight},\ref{eq-LaplaceLeft}).        
\end{itemize}
We may choose $\phi_{x,y}^{z}=\delta_{x,y}\delta_{x,z} z\,f(z)$ with 
$f : P^{1}\rightarrow \Bbbk$ a weight function.
\myeenv
This theorem does not produce the most general Laplace pairings. In the above
form $\sfa$ is Frobenius. But one may easily construct derived Laplace 
pairings by post-composing with a 1-cocycle.

%-----------------------------------------------------------------------
\subsection{Hopf algebra deformations (circle products)}\label{sec-CircleProducts}
%-----------------------------------------------------------------------
Rota and Stein gave in~\cite{rota:stein:1994a} a deformation theory
and a list of examples for supersymmetric letter place algebras showing
how to deform the multiplication in a graded connected super-commutative
Hopf algebra. One could equivalently deform the coproduct. Using letter
place algebras allows in principle a development over any base ring,
including the posibility of additive torsion. Rota and Stein's main
tool is the so-called \emph{circle product}, also called \emph{cliffordization},
which is a twist in terms of Hopf algebra theory. In~\cite{rota:stein:1994a}
the Laplace pairings are given by four conditions (a), (b) and (c), (d).
The conditions (a), (b) are the Laplace expansion laws~\eqref{eq-LaplaceLeft}
and~\eqref{eq-LaplaceRight}. We have seen above that Rota and Stein's
fifth condition (e) is equivalent to a mixed bialgebra
law~\eqref{grph-mixedBialgebra} and hence by
Proposition~\ref{prop-FrobeniusBialgebra} to a Frobenius Laplace pairing.
Conditions (c) and (d) (which we have not displayed) are implied by
condition (e), however if the Laplace pairing is not Frobenius they may
fail to hold (check the Grassmann-Clifford example
in~\cite{rota:stein:1994a}). We do not employ the conditions (c) and (d)
here, and note that Rota and Stein did not use them for deriving their
results either.

\mybenv{Definition}\textbf{[circle product]}
Let $\HB$ be a graded commutative Hopf algebra and $\sfc$ an unital
normalized Laplace pairing which is then also a 2-cocycle $\partial\sfc=\sfe$,
the \emph{circle product} is defined as the convolution
$\tikz{\node[circle,draw,inner sep=2pt] at (0,0) {};}=\sfc\conv\mul$
or graphically     
\begin{align}
   \begin{pic}
      \node (i1) at (0,1) {};
      \node (i2) at (1,1) {};
      \node[circle,draw,inner sep=2pt,fill=white] (a) at (0.5,0) {};
      \node (o) at (0.5,-1) {};
\begin{pgfonlayer}{background}
      \draw[thick,out=270,in=180] (i1.center) to (a.center);
      \draw[thick,out=270,in=0] (i2.center) to (a.center);
      \draw[thick] (a.center) to (o.center);
\end{pgfonlayer}
\end{pic}
\cong
   \begin{pic}
      \node (i1) at (0.25,1) {};
      \node (i2) at (1,1) {};
      \node (m1) at (0.25,0.75) {};
      \node (m2) at (1,0.75) {};
      \node[rectangle,draw,inner sep=3pt,fill=white] (a1) at (0.25,0) {$\sfc$ };
      \node[circle,draw,inner sep=1pt,fill=white] (a2) at (1,0) {$\mul$ };
      \node (d) at (0.6125,-0.75) {};
      \node (o) at (0.6125,-1) {};
\begin{pgfonlayer}{background}
      \draw[thick] (i1.center) to (m1.center);
      \draw[thick] (i2.center) to (m2.center);
      \draw[thick,out=180,in=90] (m1.center) to (a1.north west);
      \draw[thick,out=  0,in=180] (m1.center) to (a2.center);
      \draw[thick,out=180,in=90] (m2.center) to (a1.north east);
      \draw[thick,out=  0,in=0] (m2.center) to (a2.center);
      \draw[thick,out=270,in=180] (a1.center) to (d.center);
      \draw[thick,out=270,in=  0] (a2.center) to (d.center);
      \draw[thick] (d.center) to (o.center);
\end{pgfonlayer}
   \end{pic}
\cong
   \begin{pic}
      \node (i1) at (0,1) {};
      \node (i2) at (1.5,1) {};
      \node[rectangle,draw,fill=white] (a) at (0.75,0) {$\sfc\conv\mul$ };
      \node (o) at (0.75,-1) {};
\begin{pgfonlayer}{background}
      \draw[thick,out=270,in=135] (i1.center) to (a.north west);
      \draw[thick,out=270,in=45] (i2.center) to (a.north east);
      \draw[thick] (a.center) to (o.center);
\end{pgfonlayer}
\end{pic}
\end{align}
\vskip-3ex
\myeenv

\mybenv{Proposition}
The circle product is associative, unital with unit $\eta$.
\myeenv

\mybenv{Theorem}
If $\sfc$ is Frobenius (fulfils the Rota-Stein condition (e)), then
the circle product turns
$(\HB,\tikz{\node[circle,draw,inner sep=2pt] at (0,0) {};},\eta,\Delta,\epsilon)$
into graded connected Hopf algebra.
\myeenv
Among other examples, Rota and Stein show how a Grassmann algebra can
be deformed into a Clifford algebra (not Hopf), and how from a cofree
cogenerated Hopf algebra $\H_\sqcup$ one can derive the symmetric
function Hopf algebra in one generator (alphabet). In what follows we
want to show that the subgroup of Laplace pairings in the monoid of
2-cocycles allows the parameterization of such products, Hopf or not,
in a very neat and efficient way. To do so we are going to generalize
the Rota-Stein circle product to allow derived Laplace pairings and multiple
convolutions of them and call such products higher derived hash products.
%-----------------------------------------------------------------------
\section{Hash products}\label{sec-HashProd}
\subsection{Higher derived hash products}
%-----------------------------------------------------------------------
In this section we extend the Rota-Stein deformation process by exploiting
the structure of the Laplace, or Frobenius Laplace subgroups of derived
Laplace pairings as a means of parameterizing the 2-cocycles used in the
product deformation. Choosing one or more elements of these subgroups
as generators, we obtain new associative products, which we term hash products.
In the Frobenius case these still form Hopf algebras with the
\emph{undeformed} coproduct; in the general cases the coproduct needs
to be changed to satisfy the bialgebra law~\eqref{grph-bialgebra-and-antipode}.
These higher derived hash products will then be discussed in the
application Section~\ref{sec-Applications}, where we study product
decompositions of characters or restricted subgroups of the general
linear group. The ambient Hopf algebra there is that of symmetric functions,
and the Frobenius Laplace pairing is inner multiplication, which serves
in derived multiplications as a generator for the deformations.

%-----------------------------------------------------------------------
In the previous section we studied the monoid of $k$-cochains and
some of its subgroups. For $k=2$ we call them pairings, and adding
more structure we get the inclusions:
\begin{itemize}
\item[\relax] Monoid of pairings $\mathsf{P}$
\item[$\supset$] subgroup $\mathsf{IP}$ of invertible pairings
\item[$\supset$] subgroup $\mathsf{C}$ of $2$-cocycles
\item[$\supset$] subgroup $\mathsf{L}$ of Laplace pairings, which we may restrict
     to unital normalized pairings.
\item[$\supset$] subgroup $\mathsf{F}$ of Frobenius Laplace pairings, which also
     may be restricted to the unital normalized case.
\end{itemize}

It is not clear if the Frobenius Laplace subgroup is maximal in the
Laplace subgroup of pairings or how it lies inside the Laplace subgroup.
However note that given a Frobenius Laplace pairing $\sfa$, the derived
pairing $\sfa_\phi=\phi\circ\sfa$ is in general only Laplace.

%-----------------------------------------------------------------------
Recall that we have already introduced \emph{derived Laplace pairings}
$\sfa_\phi := \phi\circ\sfa$ with $\sfa$ a (Frobenius) Laplace pairing
and $\phi$ a 1-cocycle and used their convolutively generated subgroup
of pairings to produce new deformations. This motivates the following

\mybenv{Definition}\textbf{[Higher derived hash product]}
Let $\sfa_{i,\phi_i}$, $1\leq i\leq k$ be derived Laplace pairings and
$\mul_{\phi_{k+1}}=\phi_{k+1}\circ\mul$ a 
\emph{derived} multiplication
of the ambient Hopf algebra $\HB$, for 1-cocycles $\phi_i$, $1\leq k \leq k+1$.
We define the \emph{higher derived hash product}
$\#_{\sfa_{1,\phi_1},\ldots,\sfa_{k,\phi_k},\phi_{k+1}}$ to be the
convolution product $\sfc = \sfa_{1,\phi_1}\conv\ldots\conv\sfa_{k,\phi_{k}}
\conv(\phi_{k+1}\circ \mul)$
\begin{align}
   \begin{pic}
      \node (i1) at (0,1) {};
      \node (i2) at (1,1) {};
      \node[circle,draw,inner sep=0pt,fill=white] (a) at (0.5,0) {$\#_{\ldots}$ };
      \node (o) at (0.5,-1) {};
\begin{pgfonlayer}{background}
      \draw[thick,out=270,in=180] (i1.center) to (a.center);
      \draw[thick,out=270,in=0] (i2.center) to (a.center);
      \draw[thick] (a.center) to (o.center);
\end{pgfonlayer}
\end{pic}
\cong
   \begin{pic}
      \node (i1) at (0.25,1.25) {};
      \node (i2) at (2.25,1.25) {};
      \node (m1) at (0.25,1) {};
      \node (m2) at (2.25,1) {};
      \node[rectangle,draw,inner sep=3pt,fill=white] (a1) at (0.25,0)
          {$\sfa_{1,\phi_1}\conv\ldots\conv\sfa_{k,\phi_{k}}$ };
      \node[circle,draw,inner sep=1pt,fill=white] (a2) at (2.5,0.3) {$\mul$ };
      \node[rectangle,draw,fill=white] (phi) at (2.5,-0.3) {$\phi_{k+1}$ };
      \node (d) at (1.5,-1) {};
      \node (o) at (1.5,-1.25) {};
\begin{pgfonlayer}{background}
      \draw[thick] (i1.center) to (m1.center);
      \draw[thick] (i2.center) to (m2.center);
      \draw[thick,out=180,in=90] (m1.center) to (a1.north west);
      \draw[thick,out=  0,in=135] (m1.center) to (a2.center);
      \draw[thick,out=180,in=90] (m2.center) to (a1.north east);
      \draw[thick,out=  0,in=0] (m2.center) to (a2.center);
      \draw[thick,out=270,in=180] (a1.center) to (d.center);
      \draw[thick] (a2.center) to (phi.north);
      \draw[thick,out=270,in=  0] (phi.south) to (d.center);
      \draw[thick] (d.center) to (o.center);
\end{pgfonlayer}
   \end{pic}
\end{align}
\myeenv

We usually simplify the notation for hash products to $\#_{\mathcal{L}}$
or more specifically to $\#_{\phi_1,\ldots,\phi_{k+1}}$ exhibiting either
a set $\mathcal{L}$ of derived Laplace pairings, or just the set of
$1$-cocycles of the derived pairings, if the Laplace pairings are
understood from the context. The identity morphism is denoted as $1$.
\mybenv{Proposition}
The higher derived hash product is associative.
\myeenv
This is immediate, as we constructed the subgroup of Laplace pairings to
be 2-cocycles, and deformations by a 2-cocycle are associative by general
Hopf algebra theory.

\mybenv{Theorem}\label{thm-HashAsHopf}
The coalgebra $(\HB,\comul,\epsilon)$ together with the derived higher
hash product and its unit $\eta$
$(\HB,\#_{\sfa_{1,\phi_1},\ldots,\sfa_{k,\phi_k},\phi_{k+1}},\eta,\Delta,\epsilon)$
is a commutative connected Hopf algebra, if and only if the derived higher
pairing $\sfc=\prod_{\conv,i}\sfa_{i,\phi_i}$ is Frobenius (fulfils
condition (e)~\eqref{grph-mixedBialgebra}).
\myeenv

\noindent
\textbf{Proof:}
If the derived higher pairing $\sfc=\prod_{\conv,i}\sfa_{i,\phi_i}$
is a Frobenius Laplace pairing, then we can adopt the Rota-Stein argument
or prove the result directly using the mixed bialgebra
property~\eqref{grph-mixedBialgebra}, which is implied by the Frobenius
condition. The same proof shows that if the mixed bialgebra property is
not available, then the pair ($\#_{\mathcal{L}},\Delta$) cannot be a 
bialgebra. In consequence it cannot be a Hopf algebra either. In the
Hopf algebra case the antipode is given by the Milnor--Moore inverse of
$\Id$ or by the Schmitt formula~\cite{schmitt:1987a}.
\qed

%-----------------------------------------------------------------------
\subsection{Hash products, special cases, and Heisenberg product}
%-----------------------------------------------------------------------
We will have occasion to deal with special cases of our higher derived
hash products $\#_{a_{1,\phi_1},\ldots,\sfa_{k,\phi_k},\phi_{k+1}}$
and find it convenient to name them.
\begin{itemize}
\item If in a higher derived hash product all $\phi_i$ are identity
      morphisms $\Id$, then we call it a \emph{higher hash product}.
      The algebraic closure property of Proposition~\ref{prop-Closure}
      shows that Frobenius Laplace pairings parameterize Hopf algebra
      isomorphisms $(\mul,\Delta)\rightarrow (\#_{\Id,\mathcal{L}},\Delta)$
      where only the product map is deformed.
\item A higher derived hash product involving only a single derived
      Laplace pairing is called a \emph{derived hash product} and takes
      the form $\#_{\phi_1,\phi_2} = (\phi_1\circ \sfa)\conv (\phi_2\circ\mul)$.
      A derived hash product can interpolate between the pairing $\sfa$ and 
      the multiplication $\mul$ through suitable choices of $\phi_1$ and $\phi_2$. 
      Restricting to minimal or maximal grades, one obtains either the pairing $\sfa$ 
      or the multiplication $\mul$.
\item A higher derived hash product for a single (Frobenius) Laplace
      pairing is simply called a \emph{hash product} and is of the form
      $\#_{\Id,\Id}\sfc=\sfa\conv\mul$, that is it is 
      a Rota-Stein circle product. In the case that $\sfc$ is Frobenius
      the circle product and the original comultiplication of the
      ambient Hopf algebra form a Hopf algebra.
\end{itemize}

%--NSYM-----------------------------------------------------------------
Looking at noncommutative symmetric functions forming the noncommutative
Hopf algebra $\NSym$, one can use the multiplication of the
Solomon descent algebra as a Laplace pairing. This situation was studied by
Aguiar et al.~\cite{aguiar:ferrer:moreira:2004a,aguiar:ferrer:moreira:2004b},
where the equivalent of our present (commutative) hash product was termed a
\emph{Heisenberg product} in the noncommutative situation. This product
is interpolating in the sense that projecting on lowest or highest grades
produces either the symmetric group or Solomon descent product (or 
outer and inner product in our case). This shows that hash products can
be generalized to a noncommutative setting, and also opens the way to
the study of higher derived noncommutative hash products.
In previous work~\cite{fauser:jarvis:king:wybourne:2005a} we showed how
plethystic branchings could be used to compute character decompositions
of non classical groups. The present method does not use plethystic
techniques for the product deformation. In a noncommutative setting
plethysms are a cumbersome operation and as of now not well understood,
but see~\cite{malvenuto:reutenauer:1998a}. However, in this work we stay
commutative as our target applications are in the field of group
characters and (ordinary) symmetric functions.

We will see that various specialisations of the higher derived hash
products are needed to describe the applications for group characters
to be discussed in the next Section. In the course of examination it
will also become clear why certain deformations are, or are not, Hopf
algebra deformations if the coproduct remains unchanged (as is the case
in the Rota--Stein setting).
%-----------------------------------------------------------------------
\subsection{Hopf algebra morphisms}\label{sec-HashProdComul}

Rota and Stein~\cite{rota:stein:1994a,rota:stein:1994b} gave a deformation
theory for Hopf algebras by deforming the product only. That is they
constructed a cliffordization process mapping structure maps $(\mul,\comul)${}
to $(\tikz{\node[circle,draw,inner sep=2pt] at (0,0) {};},\comul)$
where the circle product $\tikz{\node[circle,draw,inner sep=2pt] at (0,0) {};}$
is given by a convolution with a Laplace pairing with respect to the
underlying convolution algebra based on the original Hopf algebra maps
$(\mul,\comul)$. From the theory of Hopf algebra deformations and the
fact that a Laplace pairing is always a 2-cocycle,
Proposition~\ref{prop-laplaceImpliesCoycle}, it is clear that the
circle product is always associative, but one cannot hope in general
that the new pair of structure maps is again a bialgebra, much less Hopf.
We have shown, that the condition (e) of Rota-Stein is equivalent to
the fact that the Laplace pairing actually defines a (super)
commutative Frobenius algebra, Proposition~\ref{prop-FrobeniusBialgebra}.
If the pairing is not Frobenius there seem to be at least two
ways to overcome this problem.

In \cite{rota:stein:1994a,rota:stein:1994b} Rota and Stein examined
certain examples of circle products. The motivating example was to deform
a (supercommutative) Grassmann Hopf algebra $\bigwedge(V)$ over a space
$V$ into a Clifford algebra $C\!\ell(V,Q)$. As the pairings (symmetric bilinear
forms $B$) given by the polarization of a quadratic form $Q$ are in
general not Frobenius, the Clifford algebras are no longer Hopf
algebras, but see~\cite{fauser:oziewicz:2001a}. Then they extend the
ambient Hopf algebra to be $\H=\mathrm{Hilb}[V]\otimes\mathrm{Super}[V]$
with an extended Laplace pairing. For $w,w^{\prime}\in\mathrm{Hilb}[V]$
and $v,v^{\prime}\in\mathrm{Super}[V]$~\cite{rota:stein:1990b} the
Laplace pairing and circle product are given as
\begin{align}
  \sfa = \la w\ot v\mid w^{\prime}\ot v^{\prime}\ra
   &= \pm \epsilon(ww^{\prime}) \la v\mid v^{\prime}\ra
      \ot 1
\nn
   ( w\ot v)\, \tikz{\node[circle,draw,inner sep=2pt] at (0,0) {};} \,
   ( w^{\prime}\ot v^{\prime})
   &= \pm \epsilon(ww^{\prime}) \la v_{(1)}\mid v^{\prime}_{(1)}\ra
      \ot v_{(2)}v^{\prime}_{(2)}
\end{align}
where $\epsilon$ is a $\mathrm{Hilb}[V]$ valued form on $\mathrm{Hilb}[V]$
and the pairing $\la v_{(1)} \mid v^{\prime}_{(1)}\ra$ is also $\mathrm{Hilb}[V]$
valued, hence a map $\la -\mid -\ra : \mathrm{Super}[V]\ot \mathrm{Super}[V]
\rightarrow \mathrm{Hilb}[V]$. This provides a \emph{free Clifford algebra}
and the deformation by the circle product remains a Hopf algebra
for the original comultiplication. This process substantially enlarges
the underlying ambient Hopf algebra.

On the other hand, when we deform group characters to subgroup characters
in Section~\ref{sec-Applications} we know that we are still dealing with
a Hopf algebra of characters. In fact in the stable limit we are dealing
with the same Hopf algebra $\HB=\Sym$ of symmetric functions. It is hence
undesirable to lose the Hopf algebra structure in the deformation process.
We have developed, using plethystic
branchings~\cite{fauser:jarvis:2003a,fauser:jarvis:king:wybourne:2005a,%
fauser:jarvis:king:2005a,fauser:jarvis:2006a,fauser:jarvis:king:2007c,%
fauser:jarvis:king:2010a,fauser:jarvis:king:2010c},
a deformation theory of universal characters which is based on
isomorphisms between the underlying modules. For the product we use exactly
the same deformation theory as developed in Section~\ref{sec-ComHopf},
but we also need to deform the coproducts. A case study for orthogonal
and symplectic groups is given in~\cite{fauser:jarvis:king:2007c}.
Using some notation explained below in Section~\ref{sec-Applications},
we assume that we have a module map on $\HB$ given by 
$\Phi : \HB\rightarrow\HB :: x\mapsto x_{(1)} \la M_\pi\mid x_{(2)}\ra$,
with $M_\pi=M_\pi(1)$ where $M_{\pi}(t)\in\HB\lsqb t\rsqb$ is an infinite
series of symmetric functions defined later. The inverse map involves
$L_\pi$, the series inverse to $M_\pi$. Then the coproduct is deformed as
\begin{align}
   \comul_{\pi}(x) 
     &= x_{(1)}\ot x_{(2)} \la M_{\pi}\mid x_{(3)}\ra
      = \la M_{\pi}\mid x_{(1)}\ra x_{(2)}\ot x_{(3)} 
\end{align}
And regardless of the type of the Laplace pairing in use, the map $\Phi$
induced by $M_{\pi}$ is a Hopf algebra morphism for the structure maps
$\Phi : (\mul,\comul) \mapsto
(\, \tikz{\node[circle,draw,inner sep=2pt] at (0,0) {};}_{\pi}\,,\comul_{\pi}) $.
We use $\simeq$ to denote the \emph{identification} of basis elements
in the two bases $\{-\}$ and $[-]$ of indecomposables of the group and
its restricted subgroup. In graphical notation our method of deforming
Hopf algebras reads for an inverse pair of series $M_{\pi}(t)L_{\pi}(t)=1$
inducing the transformations between bases $\{$-$\}$ and $[$-$]$.
\begin{align}\label{grph-BasisChange}
   \begin{pic}
      \node (i1) at (0.5,1) {};
        \node[yshift=5pt] at (i1.north) {$\{$-$\}$ };
      \node (u) at (0.5,0.25) {};
      \node[rectangle,fill=white] (m1) at (0,-0.25) {$\simeq$ };
      \node[circle,draw,thick,fill=white,inner sep=2pt] (M) at (1,-0.25) {};
        \node[xshift=8pt] at (M.east) {$M_{\pi}$ };
      \node (o1) at (0,-1) {};
      \node[yshift=-5pt] at (o1.south) {$[$-$]$};
\begin{pgfonlayer}{background}
      \draw[thick] (i1.center) to (u.center);
      \draw[thick,out=0,in=90] (u.center) to (M.center);
      \draw[thick,out=180,in=90] (u.center) to (m1.center);
      \draw[thick] (m1.center) to (o1.center);
\end{pgfonlayer}
   \end{pic}
;&&&&
   \begin{pic}
      \node (i1) at (0.5,1) {};
        \node[yshift=5pt] at (i1.north) {$[$-$]$ };
      \node (u) at (0.5,0.25) {};
      \node[rectangle,fill=white] (m1) at (0,-0.25) {$\simeq$ };
      \node[circle,draw,thick,fill=white,inner sep=2pt] (M) at (1,-0.25) {};
        \node[xshift=8pt] at (M.east) {$L_{\pi}$ };
      \node (o1) at (0,-1) {};
      \node[yshift=-5pt] at (o1.south) {$\{$-$\}$};
\begin{pgfonlayer}{background}
      \draw[thick] (i1.center) to (u.center);
      \draw[thick,out=0,in=90] (u.center) to (M.center);
      \draw[thick,out=180,in=90] (u.center) to (m1.center);
      \draw[thick] (m1.center) to (o1.center);
\end{pgfonlayer}
   \end{pic}
;&&&&
   \begin{pic}
      \node (i1) at (0.5,1) {};
        \node[yshift=5pt] at (i1.north) {$\{$-$\}$ };
      \node (u) at (0.5,0.25) {};
      \node[rectangle,fill=white] (m1) at (0,-0.25) {$\simeq^{2}$ };
      \node (m2) at (1,-0.25) {};
      \node[circle,draw,thick,fill=white,inner sep=2pt] (M) at (0.5,-0.75) {};
      \node[circle,draw,thick,fill=white,inner sep=2pt] (L) at (1.5,-0.75) {};  
        \node[xshift=8pt] at (M.east) {$L_{\pi}$ };
        \node[xshift=8pt] at (L.east) {$M_{\pi}$ };
      \node (o1) at (0,-1) {};
        \node[yshift=-5pt] at (o1.south) {$\{$-$\}$};
\begin{pgfonlayer}{background}
      \draw[thick] (i1.center) to (u.center);
      \draw[thick,out=0,in=90] (u.center) to (m2.center);
      \draw[thick,out=180,in=90] (u.center) to (m1.center);
      \draw[thick] (m1.center) to (o1.center);
      \draw[thick,out=180,in=90] (m2.center) to (M.center);
      \draw[thick,out=0,in=90] (m2.center) to (L.center);
\end{pgfonlayer}
   \end{pic}
\cong
   \begin{pic}
      \node (i1) at (0,1) {};
        \node[yshift=5pt] at (i1.north) {$\{$-$\}$ };
      \node (o1) at (0,-1) {};
        \node[yshift=-5pt] at (o1.south) {$\{$-$\}$};
      \draw[thick] (i1.center) to (o1.center);
   \end{pic}              
\end{align}
This leads to the deformed comultiplication $\comul_{\pi}$ (used in~\cite{fauser:jarvis:king:2007c})
\begin{align}
   \begin{pic}
      \node (i1) at (0.5,1) {};
      \node (u) at (0.5,0.25) {};
      \node[circle,draw,thick,fill=white,inner sep=1pt] (c) at (0.5,0) {$\pi$};
      \node (o1) at (0,-1) {};
      \node (o2) at (1,-1) {};
\begin{pgfonlayer}{background}
      \draw[thick] (i1.center) to (c.center);
      \draw[thick,out=180,in=90] (c.center) to (o1.center);
      \draw[thick,out=0,in=90] (c.center) to (o2.center);
\end{pgfonlayer}      
   \end{pic}
&:\cong
   \begin{pic}
      \node (i1) at (2,1.3) {};
      \node (p1) at (2,1) {};
      \node (p2) at (1.25,0.5) {};
      \node[circle,draw,thick,fill=white,inner sep=2pt] (L) at (2.75,0.5) {};
         \node[yshift=-10pt] at (L.east) {$L_{\pi}$ };
      \node (p3) at (0.5,0) {};   
      \node (p4) at (2,0) {};
      \node[circle,draw,thick,fill=white,inner sep=2pt] (M1) at (1,-0.5) {};
      \node[circle,draw,thick,fill=white,inner sep=2pt] (M2) at (2.5,-0.5) {};  
        \node[yshift=-10pt] at (M1.east) {$M_{\pi}$ };
        \node[yshift=-10pt] at (M2.east) {$M_{\pi}$ };
      \node (o1) at (0,-1.1) {};
      \node (o2) at (1.5,-1.1) {};
\begin{pgfonlayer}{background}
      \draw[thick] (i1.center) to (p1.center);
      \draw[thick,out=0,in=90] (p1.center) to (L.center);
      \draw[thick,out=180,in=90] (p1.center) to (p2.center);
      \draw[thick,out=180,in=90] (p2.center) to (p3.center);
      \draw[thick,out=0,in=90] (p2.center) to (p4.center);
      \draw[thick,out=180,in=90] (p3.center) to (o1.center);
      \draw[thick,out=0,in=90] (p3.center) to (M1.center);
      \draw[thick,out=180,in=90] (p4.center) to (o2.center);
      \draw[thick,out=0,in=90] (p4.center) to (M2.center);
\end{pgfonlayer}      
   \end{pic}
\cong
   \begin{pic}
      \node (i1) at (0.5,1) {};
      \node (u) at (0.5,0.25) {};
      \node (m1) at (0,-0.25) {};
      \node (m2) at (1,-0.25) {};
      \node[circle,draw,thick,fill=white,inner sep=2pt] (M) at (1.5,-0.75) {};
         \node[xshift=8pt] at (M.east) {$M_{\pi}$ };
      \node (o2) at (0.5,-1) {};
      \node (o1) at (0,-1) {};
\begin{pgfonlayer}{background}
      \draw[thick] (i1.center) to (u.center);
      \draw[thick,out=0,in=90] (u.center) to (m2.center);
      \draw[thick,out=180,in=90] (u.center) to (m1.center);
      \draw[thick] (m1.center) to (o1.center);
      \draw[thick,out=0,in=90] (m2.center) to (M.center);
      \draw[thick,out=180,in=90] (m2.center) to (o2.center);
\end{pgfonlayer}
   \end{pic}
\cong
   \begin{pic}
      \node (i1) at (1,1) {};
      \node (u) at (1,0.25) {};
      \node (m1) at (0.5,-0.25) {};
      \node (m2) at (1.5,-0.25) {};
      \node[circle,draw,thick,fill=white,inner sep=2pt] (M) at (0,-0.75) {};
         \node[xshift=8pt] at (M.east) {$M_{\pi}$ };
      \node (o1) at (1,-1) {};
      \node (o2) at (1.5,-1) {};
\begin{pgfonlayer}{background}
      \draw[thick] (i1.center) to (u.center);
      \draw[thick,out=180,in=90] (u.center) to (m1.center);
      \draw[thick,out=0,in=90] (u.center) to (m2.center);
      \draw[thick,out=0,in=90] (m1.center) to (o1.center);
      \draw[thick] (m2.center) to (o2.center);
      \draw[thick,out=180,in=90] (m1.center) to (M.center);
\end{pgfonlayer}
   \end{pic}   
\end{align}
However, in what follows we are more interested in the application
of the theory developed in Section~\ref{sec-ComHopf} and the deformation
of multiplications.
%-----------------------------------------------------------------------
\section{Applications to tensor product decompositions of group characters}
\label{sec-Applications}
%-----------------------------------------------------------------------
As mentioned above, in previous work~\cite{fauser:jarvis:king:wybourne:2005a}
we developed character methods for deriving tensor product decompositions
for representations of subgroups of $\GL(N)$ for large $N$, including both
classical and non-classical algebraic subgroups (for which the representations
are in general indecomposable only). The methods employed were based on
techniques for manipulating plethysms and associated Schur function series
generated by them (see below). In this Section we apply the theory of
higher derived hash products, developed so far, to the theory of group
characters, and show that many of the classical decompositions can
instead be derived in this way (without explicitly using plethysms).
The goal of this Section is however to demonstrate how several classical
decomposition formulae are in fact deformations using a higher derived
hash prodcut, and hence that the subgroup of Laplace pairings parametrizes
these decompositions. As we consider polynomial $\GL(N)$ representations
and their characters, the proper Hopf algebra to use as ambient Hopf
algebra for characters is that of symmetric functions $\Sym$. We
therefore start by introducing the relevant notations to define this
Hopf algebra, as well as the associated subgroup characters.
%-----------------------------------------------------------------------
\subsection{Symmetric function Hopf algebra}
%-----------------------------------------------------------------------
\subsubsection{Notation}

Our notation follows in large part that of~\cite{macdonald:1979a}.
Integer partitions are specified by lower case Greek letters.
If $\lambda$ is a partition of $n$ we write $\lambda\vdash n$, and
$\lambda=(\lambda_1,\lambda_2,\ldots,\lambda_n)$ is a sequence of
non-negative integers $\lambda_i$ for $i=1,2,\ldots,n$
such that $\lambda_1\geq\lambda_2\geq\cdots\geq\lambda_n\geq0$,
with $\lambda_1+\lambda_2+\cdots+\lambda_n=n$. The partition
$\lambda$ is said to be of weight $|\lambda|=n$ and length
$\ell(\lambda)$ where $\lambda_i>0$ for all $i\leq\ell(\lambda)$
and $\lambda_i=0$ for all $i>\ell(\lambda)$. In specifying $\lambda$
the trailing zeros, that is those parts $\lambda_i=0$, are often
omitted, while repeated parts are sometimes written in exponent form
$\lambda=[1^{m_1},2^{m_2},\cdots]$ where $\lambda$ contains $m_i$
parts equal to $i$ for $i=1,2,\ldots$. For each such partition define
$n(\lambda)=\sum_{i=1}^n (i-1)\lambda_i$ and
$z_\lambda=\prod_{i\geq1} i^{m_i}\, m_i!$. Note that $\vert\lambda\vert
= \sum_i i m_i$.

Each partition $\lambda$ of weight $|\lambda|$ and length $\ell(\lambda)$
defines a Young or Ferrers diagram, $F^\lambda$, consisting of $|\lambda|$
boxes or nodes arranged in $\ell(\lambda)$ left-adjusted rows of lengths
from top to bottom $\lambda_1,\lambda_2,\ldots,\lambda_{\ell(\lambda)}$
(in the English convention). The partition $\lambda'$, conjugate to
$\lambda$, is the partition specifying the column lengths of $F^\lambda$
read from left to right. The box $(i,j)\in F^\lambda$, that is in
the $i$th row and $j$th column of $F^\lambda$, is said to have content
$c(i,j)=j-i$ and hook length $h(i,j)=\lambda_i+\lambda'_j-i-j+1$.
A box on the diagonal $(k,k)$ has arm length $a_k=\lambda_k-k$ and
leg length $b_k=\lambda_k^\prime-k$ for $1\leq k \leq r$, and $\lambda$
is said to have rank $r(\lambda)=r$. This allows partitions to be presented 
in Frobenius notation
\begin{align}
\lambda
  &=\left(\begin{array}{cccc}
        a_1 & a_2 & \ldots & a_r \\
        b_1 & b_2 & \ldots & b_r
    \end{array}\right)
&&&
\lambda^\prime
  &=\left(\begin{array}{cccc}
        b_1 & b_2 & \ldots & b_r \\
        a_1 & a_2 & \ldots & a_r
    \end{array}\right)
\end{align}
with $a_1> a_2 > \ldots a_r\geq 0$ and $b_1> b_2 > \ldots b_r\geq 0$.
By way of illustration, if the partition is
$\lambda=(4,2,2,1,0,0,0,0,0)=(4,2,2,1)=[1,2^2,4]$
then $|\lambda|=9$, $\ell(\lambda)=4$, $\lambda'=(4,3,1^2)$,
\begin{align}
F^\lambda=F^{(4,2^2,1)}
&=
\begin{pic}
   \node[rectangle,draw,inner sep=0.125cm] at (0,0) {};
   \node[rectangle,draw,inner sep=0.125cm] at (0.25,0) {};
   \node[rectangle,draw,inner sep=0.125cm] at (0.5,0) {};
   \node[rectangle,draw,inner sep=0.125cm] at (0.75,0) {};
   \node[rectangle,draw,inner sep=0.125cm] at (0,-0.25) {};
   \node[rectangle,draw,inner sep=0.125cm] at (0.25,-0.25) {};
   \node[rectangle,draw,inner sep=0.125cm] at (0,-0.5) {};
   \node[rectangle,draw,inner sep=0.125cm] at (0.25,-0.5) {};
   \node[rectangle,draw,inner sep=0.125cm] at (0,-0.75) {};
\end{pic}
&&\textrm{and}
&
F^{\lambda'}=F^{(4,3,1^2)} &=
\begin{pic}
   \node[rectangle,draw,inner sep=0.125cm] at (0,0) {};
   \node[rectangle,draw,inner sep=0.125cm] at (0.25,0) {};
   \node[rectangle,draw,inner sep=0.125cm] at (0.5,0) {};
   \node[rectangle,draw,inner sep=0.125cm] at (0.75,0) {};
   \node[rectangle,draw,inner sep=0.125cm] at (0,-0.25) {};
   \node[rectangle,draw,inner sep=0.125cm] at (0.25,-0.25) {};
   \node[rectangle,draw,inner sep=0.125cm] at (0.5,-0.25) {};
   \node[rectangle,draw,inner sep=0.125cm] at (0,-0.5) {};
   \node[rectangle,draw,inner sep=0.125cm] at (0,-0.75) {};
\end{pic}
\end{align}
The content and hook lengths of $F^\lambda$ are specified by
\begin{align}
&
\begin{pic}
   \node[rectangle,draw,minimum width=0.5cm] at (0,0)    {$0$};
   \node[rectangle,draw,minimum width=0.5cm] at (0.5,0)  {$1$};
   \node[rectangle,draw,minimum width=0.5cm] at (1.0,0)  {$2$};
   \node[rectangle,draw,minimum width=0.5cm] at (1.5,0)  {$3$};
   \node[rectangle,draw,minimum width=0.5cm,inner sep=2.5pt] at (0,-0.5) {$\overline{1}$};
   \node[rectangle,draw,minimum width=0.5cm] at (0.5,-0.5) {$0$};
   \node[rectangle,draw,minimum width=0.5cm,inner sep=2.5pt] at (0,-1) {$\overline{2}$};
   \node[rectangle,draw,minimum width=0.5cm,inner sep=2.5pt] at (0.5,-1) {$\overline{1}$};
   \node[rectangle,draw,minimum width=0.5cm,inner sep=2.5pt] at (0,-1.5) {$\overline{3}$};
\end{pic}
&&\textrm{and}&&
\begin{pic}
   \node[rectangle,draw,minimum width=0.5cm] at (0,0)    {$7$};
   \node[rectangle,draw,minimum width=0.5cm] at (0.5,0)  {$5$};
   \node[rectangle,draw,minimum width=0.5cm] at (1.0,0)  {$2$};
   \node[rectangle,draw,minimum width=0.5cm] at (1.5,0)  {$1$};
   \node[rectangle,draw,minimum width=0.5cm] at (0,-0.5) {$4$};
   \node[rectangle,draw,minimum width=0.5cm] at (0.5,-0.5) {$2$};
   \node[rectangle,draw,minimum width=0.5cm] at (0,-1) {$3$};
   \node[rectangle,draw,minimum width=0.5cm] at (0.5,-1) {$1$};
   \node[rectangle,draw,minimum width=0.5cm] at (0,-1.5) {$1$};
\end{pic}
\end{align}
where $\overline{m}=-m$ for all $m$.
In addition, $n(4,2^2,1)=0\cdot4+1\cdot2+2\cdot2+3\cdot1=9$ and
$z_{(4,2^2,1)}=4\cdot2^2\cdot1\cdot1!\cdot2!\cdot1!=32$.

The importance of using integer partitions to label vector spaces comes
from the fact, shown by I. Schur, that the irreducible finite dimensional
co- or contravariant representations of the general linear group $\GL(N)$
are labelled by integer partitions. The same partition labelling applies
to irreducible representations of the symmetric groups $S_n$, as was shown
by Frobenius.
%-----------------------------------------------------------------------
\subsubsection{Tensor product decomposition of irreducible
              (indecomposable) representations}

In what follows we will study the large $N$ (algebraic) limit of finite
dimensional representations and especially their characters. Let $V$ be
a (complex) vector space of dimension $N$. One studies the tensor
algebra $T(V)$, which is a graded Hopf algebra
$T(V) = \oplus_n V^{\ot^n}$ with a left $\GL(N)$ diagonal action. A
decomposition of this space into $\GL(N)$ irreducibles produces vector
spaces $V^\lambda \subset V^{\ot^n}$ with $n=\vert\lambda\vert$ the
weight of the partition and $\dim(V^{\lambda})=M$. Such spaces are
highest weight representations.

Then one has additionally the centralizer (algebra) of the action 
of $\GL(N)$ on $V^{\ot^n}$, which in this general linear case
amounts to the (group algebra of the) symmetric group $S_n$ acting from
the right on $V^{\ot^n}$ by permuting factors. Irreducible $S_n$ modules
$S^\lambda$ are labelled also by partitions.

Schur-Weyl duality states that the tensor algebra $T(V)$, $\dim(V)=N$
with a $\GL(N)$-left action and an $S_\infty$-right action (seen as
union of the $S_n$'s acting on $V^{\ot^n}$ from the right) decomposes
into
\begin{align}
  T(V) &= \sum_n \sum_{\lambda\vdash n} V^\lambda\ot S^\lambda\,.
\end{align}
It follows that the explicit decomposition of $V^{\ot^n}$, ignoring the
right $S_n$ action for the moment, into irreducible components with
respect to the left $\GL(N)$ action takes the form
\begin{align}
   V^{\ot^n} &= \sum_{\lambda\vdash n}\, \chi^\lambda(1)\, V^\lambda\,,
\end{align}
where the multiplicities, $\chi^\lambda(1)$, are the dimensions of 
the irreducible representations $S^\lambda$ of $S_n$ obtained by evaluating
their  characters $\chi^\lambda$ at the identity.

As the irreducibles span the whole space of representations, one is more
generally interested in finding such a decomposition for a tensor product
of irreducible representation spaces including multiplicities
$c^\lambda_{\mu,\nu}$
\begin{align}
   V^\mu \ot V^\nu &= \sum_{\lambda} c^\lambda_{\mu,\nu} V^\lambda
\end{align}
One needs only the isoclass of the vector space to compute this
decomposition, hence one can reduce the complexity of the problem
by using characters.

The character $\textrm{ch} V^\lambda(g)$ of a representation $V^\lambda${}
of dimension $M$ is given by the trace of a homomorphism
$\rho: \GL(N) \rightarrow \GL(M)$ inducing the diagonal $\GL(N)$ action
into $\GL(M)$. Hence for $g\in \GL(N)$ with representation $\pi^\lambda(g)${}
on an irreducible component $V^\lambda$ one gets
\begin{align}
 \textrm{ch} V^\lambda(g) = \textrm{Tr}_{V^\lambda}(\pi^\lambda(g)) &= s_\lambda(g)
\end{align}
where $s_\lambda$ is the character of $V^\lambda$ written as a (universal)
polynomial, the Schur polynomial, usually referred to as a Schur function.
This polynomial actually depends, due to the properties of the trace, only
on the invariants (latent roots or eigenvalues) of $g$. Moreover it is
invariant under a permutation of these invariants. We will denote the
invariants by a set of indeterminates $x_1,\ldots,x_N$, hence arriving
at a polynomial $s_\lambda(x_1,\ldots, x_N)$,
which for a sufficiently large number $N$ of indeterminates is universal.
The decomposition coefficients $c^\lambda_{\mu,\nu}$ appear as the
Littlewood Richardson coefficient in the polynomial ring of characters
\begin{align}
  \mul(s_\mu \ot s_\nu)
    &= s_\mu \cdot s_\nu
     = \sum_\lambda c^\lambda_{\mu,\nu} s_\lambda
\end{align}
We will see that this ring is freely generated by certain characters
$s_\lambda=s_\lambda(x_1,\ldots,)$ in a denumerably infinite number of
indeterminates $X=\{x_i\}_i^\infty$, which we call letters of an alphabet
$X$, and is a universal bicommutative Hopf algebra, fitting into the
theory developed in sections~\ref{sec-ComHopf} and~\ref{sec-HashProd}.

%-----------------------------------------------------------------------
\subsubsection{Restricted groups}

Before we do this we discuss briefly how this method can be used to
arrive at decompositions of irreducible or indecomposable characters
of so called restricted groups. The general linear and unitary groups
behave similarly in terms of their character theory. However, classical
groups such as the orthogonal or symplectic groups need further
treatment. A restricted group, in our sense, is an algebraic
subgroup of the general linear group defined by polynomial equations.
Among these are the classical groups
\begin{align}
\GL(N) &= \{ X \mid X \in \mathrm{Mat}(N,\mathbb{C}),~\det X \not= 0 \}
   &&& s_\lambda &= \{\lambda\}
  \nonumber \\
\SL(N) &= \{ X \mid X \in \GL(N),~\det X - 1 \}
   &&& s_\lambda &= \{\lambda\}
  \nonumber \\
\O(N) &= \{ X \mid X \in \GL(N),~~X^tX -\Id \}
   &&& o_\lambda &= [\lambda ]
  \nonumber \\
\Sp(N) &= \{ X \mid X\in \GL(N),~~XJX -J,~~J^t+J \}
   &&& sp_\lambda &= \la\lambda\ra
  \nonumber \\
\U(N) &= \{ X \mid X \in \GL(N),~~ X^{t*}X-\Id \}
   &&& s_\lambda &= \{\lambda\}
\end{align}
where $t$ is (matrix) transposition, $J$ is an antisymmetric matrix
(two form), and $*$ the star involution inherited from conjugation in
$\mathbb{C}$.
The right column gives the irreducible (for odd symplectic groups
indecomposable, see Proctor~\cite{proctor:1988a}) characters and
their notation using Littlewood's~\cite{littlewood:1940a}
(partly ambiguous but typographically convenient) bracket notation,
see also Section~\ref{subsec:OrthogonalSymplectic} below.

Another example of a restricted group is the symmetric group
$S_N \subset \GL(N)$ of permutation matrices. Here $S_N$ is seen as
a subgroup, by virtue of its representation by means of $N\times N$
permutation matrices, not as the (Weyl) group permuting tensor factors.
These symmetric groups $S_N$ stabilize symmetric tensors $R_i$, $R_{ij}$
and $R_{ijk}$ of symmetry type $(1)$, $(2)$ and
$(3)$~\cite{littlewood:1944a,littlewood:1958b}. 
Studying the union of all symmetric groups $S_N$ and taking an inductive
limit with respect to $N$ allows one to remove the $N$-dependence of
characters through a consideration of reduced characters, that are
unfortunately also written as $\la\lambda\ra$, see
Section~\ref{ssubsect:StableSymCharacters}.
The tensor product of reduced characters will be related to a
deformation of the inner product of symmetric functions.

Looking at the orthogonal groups, one notices that due to the existence
of an invariant tensor $g_{ij}=g_{ji}$ of symmetry type $(2)$ one
can extract traces with respect to $g_{ij}$ and all its concomitants
(freely generated algebraic products). For a detailed exposition of the
Hopf algebraic character theory for some classical groups
see~\cite{fauser:jarvis:king:2007c}. As an example a tensor
$V^{\ot^2}=V\ot V$ of rank 2 can be decomposed, under the action of $S_2$,
into symmetric and antisymmetric parts, we write $\cdot$ for character
multiplication
\begin{align}
  V\ot V &= V\vee V + V\wedge V;
  &&&
  [1]\cdot [1] \cong \{1\} \cdot \{1\} &= \{2\} + \{11\}
  \cong ([2] + [0]) + [11]
\end{align}
where the symmetric part $\{2\}$ can be further decomposed in the
orthogonal case into a trace free part $[2]$ and a trace part $[0]$.

In~\cite{fauser:jarvis:king:wybourne:2005a} we developed a theory of
group branchings which allows the tensor product decomposition formula
to be produced in such cases where a character branching is induced
by a (plethystic) Schur function series. Such a series induces a module
map $\Phi : \Sym[X] \longrightarrow \Sym[X]$ and the deformed,
in a Hopf algebraic sense, product of characters paralleling the tensor
product decomposition reads (using Sweedler notation for the
comultiplication $\Delta(A)=A_{(1)}\otimes A_{(2)}$ and the coboundary
operator $\partial$ from Definition~\ref{def-Sweedler-cohomology}, for
details we refer to~\cite{fauser:jarvis:king:wybourne:2005a})
\begin{align}\label{eq-PlethysticBranching}
A \circ B
  &= \Phi( \Phi^{-1}(A) \Phi^{-1}(B))
   = (\partial\Phi)(A_{(1)},B_{(1)})\, A_{(2)}B_{(2)}
\end{align}
Note that such general subgroup characters are in general only
indecomposable and not irreducible. However, the Hopf algebra
theory is independent of this property. While this process allows a direct
approach to such decompositions, it depends heavily on the notion of
plethysm. In the setting of noncommutative symmetric functions however,
the formulation using plethysm is unavailable as currently understood.
A main aim of the present work is to provide an alternative way to use
the subgroup of Laplace pairings to parameterize tensor product
decompositions without making direct use of plethystic module maps.

%-----------------------------------------------------------------------
\subsubsection{The universal Hopf algebra of symmetric functions}

Character polynomials $f$ are elements of a polynomial ring over the
integers generated by the commutative indeterminates or letters
$X^N=\{x_i\}_{i=1}^N = x_1+x_2+\ldots + x_N$, hence
$f(x_1,\ldots,x_N)\in \mathbb{Z}[x_1,\ldots,x_N]$.
However, a permutation $\pi\in S_N$ of the root variables leaves a
character invariant, and one is interested in the subring
$\Sym[X^N]=\mathbb{Z}[x_1,\ldots,x_N]^{S_N}$ of symmetric polynomials
$\pi\cdot f = f$. The first fundamental theorem of invariant theory
states that there is a basis of polynomial generators $e_n = s_{1^n}$
($e_n$ is the character of the antisymmetric $n$-th power $\wedge^n V$)
which freely generates this ring of symmetric functions
\begin{align}
  \Sym[X^N] 
    &= \mathbb{Z}[x_1,\ldots, x_N]^{S_N} = \mathbb{Z} [e_1,\ldots, e_N]
\end{align}
Letting $N$ tend to infinity, in the inductive limit, Schur polynomials
have the important stability property that $s_\lambda(x_1,\ldots,x_N,0,\ldots,0)
= s_\lambda(x_1,\ldots,x_N)$ for $N\geq \ell(\lambda)$. This stability
makes these polynomials universal. For small $N$ one encounters so-called
modification rules, which emerge due to the fact that certain characters
are zero (for example characters of the $k$-th exterior power $\wedge^k V$
of an $N$-dimensional space $V$ if $k>N$). The Hopf algebra development
used here assumes free modules, and modifications induce syzygies. For this
reason we work in the $N\rightarrow\infty$ limit and our characters are
formal universal characters.

It is well known that the ring of symmetric functions $\Sym$ has a Hopf
algebra structure~\cite{geissinger:1977a,zelevinsky:1981b,thibon:1991a,thibon:1991b,fauser:jarvis:2003a}.
Moreover, $\Sym$ is graded by the weight of partitions
$\Sym = \oplus_n \Sym^n$ and can be shown to be the universal bicommutative
self dual connected graded Hopf algebra. Self duality embodies Schur's lemma,
that between two isoclasses of vector spaces there exists either one isomorphism
or no map at all. This induces the Schur-Hall scalar product on $\Sym$
rendering Schur polynomials $s_\lambda$ orthonormal. The Hopf algebra is
defined, especially if using graphical calculus, in a basis free manner.
However, choosing the Schur function basis of irreducible $\GL$-characters,{}
we can summarize its structure as follows:
\mybenv{Theorem}\label{thm-symfuncHopfAlg}
The ring of symmetric functions $\Sym$ spanned by
$\{s_\lambda\}_{\lambda\in\mathcal{P}}$ together with the Schur-Hall
scalar product $\la -\mid -\ra : \Sym \times \Sym \rightarrow \mathbb{Z}$
is a bicommutative self dual graded connected Hopf algebra with the
following structure maps
\begin{align}\label{eq-SymHopf}
s_\mu \cdot s_\nu
   &= \sum_\lambda c^\lambda_{\mu,\nu}s_\lambda
   &&&&
   \textrm{(outer) multiplication}
\nn
\eta(1) &= s_{(0)}
   &&&&
   \textrm{unit~} \eta : \mathbb{Z}\rightarrow \Sym
\nn
\Delta(s_\lambda)
   &= \sum_{\mu,\nu} c_\lambda^{\mu,\nu} s_\mu\ot s_\nu
   &&&&
   \textrm{(outer) comultiplication}
\nn
\epsilon(s_\lambda)
   &=\delta_{\lambda,(0)}
   &&&&
   \textrm{counit~} \epsilon : \Sym\rightarrow \mathbb{Z}
\nn
\antip(s_\lambda)
   &= (-1)^{\vert\lambda\vert}s_{\lambda^\prime}
   &&&&
   \textrm{antipode}
\nn
\la s_\mu\mid s_\nu \ra
   &=\delta_{\mu,\nu}
   &&&&
   \textrm{Schur-Hall scalar product}
\end{align}
Grading is by weight of the partitions. Self duality implies that
$\la \Delta(s_\lambda) \mid s_\mu \ot s_\nu\ra = \la s_\lambda \mid
s_\mu\cdot s_\nu\ra$ and hence the numerical identity of the
coefficients $c^\lambda_{\mu,\nu}=c_\lambda^{\mu,\nu}$. The component
$\Sym^0=\mathbb{Z}$ hence $\Sym = \mathbb{Z} + \Sym^+$ with
$\Sym^+=\ker\epsilon$.
\myeenv
We define a \emph{skew Schur function} $s_{\mu/\nu}$ to be the adjoint
of multiplication by $s_\nu$ as $\la s_{\mu/\nu} \mid s_\lambda\ra
=\la s_\mu\mid s_\nu\cdot s_\lambda\ra$. As an operator we write
$s_\nu^\perp(s_\mu)=s_{\mu/\nu}$ and in Littlewood notation we have
$s_{\mu/\nu}=\{\mu/\nu\}$. Using these notions, the coproduct has
different forms which we use interchangeably
\begin{align}
\Delta(s_\lambda)
 &= s_{\lambda_{(1)}}\ot s_{\lambda_{(2)}}
  = \sum_{\mu,\nu} c^\lambda_{\mu,\nu} s_\mu \ot s_\nu
  = \sum_{\eta} s_{\mu/\eta} \ot s_{\eta}
  = \sum_{\eta} s_{\eta} \ot s_{\mu/\eta}
\end{align}
with the obvious similar forms in Littlewood bracket notation.

Sometimes it is convenient to extend the ring of symmetric functions
to $\widehat{\Sym}=\Sym\lsqb t\rsqb$, the ring of formal power series with
coefficients in $\Sym$. These series play an important role in the
theory of restricted
groups~\cite{littlewood:1940a,littlewood:1958b,murnaghan:1958a,%
fauser:jarvis:2003a,fauser:jarvis:king:wybourne:2005a,fauser:jarvis:king:2007c}.
We need especially the series
\begin{align}
\label{eq-Mser}
M(t)
   &:= \prod_i \frac{1}{1-x_it}  = \sum_{n\geq 0} h_n t^n   &&&& h_n=s_{(n)} 
\\
\label{eq-Lser}
L(t)
   &:= \prod_i {1-x_it}  = \sum_{n\geq 0} (-1)^n\, e_n t^n   &&&& e_n=s_{(1^n)}  
\end{align}
where the characters $h_n$ are complete symmetric functions.
and the $e_n$ are elementary symmetric functions.

Using the Frobenius notation for partitions, one can define the sets
of partitions, $\mathcal{P}$ all partitions, $\mathcal{D}=2\mathcal{P}$
all parts even, $\mathcal{B}$ conjugates of $\mathcal{D}$ and further
using
\begin{align}\label{eq-mathcalP}
\mathcal{P}_n
  &=\left\{
       \left(\begin{array}{cccc}
              a_1 & a_2 & \ldots & a_r \\
              b_1 & b_2 & \ldots & b_r
       \end{array}\right)
       \right\vert
              a_k-b_k=n\hskip2ex\textrm{for all}\hskip2ex
              \left.\begin{array}{rcl}
              r & = & 0,1,2,3,\ldots \\
              k & = & 1,2,\ldots,r\end{array}
              \right\}
\end{align}
we get $\mathcal{A}=\mathcal{P}_{-1}$, $\mathcal{C}=\mathcal{P}_{1}$
and $\mathcal{E}=\mathcal{P}_{0}$, the set of all self-conjugate partitions.
Using these sets of partitions we define for later use also the
series
\begin{align}
\label{eq-ABser}
A(t)
   &:= \prod_{i<j} \frac{1}{1-x_ix_jt}
     = \sum_{\alpha\in \mathcal{A}} 
       (-t)^{\vert\alpha\vert/2} \{\alpha\}
&&&
B(t)
   &:= \prod_{i<j} (1-x_ix_jt)
     = \sum_{\beta\in \mathcal{B}} 
       t^{\vert\beta\vert/2} \{\beta\}
\\
\label{eq-CDser}
C(t)
   &:= \prod_{i\leq j} \frac{1}{1-x_ix_jt}
     = \sum_{\gamma\in \mathcal{C}} 
       (-t)^{\vert\gamma\vert/2} \{\gamma\}
&&&
D(t)
   &:= \prod_{i\leq j} (1-x_ix_jt)
     = \sum_{\delta\in \mathcal{D}} 
       t^{\vert\delta\vert/2} \{\delta\}
\end{align}
As seen from the product expansion we have $M(t)L(t)=1$, $A(t)B(t)=1$,
$C(t)D(t)=1$.

Dually to a series we define linear forms on $\Sym$ by using the
Schur-Hall scalar product with lower case letters $\epsilon^1(f) = m(f) =
\la M(1) \mid f\ra$ and $l(f) = \la L(1) \mid f\ra$. The reason for
introducing the notion $\epsilon^1$ is twofold. It will serve as a counit for
another Frobenius comultiplication, and secondly applying linear forms
is equivalent to specializing the indeterminates. The Hopf algebra counit
$\epsilon=\epsilon^0$ acts as
$\epsilon(s_\lambda(x_1,\ldots,x_N)) = s_\lambda(0,\ldots,0) = \delta_{\lambda,(0)}$
can be seen to specialize all $x_i=0$, while $\epsilon^1(s_\lambda(x_1,\ldots,x_N))
= s_\lambda(1,0,\ldots,0) = \sum_{n\geq0}\delta_{\lambda,(n)}$ 
specializes $x_1=1$, and all other $x_i=0$.
More generally $\epsilon^d(s_\lambda(x_1,\ldots,x_N)=s_\lambda(1^d,0,\ldots)$
with $d$ ones provides the dimension of the $\GL(d)$-module $V^\lambda$
for $\dim(V)=d$.

As the character multiplication is also called `outer product' we call
this Hopf algebra $\Sym$ the outer Hopf algebra and will use it as the
\emph{ambient Hopf algebra} in the sense of the theory from
Section~\ref{sec-ComHopf}.

The outer Hopf algebra structure is tied to the $\GL(N)$ aspect of the
theory. We get another multiplication from that of irreducible $S_n$
characters $\chi^\lambda$ for each $n$. This takes the form
$\chi^\mu\,\chi^\nu = \sum_\lambda g^\lambda_{\mu,\nu}\ \chi^\lambda$,
where $\lambda,\mu,\nu$ are all partitions of $n$, and the non-negative
integers $g^\lambda_{\mu,\nu}$ are known as Kronecker coefficients.
Using the isometric Frobenius characteristic map 
$\textsf{ch}: R(S_n) \rightarrow \Sym^n ::\chi^\lambda \mapsto s_\lambda$, 
see~\cite{macdonald:1979a}, the so-called \emph{inner product} of Schur
functions is defined as
\begin{align}
s_\mu * s_\nu
   &= \sum_\lambda g^\lambda_{\mu,\nu} s_\lambda
   &&&&
   \textrm{inner multiplication}
\end{align}
The Kronecker coefficients $g^\lambda_{\mu,\nu}$ can be computed by the
Murnaghan-Nakayama formula~\cite{sagan:1991a}. We still can dualize the
multiplication to obtain a comultiplication $\delta$ on the coalgebra
on the dual space $\la \delta(s_\lambda)\mid s_\mu\ot s_\nu\ra :=
\la s_\lambda \mid s_\mu * s_\nu\ra$. 

There is no reason to assume that this comultiplication should fulfill
a bialgebra law~\eqref{grph-bialgebra-and-antipode} or be antipodal.
In fact in~\cite{fauser:jarvis:2003a} we demonstrated that the module
$\Sym$ together with the inner multiplication and inner comultiplication
cannot be a bialgebra and is not antipodal. However, inner multiplication
and inner comultiplication interact nicely with the ambient outer Hopf
algebra.
\mybenv{Theorem}
The inner multiplication is a Laplace pairing for the outer Hopf
algebra $\Sym$ (set $\sfa=*$ and $\mul=\cdot$ in~\eqref{eq-LaplaceRight}
and~\eqref{eq-LaplaceLeft}).
\myeenv
Noting from~\cite{fauser:jarvis:2003a} that the module $\Sym$ together
with the \emph{inner} multiplication and \emph{outer} comultiplication
does form a mixed bialgebra, one has by (the dual of)
Proposition~\ref{prop-FrobeniusBialgebra} the following
\mybenv{Corollary}
The inner multiplication $*$ with unit $M^n(1)=s_{(n)}$ and the inner
comultiplication $\delta$ with counit $\epsilon^1$ form for each
degree $n$ a Frobenius algebra on the module $\Sym^n$. Maps between
different degrees are zero.
\myeenv
Using series notation to collect all degrees, we denote the structure
maps of the inner Frobenius algebra as follows
\begin{align}
s_\mu * s_\nu
   &= \sum_\lambda g^\lambda_{\mu,\nu} s_\lambda
   &&&&
   \textrm{inner multiplication}
\nn
M(1) * s_\mu
   &= s_\mu
   &&&&
   \textrm{unit~} M(1) : \mathbb{Z} \rightarrow \oplus_n\Sym^n
\nn
\delta(s_\lambda)
   &= s_{\lambda[1]} \ot s_{\lambda[2]}
    = \sum_{\mu,\nu} g^\lambda_{\mu,\nu} s_{\mu}\otimes s_{\nu}
   &&&&
   \textrm{inner comultiplication}
\nn
\epsilon^1(s_\mu)
   &= \sum_{n\geq0}\,\delta_{\mu,(n)}
   &&&&
   \textrm{counit~} \epsilon^1 : \oplus_n \Sym^n \rightarrow \mathbb{Z}
\end{align}
where the Sweedler indices are now enclosed in rectangular brackets.
%-----------------------------------------------------------------------
\subsection{$\GL$-tensor product decompositions}
%-----------------------------------------------------------------------
\subsubsection{Tensor product decompositions of polynomial characters}

In the following subsections we provide a list of applications of
higher derived hash products to certain group-subgroup branching
processes. We will see that our theory covers a wide variety of these
cases and that we need indeed the generality of higher derived hash
products. We start with the simplest case of $\GL(N)$ polynomial
characters and their decomposition for various subgroup branchings
with $\GL(M)$ subgroups.

As mentioned, from now on the ambient Hopf algebra $\HB$ is understood
to be the Hopf algebra of symmetric functions $\Sym$. If interpreted as
the Hopf algebra of (universal polynomial) $\GL$ group characters,
we call it $\CGL$. As the Schur polynomials form the irreducible
$\GL$ characters, in this case the isomorphism is trivial. The
tensor product decomposition
\begin{align}
V^\mu \otimes V^\nu
  &= \oplus_\lambda \oplus^{c^\lambda_{\mu,\nu}} V^\lambda
\end{align}
provides the multiplicative structure of the character Hopf algebra
$\CGL=\Sym$. We identify the outer multiplication as the multiplication
used in $\CGL$ and for later use we identify the inner multiplication
$*=\sfa$ as a Frobenius Laplace pairing $\sfa$ on $\CGL$.
\begin{align}
\mul(s_\mu \otimes s_\nu) = s_\mu s_\nu
  &= \sum_\lambda c^\lambda_{\mu,\nu} s_\lambda
&&&
\sfa(s_\mu \otimes s_\nu) = s_\mu * s_\nu
  &= \sum_\lambda g^\lambda_{\mu,\nu} s_\lambda
\end{align}
Outer and inner comultiplications define, or are identified with,
additive and multiplicative branchings, respectively, whilst outer
multiplication itself defines a further diagonal subgroup branching.
In terms of representations the operations of multiplication and
comultiplication become induction and reduction functors. This complies
with the combinatorial intuition of Rota that multiplications assemble
things and comultiplications disassemble them. We study the branchings:
\begin{align}
\label{eq-additiveGLdown}
   &\GL(N+M) \downarrow \GL(N)\times \GL(M)
   &&\Rightarrow
   & s_\lambda
   &\mapsto \comul(s_\lambda) 
   = \sum c_{\mu,\nu}^\lambda s_\mu\ot s_\nu
\\
\label{eq-multiplicativeGLup}
   &\GL(NM) \downarrow \GL(N)\times \GL(M)
   &&\Rightarrow
   & s_\lambda
   &\mapsto \hskip3.5pt \delta(s_\lambda) 
   = \sum g_{\mu,\nu}^\lambda s_\mu\ot s_\nu
\\
\label{eq-additiveGLup}
   &\GL(N)\times \GL(N) \downarrow \GL(N){}
   &&\Rightarrow
   & s_\mu\ot s_\nu
   &\mapsto \hskip7pt s_\mu s_\nu 
   = \sum c_{\mu,\nu}^\lambda s_\lambda
%%\\
% &\red{BF: why~next~line~deleted?}\nonumber
% \\[3ex]
% %%FB
% \label{eq-multiplicativeGLdown}
%   &\GL(N)\times \GL(M) \uparrow \GL(NM)
%   &&\Rightarrow
%   & s_\mu\ot s_\nu
%   &\mapsto s_\mu * s_\nu = \sum g_{\mu,\nu}^\lambda s_\lambda 
\end{align}

These branching rules allow us to exhibit the Frobenius Laplace nature
of the action of $\mul$ and $\ast$ on $\CGL$. To this end consider the
restriction of $\GL(KM+KN)$ to its subgroup $\GL(K)\times\GL(M)\times\GL(N)$.{}
This may be accomplished in two ways: first via the group subgroup chain
\begin{align}\label{Eq-GL-KMN1}
  \GL(KM+KN)
    &\downarrow\GL(K)\times\GL(M+N)
     \downarrow\GL(K)\times\GL(M)\times\GL(N)
\end{align}
for which any element $z\in\CGL$ branches as follows
\begin{align}\label{Eq-CGL-zKMN1}
  z &\mapsto (\Id\ot\Delta)\circ\delta(z)\,,
\end{align}
and then via the subgroup chain
\begin{align}\label{Eq-GL-KMN2}
  \GL(& KM+KN)
    \downarrow\GL(KM)\times\GL(KN)
     \downarrow\GL(K)\times\GL(M)\times\GL(K)\times\GL(N)\nn
    &=\GL(K)\times\GL(K)\times\GL(M)\times\GL(N)
     \downarrow\GL(K)\times\GL(M)\times\GL(N)
\end{align}
for which the branching is given by
\begin{align}\label{Eq-CGL-zKMN2}
  z&\mapsto (\mul\ot\Id\ot\Id)\circ(\Id\ot\sw\ot\Id){}
       \circ(\delta\ot\delta)\circ\Delta(z)\,.
\end{align}
Since these maps must be identical for all $z\in\CGL$ it follows that
\begin{align}\label{Eq-dual-Laplace}
   (\Id\ot\Delta)\circ\delta 
     &= 
        (\mul\ot\Id\ot\Id)\circ(\Id\ot\sw\ot\Id){}
   \circ(\delta\ot\delta)\circ\Delta\,,
\end{align}
but this is just the dual of the Laplace property~\eqref{eq-LaplaceRight}
\begin{align}\label{Eq-Laplace}
  \ast\circ(\Id\ot\mul) 
  &= 
  \mul \circ(\ast\ot\ast)\circ(\Id\ot\sw\ot\Id)\circ(\Delta\ot\Id\ot\Id)\,.
\end{align}

Similarly, the restriction from $\GL(MN)\times\GL(MN)$ to the subgroup
$\GL(M)\times\GL(N)$ may proceed by two routes: first the group-subgroup
chain
\begin{align}\label{Eq-GL-MN1}
   \GL(MN)\times\GL(MN)
     &\downarrow\GL(MN)
      \downarrow\GL(M)\times\GL(M)
\end{align}
for which any element $z\in\CGL$ branches as follows
\begin{align}\label{Eq-CGL-zMN1}
   z &\mapsto \delta\circ\mul(z)\,,
\end{align}
and then via the subgroup chain
\begin{align}\label{Eq-GL-MN2}
  \GL(MN) &\times\GL(MN)
  \downarrow\GL(M)\times\GL(N)\times\GL(M)\times\GL(N)\nn
  &=
  \GL(M)\times\GL(M)\times\GL(N)\times\GL(N)\downarrow\GL(M)\times\GL(N)
\end{align}
for which the branching is given by
\begin{align}\label{Eq-CGL-zKMN2}
  z&\mapsto (\mul\ot\mul)\circ(\Id\ot\sw\ot\Id)\circ(\delta\ot\delta)\,.
\end{align}
Since once again these maps must be identical for all $z\in\CGL$ it follows that
\begin{align}\label{Eq-bialgebra}
   \delta\circ\mul{}
   &=
   (\mul\ot\mul)\circ(\Id\ot\sw\ot\Id)\circ(\delta\ot\delta)
\end{align}
but this is just the mixed bialgebra property~\eqref{grph-mixedBialgebra}
required of our Frobenius Laplace algebra. 

%---------------------------------------------------------------
\subsubsection{Stability of tensor product decompositions}
Recall the definition of the $M(t)$ series~\eqref{eq-Mser}. It is easy
to show that this series is group like $\Delta(M(t)) = M(t)\ot M(t)$ and
similarly for the inverse $L(t)$ series~\eqref{eq-Lser}. 
The $M(t)$ series includes all polynomial irreducible characters in 
the case of $\GL(1)$.

Specialising one variable $x_{N}=1$ in a Schur polynomial
$s_\lambda(x_1,...,x_N) = \{\lambda\} \in \CGL(N)$ reduces it to
$s_\lambda(x_1,...,x_{N-1},1)=\{\lambda\}_1 \in \CGL(N-1)$,
which induces the isomorphim, see~\eqref{grph-BasisChange},
\begin{align}
  \{\lambda\} &= \{\lambda/M\}_{1} 
  &&&
  \{\lambda\}_{1} &= \{\lambda/L\}
  &&&  
\end{align}
We get from this for the branching of characters
\begin{align}
\label{eq-GLminusoneDown}
  &\GL(N)\downarrow \GL(N-1)
  &&\Rightarrow
  & s_\lambda
  &\mapsto s_{\lambda/M}
    = {\sum}_{n,\mu} \langle s_\mu s_{(n)}, s_\lambda\rangle s_\mu
\\
\label{eq-GLminusoneUp}
  &\GL(N-1)\uparrow \GL(N)
  &&\Rightarrow
  & s_\lambda
  &\mapsto s_{\lambda/L}
    \hskip3pt = {\sum}_{n,\mu} (-1)^n \langle s_\mu s_{(1^n)}, s_\lambda\rangle s_\mu
\end{align}
While the characters for both groups differ $\{\lambda\}\not={}
\{\lambda\}_1 = \{\lambda/M\}$, the product rule stays the same, since
the inverse pair $M$, $L$ is group like.
\begin{align}
\mul_1(\{\mu\}_1\ot\{\nu\}_1)
  &= \mul(\{\mu/M\}\ot\{\nu/M\})
   =\{(\mu/M)\, (\nu/M)\}
\nn
  &\hskip-2ex\stackrel{\textrm{\tiny group~like}}{=}\{(\mu\nu)/M\}
\nn
  &=\{\mu\nu\}_1
\end{align}
We can hence devise the trivial hash product $\#_{\Id} = \mul_1 =
\sfe\conv\mul=\mul$ for product decompositions in $\GL(N-1)$.
\mybenv{Proposition}
The characters of $\GL(N)$ and $GL(N-1)$ have in the stable limit the
same product decomposition.
\myeenv
This is anything but surprising, as in the large $N$ limit there is no
difference between $N$ and $N-1$. For small $N$, however, there will
be a difference due to the different modification rules which are
sensitive to the exact number of indeterminates.

We close this discussion with a remark on \emph{modification rules}
needed for finite alphabets. The coproduct for $\CGL\cong\Sym$,
the $\GL$-character outer Hopf algebra, is obtained by splitting
the alphabet of a symmetric function additively. For a finite alphabet
$X^N$ expressed as a disjoint union $X^N=Y^R\cup Z^S$ with $N=R+S$,
and a function $f(x_1,\ldots,x_N)$ we get $\Delta(f)(x_1,\ldots,x_N) =
f(y_1,\ldots,y_R,z_1,\ldots,z_S) = f(y_1,\ldots,y_R)_{(1)}\ot{}
f(z_1,\ldots,z_S)_{(2)}$. In the inductive large $R$ and $S$ 
limit the splitting does not imply conditions on $Y^R$ and $Z^S$.

Modification rules for $\GL(K)$ are given by just projecting all
characters (Schur polynomials) $s_\lambda$ with partition length
$\ell(\lambda)>K$ to zero (for $K$ one of $N,R,S$ respectively),
which is possible due to the stability property of Schur polynomials.
For restricted groups, studied below, modification rules are more
complicated and we do not enter that realm. In addition we show in
Appendix~\ref{formalgroups} that the additive splitting is essentially
employing an additive formal group law in $N$ variables.

%-----------------------------------------------------------------------
\subsubsection{Mixed tensor product decomposition and rational characters}

If one intends to deal with mixed co- and contravariant representations,
that is with the $\GL(N)$ acting on a finite vector space $V$ of
dimension $N$ and its linear dual space $V^*$, one needs to extend the
characters to rational functions. The character ring of such mixed tensor
representations is freely generated by the polynomial characters and a
determinantal character $\varepsilon$. If $g\in \GL(N)$ acts on $V$, then
the contragredient representation corresponds to the action of $g^{-1}$
on $V^*$. Hence the characters on contravariant irreducible
spaces $V^{*,\lambda}$ have polynomial characters in the eigenvalues
$\overline{X}=\{\ox_i\}_{i=1}^{\infty}$ with $\ox_i= x_i^{-1}$. Multiplying
$\Sym$ by negative powers of the determinant $\varepsilon$ allows any character
$s_\lambda(\overline{X})$ to be expressed as polynomial characters times
powers of the inverse determinant. The ring of invariants has then the
structure $\Sym\lsqb \overline{\varepsilon}\rsqb$. Mixed rational
characters can hence be seen as elements of $\Sym \ot \overline{\Sym}$.
Hence in this particular subsection we \emph{change the ambient Hopf
algebra} to $\Sym^2=\Sym\ot\Sym$. The multiplication and comultiplication
on this space is given by
\begin{align}\label{eq-mul2comul2}
\mul_{\HBI{2}} &= (\mul\ot\mul)\circ ( 1 \ot \sw \ot 1)
\nn
\comul_{\HBI{2}} &= ( 1 \ot \sw \ot 1)\circ (\comul\ot\comul)
\end{align}
with obvious unit, counit and antipode. $\Sym^2$ is a connected graded
bicommutative Hopf algebra and our theory of Section~\ref{sec-ComHopf}
applies.

A further complication arises due to the fact that co- and contravariant
representations may be contracted, hence can be reduced, and the embedding
of the characters has to respect this. A similar process will be
encountered in the case of orthogonal and symplectic characters below.
For a definition of these rational characters
see~\cite{koike:1989a,king:1989a,fauser:jarvis:king:2007c}. 
The contravariant irreducible representation $V^{*,\mu}$, with
$\mu=(\mu_1,\mu_2,\ldots)$ a partition, is conveniently denoted by
$V^{\overline{\mu}}$, where $\overline{\mu} = (\ldots,-\mu_2,-\mu_1)${}
is its highest weight. Its character is given by the Schur polynomial
$s_\mu(\overline{X}) = s_{\overline{\mu}}(X)$. The mixed, co- and
contravariant representation  $V^\lambda\ot V^{*,\mu} = 
V^\lambda\ot V^{\overline{\mu}}$, specified by a pair of partitions
$\lambda$ and $\mu$, is in general reducible. It possesses an irreducible
constituent $V^{\lambda;\overline{\mu}}$ of highest weight 
$(\lambda;\overline{\mu})= (\lambda_1,\lambda_2,\ldots,0,
\ldots,0,\ldots,-\mu_2,-\mu_1)$. The corresponding rational characters
in $\Sym^2$ are variously denoted by 
$s_\lambda(X)\, s_\mu(\overline{X}) = s_\lambda(X)\, s_{\overline{\mu}}(X)
= s_\lambda s_{\overline{\mu}} = \{\lambda\}\ot\{\mu\}$
and $s_{(\lambda;\overline{\mu})}(X) = s_{(\lambda;\overline{\mu})}
=\{\lambda;\overline{\mu}\}$ in the reducible and irreducible cases,
respectively. The map from one to the other is provided by the
isomorphism
\begin{align}
\{\lambda\} \ot \{\overline{\mu}\}
  &=\sum_{\zeta\in\mathcal{P}} \{\lambda/\zeta; \overline{\mu/\zeta}\}
   =\{\lambda/M_{[1]}; \overline{\mu/M_{[2]}} \}
\end{align}
where we have used the notation $M=M(1)$ and the fact that
$\delta(M) = M_{[1]}\ot M_{[2]} = \sum_{\zeta\in\mathcal{P}} \zeta\ot\zeta$.
The skewing with respect to all partitions $\zeta\in\mathcal{P}$
(defined above in~\eqref{eq-mathcalP}) removes all possible contractions between co-
and contravariant spaces. By making use of the antipode~\eqref{eq-SymHopf}
one arrives at the inverse map
\begin{align}
\{\lambda;\overline{\mu}\}
  &=\sum_{\zeta\in\mathcal{P}} (-1)^{|\zeta|} \{\lambda/\zeta\}\, \{\overline{\mu/\zeta'}\}
   =\{\lambda/M_{[1]};\overline{\mu/\antip(M_{[2]})} \}
\end{align}

In graphical terms these isomorphisms read
\begin{align}\label{grph-ratIsos}
&
\begin{pic}
  \node (i1) at (0.9,1.5) {};
  \node (i2) at (1.1,1.5) {};
  \node (j1) at (0.9,0.7) {};
  \node (j2) at (1.1,0.7) {};
  \node (u1) at (0.5,0.25) {};
  \node (u2) at (1.5,0.25) {};
  \node[circle,inner sep=2pt,fill=black] (s) at (1,-0.2) {};
  \node (o1) at (0,-1.5) {};
  \node (o2) at (2,-1.5) {};
  \node[xshift=2pt,yshift=3pt] at (i1.north) {$\{\lambda;\overline{\mu}\}$};
  \node[yshift=-5pt] at (o1.south) {$\{\lambda/M_{[1]}\}$};
  \node[yshift=-5pt] at (o2.south) {$\{\overline{\mu/\antip(M_{[2]})}\}$};
  \node at (0.86,-1.8) {$\ot$};
%%%%%%%
  \draw[thick] (-0.2,0.7) rectangle (2.2,-0.6);
  \draw[thick] (i1.center) to (j1.center);
  \draw[thick] (i2.center) to (j2.center);
  \draw[thick,out=270,in=90] (j1.center) to (u1.center);
  \draw[thick,out=270,in=90] (j2.center) to (u2.center);
  \draw[thick,out=0,in=180] (u1.center) to (s.center);
  \draw[thick,out=180,in=0] (u2.center) to (s.center);
  \draw[thick,out=180,in=90] (u1.center) to (o1.center);
  \draw[thick,out=  0,in=90] (u2.center) to (o2.center);
\end{pic}
;&&&
\begin{pic}
  \node (i1) at (0.9,-1.5) {};
  \node (i2) at (1.1,-1.5) {};
  \node (j1) at (0.9,-0.7) {};
  \node (j2) at (1.1,-0.7) {};
  \node (u1) at (0.5,0.3) {};
  \node (u2) at (1.5,0.3) {};
  \node (s) at (1,-0.25) {};
  \node (o1) at (0,1.5) {};
  \node (o2) at (2,1.5) {};
  \node[xshift=2pt,yshift=-5pt] at (i1.south) {$\{\lambda/M_{[1]};\overline{\mu/M_{[2]}}\}$};
  \node[yshift=3pt] at (o1.north) {$\{\lambda\}$};
  \node[yshift=3pt] at (o2.north) {$\{\overline{\mu}\}$};
  \node at (1,1.75) {$\ot$};
%%%%%%%
  \draw[thick] (-0.2,0.6) rectangle (2.2,-0.7);
  \draw[thick] (i1.center) to (j1.center);
  \draw[thick] (i2.center) to (j2.center);
  \draw[thick,out=90,in=180] (j1.center) .. controls +(0,0.2) and +(-1,-0.4) .. (u1.center);
  \draw[thick,out=90,in=0] (j2.center) .. controls +(0,0.2) and +(1,-0.4) .. (u2.center);
  \draw[thick,out=0,in=180] (u1.center) to (s.center);
  \draw[thick,out=180,in=0] (u2.center) to (s.center);
  \draw[thick,out=90,in=270] (u1.center) to (o1.center);
  \draw[thick,out=90,in=270] (u2.center) to (o2.center);
\end{pic}
\end{align}
where double lines represent mixed characters and the two separated
lines represent elements in $\Sym^2=\Sym\ot\Sym$. Operations inside the
box are in $\Sym$. The contractions are done with respect to the
Schur-Hall scalar product and its convolutive inverse defined by applying
the antipode.
\begin{align}\label{grph-sym2scalar}
   \begin{pic}
      \node (i1) at (0,0.5) {};
      \node (i2) at (1,0.5) {};
      \node[circle,fill=black,inner sep=2pt] (u) at (0.5,-0.5) {};
      \draw[thick,out=270,in=180] (i1.center) to (u.center);
      \draw[thick,out=270,in=  0] (i2.center) to (u.center);
   \end{pic}
&:\cong
   \begin{pic}
      \node (i1) at (0,0.5) {};
      \node (i2) at (1,0.5) {};
      \node[circle,draw,thick,inner sep=1pt,fill=white] (s1) at (0,0) {$\antip$};
      \node (s2) at (1,0) {};
      \node (u) at (0.5,-0.75) {};
\begin{pgfonlayer}{background}
      \draw[thick] (i1.center) to (s1.center);
      \draw[thick,out=270,in=180] (s1.center) to (u.center);
      \draw[thick] (i2.center) to (s2.center);
      \draw[thick,out=270,in=  0] (s2.center) to (u.center);
\end{pgfonlayer}
   \end{pic}
\cong
   \begin{pic}
      \node (i1) at (0,0.5) {};
      \node (i2) at (1,0.5) {};
      \node[circle,draw,thick,inner sep=1pt,fill=white] (s2) at (1,0) {$\antip$};
      \node (s1) at (0,0) {};
      \node (u) at (0.5,-0.75) {};
\begin{pgfonlayer}{background}
      \draw[thick] (i1.center) to (s1.center);
      \draw[thick,out=270,in=180] (s1.center) to (u.center);
      \draw[thick] (i2.center) to (s2.center);
      \draw[thick,out=270,in=  0] (s2.center) to (u.center);
\end{pgfonlayer}
   \end{pic}
;&&&
   \begin{pic}
      \node (i1) at (0.25,0.5) {};
      \node (i2) at (0.75,0.5) {};
      \node (u1) at (0.25,0.3) {};
      \node (u2) at (0.75,0.3) {};
      \node (m1) at (0,0) {};
      \node (m4) at (1,0) {};
      \node[circle,fill=black,inner sep=2pt] (d1) at (0.25,-0.3) {};
      \node (d2) at (0.75,-0.3) {};
      \draw[thick] (i1.center) to (u1.center);
      \draw[thick] (i2.center) to (u2.center);
      \draw[thick,out=180,in=90] (u1.center) to (m1.center);
      \draw[thick,out=  0,in=180] (u1.center) to (d2.center);
      \draw[thick,out=180,in=  0] (u2.center) to (d1.center);
      \draw[thick,out=  0,in=90] (u2.center) to (m4.center);
      \draw[thick,out=270,in=180] (m1.center) to (d1.center);
      \draw[thick,out=270,in=  0] (m4.center) to (d2.center);
   \end{pic}
&\cong
   \begin{pic}
      \node (i1) at (0,0.5) {};
      \node (i2) at (0.5,0.5) {};
      \node[circle,inner sep=2pt,draw,thick] (eps1) at (0,-0.3) {};
      \node[xshift=3pt] at (eps1.east) {$\epsilon$};
      \node[circle,inner sep=2pt,draw,thick] (eps2) at (0.5,-0.3) {};
      \node[xshift=3pt] at (eps2.east) {$\epsilon$};
      \draw[thick] (i1.center) to (eps1);
      \draw[thick] (i2.center) to (eps2);
   \end{pic}
;&&&
m^2:=
   \begin{pic}
      \node (i1) at (0,0.6) {};
      \node (i2) at (0.5,0.6) {};
      \node (i3) at (1,0.6) {};
      \node (i4) at (1.5,0.6) {};
      \draw[thick,out=270,in=270] (i1) to (i4);
      \draw[thick,out=270,in=270] (i2) to (i3);
   \end{pic}
\end{align}
The right most tangle depicts the scalar product
$m^{2}=\epsilon^1\circ*\circ(1\ot\epsilon^{1}\ot 1)\circ(1\ot *\ot 1)$
on $\Sym^2$ obtained from `bending up' two lines. Using these tools we
can show by graphical manipulations that the formula which governs
the (outer) products of such rational characters, is in fact given by
a derived hash product:
\mybenv{Theorem}
The product formula for rational tensor characters~\cite{fauser:jarvis:king:2007c}
\begin{align}\label{eq-ratGLalg}
 \{\kappa;\overline{\lambda}\}\cdot
 \{\mu;\overline{\nu}\}
 &=
 \sum_{\sigma,\tau\in\mathcal{P}}
 \{(\kappa/\sigma)\cdot (\mu/\tau);\overline{(\lambda/\tau)\cdot(\mu/\sigma)}\}
\end{align}
is given by a derived hash product $\#_{m^2,1^2}$ on $\Sym^2$
\begin{align}\label{eq-ratGLhash}
 \{\kappa;\overline{\lambda}\}\cdot
 \{\mu;\overline{\nu}\}
 &=
 \{ (\kappa;\overline{\lambda}) \#_{m^2,1^2} (\mu;\overline{\nu}) \}
\end{align}
where $m^2$ is the scalar product (derived Laplace
pairing~\eqref{grph-sym2scalar}) on $\Sym^2$ and $1^2$ the identity.
\myeenv

\noindent
\textbf{Proof:} The convolution in use is that of $\Sym^2$ with product
and coproduct defined in~\eqref{eq-mul2comul2}. Using the
isomorphisms~\eqref{grph-ratIsos} we can act with the multiplication
from $\Sym^2\ot\Sym^2$, then apply the inverse
of~\eqref{grph-ratIsos}. Reorganizing the tangle in several steps using
the bialgebra law~\eqref{grph-bialgebra-and-antipode}, biassociativity
and bicommutativity cancels the two inverse scalar
products~\eqref{grph-sym2scalar} and leaves one with the pairing
$m^2$ from~\eqref{grph-sym2scalar}. A further reorganization of the
tangle gives a graphical representation of the derived hash
product $\#_{m^2,1^2}$ on $\Sym^2$. Inserting the definitions in terms
of partitions from Theorem~\ref{thm-symfuncHopfAlg} yields the
algebraic form of the result~\eqref{eq-ratGLalg} and shows the abstract
form~\eqref{eq-ratGLhash}.
\qed

The reader is, however, encouraged to draw the respective tangles and
check this result graphically. The nature of forming hash products
becomes clearer in the next example showing how to handle orthogonal
and symplectic group characters using $\Sym$ itself.
%-----------------------------------------------------------------------
\subsection{Orthogonal and symplectic characters}\label{subsec:OrthogonalSymplectic}
%-----------------------------------------------------------------------

Orthogonal and symplectic groups provide another interesting case. We
denote the irreducible, or for odd symplectic groups indecomposable,
characters by $s_\lambda=\{\lambda\}\in\CGL$ and $os_\lambda=[\lambda]\in\CO$,
$sp_\lambda = \langle\lambda\rangle\in\CSp$, and we utilize the $C,D$
series~\eqref{eq-CDser} and $A,B$ series~\eqref{eq-ABser} to provide the
isomorphisms between these characters. Then the branchings read
\begin{align}
  &\GL(N)\downarrow \mathsf{O}(N)
  &&\Rightarrow
  & s_\lambda
  &\mapsto [\lambda/D] = os_{\lambda/D}
    = {\sum}_{n,\delta\in\mathcal{D}} 
      \langle s_\mu s_{\delta}, s_\lambda\rangle s_\mu
\\
  &\mathsf{O}(N)\uparrow \GL(N)
  &&\Rightarrow
  & os_\lambda
  &\mapsto \{\lambda/C\} = s_{\lambda/C}
    = {\sum}_{n,\gamma \in\mathcal{C}} (-1)^{\vert\gamma\vert/2}
             \langle s_\mu s_{\gamma}, s_\lambda\rangle s_\mu
\\[2ex]
  &\GL(N)\downarrow \mathsf{Sp}(N)
  &&\Rightarrow
  & s_\lambda
  &\mapsto \langle \lambda/B\rangle = s_{\lambda/B}
    = {\sum}_{n,\beta\in\mathcal{B}} \langle s_\mu s_{\beta}, s_\lambda\rangle s_\mu
\\
  &\mathsf{Sp}(N)\uparrow \GL(N)
  &&\Rightarrow
  & sp_\lambda
  &\mapsto \{\lambda/A\}= s_{\lambda/A}
    = {\sum}_{n,\alpha \in\mathcal{A}} (-1)^{\vert\alpha\vert/2}
             \langle s_\mu s_{\alpha}, s_\lambda\rangle s_\mu
\end{align}
The problem in the decomposition of products of these characters is now,
that both characters, say $[\mu]$ and $[\nu]$ in the orthogonal case,
are fully reduced. That is all traces with respect to the tensor $g_{ij}$
have been removed. However, in multiplying them (tensoring the
representations) one has to remove the traces between the two characters
(representations). This is essentially the way the following branching
result was first obtained by Littlewood. This reflects the fact that the
series $A$, $B$, $C$, and $D$ are not group like, e.g. 
$\Delta(D) = (D\ot D)\,\Delta^\prime(D)$, with the proper cut part of this
coproduct eliminating mixed traxes between $[\mu]$ and $[\nu]$.
\mybenv{Theorem}\textbf{[Newell-Littlewood]} The product decompositions
for orthogonal and symplectic characters are given by
\begin{align}
  [\mu ]\,\cdot\, [\nu ]
  &= \sum_{\zeta} [ (\mu/\zeta)\,\cdot\, (\nu/\zeta) ]
  &&\textrm{and}&&
  \langle\mu \rangle\,\cdot\, \langle\nu\rangle
  = \sum_{\zeta} \langle (\mu/\zeta)\,\cdot\, (\nu/\zeta) \rangle
\end{align}
where the sums are over all partitions $\zeta\in\mathcal{P}$.
\myeenv
We set $\mul_2([\mu ]\ot [\nu ]) = [\mu]\,\cdot\, [\nu]$ and
$\mul_{11}(\langle\mu \rangle\ot \langle\nu\rangle)
= \langle\mu \rangle\,\cdot\, \langle\nu\rangle$ for the two products. For the
origin of this naming convention see~\cite{fauser:jarvis:king:wybourne:2005a}.
We need to show that there is a higher derived hash product $\#_{\mathcal{L}}$
on $\CGL=\Sym$, such that $\mul_2([\mu]\ot[\nu]) = [\mu \#_{\mathcal{L}} \nu]$
and an identical product for the symplectic case.
\mybenv{Proposition}
The Newell-Littlewood product decomposition for orthogonal and symplectic
group characters is given by the derived hash products
\begin{align}
 \#_{m,1} ([a] \ot [b])
   &= [a \#_{m,1} b]
    = \la M\mid a_{(1)} * b_{(1)}\ra\, [ a_{(2)}b_{(2)} ]
    = \mul_2([a] \ot [b])
   \\
 \#_{m,1} (\langle a\rangle \ot \langle b\rangle)
   &= \langle a \#_{m,1} b\rangle
    = \la M\mid a_{(1)} * b_{(1)}\ra\, \langle a_{(2)}b_{(2)}\rangle{}
    = \mul_{1,1} ([a] \ot [b])
\end{align}
The $\phi_i$ maps are $m=\eta\circ\epsilon^1$ and $\Id$, and are easily
shown to be 1-cocycles.
\myeenv

\noindent
\textbf{Proof:} We first give an algebraic proof.
\begin{align}
 \#_{m,1} ([\mu] \ot [\nu])
   &= \la M\mid\mu_{(1)} * \nu_{(1)}\ra [ \mu_{(2)}\nu_{(2)} ]
   \nn
   &={\sum}_{\rho,\zeta} \la M\mid\rho * \zeta\ra [ (\mu/\rho)\, (\nu/\zeta) ]
    ={\sum}_{\rho,\zeta} \la\rho\mid\zeta\ra [ (\mu/\rho)\, (\nu/\zeta) ]
   \nn
   &={\sum}_{\zeta} [(\mu/\zeta)\, (\nu/\zeta)]
\end{align}
since $\la M\mid\rho*\zeta\ra = \la\rho\mid\zeta\ra = \delta_{\rho,\zeta}$.
The resulting expression is the right hand side of the Newell--Littlewood
formula, as required. The symplectic case is identical.
\qed

It may be instructive to see how this is obtained graphically. First we
decompose the derived hash product $\#_{m,\Id}$. The derived
pairing $a_\phi=\eta\circ\epsilon^1\circ\sfa$ is actually the Schur-Hall
scalar product, as can be seen from $\epsilon^1(f) = \la M(1), f\ra$.
Then we reorganize the tangle (recall $\sfa=*$)
\begin{align}
   \begin{pic}
      \node (i1) at (0,1) {};
      \node (i2) at (1,1) {};
      \node[circle,draw,inner sep=0pt,fill=white] (a) at (0.5,0) {$\#_{m,\Id}$ };
      \node (o) at (0.5,-1) {};
\begin{pgfonlayer}{background}
      \draw[thick,out=270,in=180] (i1.center) to (a.center);
      \draw[thick,out=270,in=0] (i2.center) to (a.center);
      \draw[thick] (a.center) to (o.center);
\end{pgfonlayer}
   \end{pic}
\cong
   \begin{pic}
      \node (i1) at (0,1) {};
      \node (i2) at (4,1) {};
      \node[rectangle,draw,fill=white] (a) at (2,0) {$(\eta\circ(\epsilon^1\circ *))\conv\mul$ };
      \node (o) at (2,-1) {};
\begin{pgfonlayer}{background}
      \draw[thick,out=270,in=165] (i1.center) to (a.north west);
      \draw[thick,out=270,in=15] (i2.center) to (a.north east);
      \draw[thick] (a.center) to (o.center);
\end{pgfonlayer}
   \end{pic}
\cong
   \begin{pic}
      \node (i1) at (0.25,1) {};
      \node (i2) at (1,1) {};
      \node (m1) at (0.25,0.75) {};
      \node (m2) at (1,0.75) {};
      \node[circle,draw,inner sep=1pt,fill=white] (a1) at (0.25,0.3) {$*$ };
      \node[circle,inner sep=2pt,draw,thick] (eps) at (0.25,-0.1) {};
      \node[xshift=4pt,yshift=2pt] at (eps.east) {$\epsilon^1$};
      \node[circle,inner sep=2pt,draw,thick] (eta) at (0.25,-0.4) {};
      \node[xshift=3pt] at (eta.east) {$\eta$};
      \node (a2) at (1,0) {};
      \node (d) at (0.6125,-0.75) {};
      \node (o) at (0.6125,-1) {};
\begin{pgfonlayer}{background}
      \draw[thick] (i1.center) to (m1.center);
      \draw[thick] (i2.center) to (m2.center);
      \draw[thick,out=180,in=180] (m1.center) to (a1.center);
      \draw[thick,out=  0,in=180] (m1.center) to (a2.center);
      \draw[thick,out=180,in= 0] (m2.center) to (a1.center);
      \draw[thick,out=  0,in=0] (m2.center) to (a2.center);
      \draw[thick,out=270,in=180] (eta.south) to (d.center);
      \draw[thick,out=270,in=  0] (a2.center) to (d.center);
      \draw[thick] (a1.south) to (eps.north);
      \draw[thick] (d.center) to (o.center);
\end{pgfonlayer}
   \end{pic}
\cong
   \begin{pic}
      \node (i1) at (0.25,1) {};
      \node (i2) at (1,1) {};
      \node (m1) at (0.25,0.75) {};
      \node (m2) at (1,0.75) {};
      \node (a1) at (0.25,0) {};
      \node (a2) at (1,0) {};
      \node (o) at (1,-1) {};
\begin{pgfonlayer}{background}
      \draw[thick] (i1.center) to (m1.center);
      \draw[thick] (i2.center) to (m2.center);
      \draw[thick,out=180,in=180] (m1.center) to (a1.center);
      \draw[thick,out=  0,in=180] (m1.center) to (a2.center);
      \draw[thick,out=180,in= 0] (m2.center) to (a1.center);
      \draw[thick,out=  0,in=0] (m2.center) to (a2.center);
      \draw[thick] (a2.center) to (o.center);
\end{pgfonlayer}
   \end{pic}
\cong
   \begin{pic}
      \node (i1) at (0.25,1) {};
      \node (i2) at (1,1) {};
      \node (m1) at (0.25,0.75) {};
      \node (m2) at (1,0.75) {};
      \node (a1) at (0.6125,0.25) {};
      \node (a2) at (0.6125,-0.25) {};
      \node (o) at (0.6125,-1) {};
\begin{pgfonlayer}{background}
      \draw[thick] (i1.center) to (m1.center);
      \draw[thick] (i2.center) to (m2.center);
      \draw[thick,out=  0,in=180] (m1.center) to (a1.center);
      \draw[thick,out=180,in=180] (m1.center) .. controls +(-0.3,0) and +(-0.9,0).. (a2.center);
      \draw[thick,out=180,in=  0] (m2.center) to (a1.center);
      \draw[thick,out=  0,in=  0] (m2.center) .. controls +(0.3,0) and +(0.9,0).. (a2.center);
      \draw[thick] (a2.center) to (o.center);
\end{pgfonlayer}
   \end{pic}
\end{align}
The last tangle was baptized the \emph{Rota sausage}\footnote{%
   By Zbigniew Oziewicz at ICCA5, Ixtapa, 1999.}
in the Rota-Stein context of cliffordization of a Grassmann Hopf
algebra~\cite{fauser:2002c}. The line connecting the two input strands
actually extracts all mixed traces of the two input characters
(orthogonal or symplectic does not matter here, the contraction tensor
could even be singular) and is a Wick type formula~\cite{fauser:2001b},
as it extracts terms with respect to a degree two covariance.
%-----------------------------------------------------------------------
\subsection{Thibon characters : multiplicative formal group}
%-----------------------------------------------------------------------

Our next example is an intermediate product decomposition needed later
in dealing with the inner product of stable symmetric function characters,
but is of interest in its own right. The characters in this example are
no longer polynomial but are in the realm of $\Sym\lsqb t\rsqb$. We
employ again the series $M=M(1)$ and $L=L(1)$.

\mybenv{Definition}\label{def-ThibonCharacter}
For all partitions $\lambda$ a Thibon character $\lla \lambda \rra$ is
given by the isomorphism
\begin{align}\label{eq-ThibonCharacter}
\lla \lambda\rra
  &= \{\lambda\,M\}
&&\mbox{with}&
\{\lambda\}
  &= \lla \lambda\,L\rra \,.
\end{align}
\vskip-1ex
\myeenv
Thibon characters were introduced
in~\cite{thibon:1991a,scharf:thibon:wybourne:1993a}. They are related in
the case of complete symmetric functions $\lla h_{\lambda}\rra$ to stable
permutation characters and to Young
polynomials~\cite{specht:1960a,kerber:1992a}.

Note that the \emph{outer} multiplication and comultiplication of
Thibon characters are given by
\begin{align}
	\lla\mu\rra \cdot\lla\nu\rra 
	   &= \lla \mu\nu\,M\rra
   &&\mbox{and}&
	\comul \lla \lambda\rra 
	   &= \lla \lambda_{(1)}\rra \ot \lla\lambda_{(2)}\rra{}
\end{align}
since   
\begin{align}
  \lla\mu\rra \cdot\lla\nu\rra
     &=\{\mu\,M\cdot\nu\,M\}=\{\mu\,\nu\,M\cdot M\}
      =\lla \mu\,\nu\,M\rra
\end{align}
and
\begin{align}
  \comul \lla \lambda\rra 
    &=\comul \{\lambda\,L\}=\comul \{\lambda\}\cdot \comul L 
     = \{\lambda_{(1)}\}\ot\{\lambda_{(2)}\}\cdot L\ot L
     = \lla \lambda_{(1)}\rra \ot \lla\lambda_{(2)}\rra\,.
\end{align}
This is the Hopf algebra structure of $\CGL^{*}(N-1)$, the dual Hopf
algebra of $\CGL(N-1)$ (in the stable limit). Comparing
with~\eqref{eq-GLminusoneDown} and~\eqref{eq-GLminusoneUp} we see that
this time the comultiplication is unchanged while the multiplication is
altered. For similar dualities in the orthogonal and symplectic case
see~\cite{fauser:jarvis:king:2007c}.

A more important problem is to calculate inner products of Thibon
characters, that is to compute the coefficients
$\widetilde{g}^\lambda_{\mu,\nu}$ in
\begin{align}
\lla\mu\rra * \lla\nu\rra
  &= \sum_{\lambda} \widetilde{g}^\lambda_{\mu,\nu} \lla\lambda\rra \,.
\end{align}
This can be done by means of
\mybenv{Theorem} \textbf{[Thibon]}\label{thm-innerThibon}
The inner product decomposition of Thibon characters is given by 
\begin{align}
  \lla \mu\rra * \lla \nu\rra
  &= \sum_{\sigma,\tau\in\mathcal{P}} 
     \lla (\sigma \ast \tau)\cdot (\mu/\sigma)\cdot (\nu/\tau) \rra \,,
\end{align}
where the sum is over all pairs of partitions, $\sigma$ and $\tau$, of
the same weight.
\myeenv

In terms of Schur function manipulations this implies that
\begin{align}
  \mul\circ (*\ot \mul)\circ 
  (\Id\ot\sw\ot\Id)\circ(\Delta\ot\Delta) (s_\mu \ot s_\nu){}
     &= \sum_{\lambda}\ \widetilde{g}^\lambda_{\mu,\nu}\ s_\lambda \,,
\end{align}
where $\mul$ and $\Delta$ denote outer multiplication and comultiplication,
respectively, and $*$ denotes inner multiplication.

It follows by interchanging algebra and coalgebra operators that
\begin{align}
  (\mul\ot \mul)\circ{}
  (\Id\ot\sw\ot\Id)\circ{}
  (\delta\ot\Delta)\circ\Delta (s_\lambda)
    &= \sum_{\mu,\nu}\ \widetilde{g}^\lambda_{\mu,\nu}
       \,s_\mu \ot s_\nu \,,
\end{align}
where $\delta$ denotes inner comultiplication. This is precisely the
stable $\CGL$ branching rule for the group-subgroup chain
\begin{align}
\GL(MN+M+N)&\supset\GL(MN)\times\GL(M+N)\cr
           &\supset(\GL(M)\times\GL(N))\times(\GL(M)\times\GL(N))\cr
           &\supset(\GL(M)\times\GL(M))\times(\GL(N)\times\GL(N))\cr
           &\supset\GL(M)\times\GL(N)\,.
\end{align}

We give here directly Thibon's Theorem~\ref{thm-innerThibon} expressed
in terms of a hash product.
\mybenv{Proposition}\label{prop-innerThibon}
The inner multiplication of Thibon characters is given by the
hash product $\#_{1,1} = \sfa\conv\mul$, where $\sfa=*$ is the
Frobenius Laplace inner multiplication. That is to say
\begin{align}
 \lla x\rra * \lla y\rra
  &= \lla (x_{(1)} * y_{(1)})\,x_{(2)}\,y_{(2)} \rra
   = \lla x\,\#_{1,1}\,y\rra
   = \#_{1,1} \lla x\rra \ot \lla y\rra
\end{align}
The coalgebra $(\Sym,\Delta,\epsilon)$ and the algebra of Thibon
characters $(\Sym,\#_{1,1},\eta)$ under the hash product $\#_{1,1}$,
with unit $\eta$ forms a Hopf algebra.
\myeenv

\noindent
\textbf{Proof:}
Expanding the characters $\lla x\rra * \lla y\rra = \{\mu\,M\} * \{\nu\,M\}$
gives a situation where we can use the left~\eqref{eq-LaplaceLeft}
and twice the right~\eqref{eq-LaplaceRight} Laplace expansions. The
result reads~\eqref{eq-cummins}
$\{ (\mu_{(1)}*\nu_{(1)}) \cdot
    (\mu_{(2)}*M)\cdot
    (\nu_{(2)}*M)\cdot (M*M) \}$,
where we have used the fact that $M$ is group like. Recalling that $M$
is, grade by grade, the unit for the inner multiplication, this reduces
to $\{ (\mu_{(1)}*\nu_{(1)})\cdot
        \mu_{(2)} \cdot \nu_{(2)} \cdot M \}
  = \lla (\mu_{(1)}*\nu_{(1)})\cdot
          \mu_{(2)}\cdot \nu_{(2)} \rra$.
That $(\Sym,\#_{1,1},\eta,\Delta,\epsilon)$ is a Hopf algebra then
follows from the fact that $*$ is Frobenius Laplace and
Proposition~\ref{prop-FrobeniusBialgebra} with a recursive antipode
$\antip_\#$ or from Theorem~\ref{thm-HashAsHopf}.
\qed

The step using multiply often the Laplace expansions is usually formulated
as a separate
\mybenv{Lemma}\textbf{[Cummins]}
\begin{align}\label{eq-cummins}
(A\cdot B)*(C\cdot D)
  &= (A_{(1)}*C_{(1)})\cdot
     (A_{(2)}*D_{(1)})\cdot
     (B_{(1)}*C_{(2)})\cdot
     (B_{(2)}*D_{(2)})
\end{align}
\vskip-1ex
\myeenv

\noindent
\textbf{Proof:}
This was first proved by Cummins~\cite{cummins:1988a} but follows from
the Laplace expansions \eqref{eq-LaplaceLeft} and \eqref{eq-LaplaceRight}
by noting that
\begin{align}
  (A\cdot B)*(C\cdot D)
    &=((A\cdot B)_{(1)}*C)\cdot((A\cdot B)_{(2)}*D)
     =((A_{(1)}\cdot B_{(1)})*C)\cdot((A_{(2)}\cdot B_{(2)})*D)
      \nn
    &=((A_{(1)}*C_{(1)})\cdot(B_{(1)}*C_{(2)}))\cdot
                      ((A_{(2)}*D_{(1)})\cdot(B_{(2)}*D_{(2)}))\,
\end{align}
as required.
\qed

Hence it turns out that the inner product of Thibon characters is realized
by the hash product. Since $\sfa=*$ is Frobenius Laplace this implies
that the underlying branching is actually a Hopf algebra isomorphism.
The reason behind this is, that the branching scheme of Thibon
characters is that of a multiplicative formal group law, as shown in
Appendix~\ref{formalgroups}. This singles the Thibon characters out
on the same footing as the $\GL$ characters for $\CGL=\Sym$. Moreover,
we have seen that the inner product alone acting on the coalgebra
$(\Sym,\comul,\epsilon)$ is only a bialgebra and cannot be Hopf as no
antipode exists in this case. We find hence that the transformation
from the outer product to the hash product is given by a transformation
of formal group laws.

%-----------------------------------------------------------------------
\subsection{Murnaghan--Littlewood stable symmetric group characters}
%-----------------------------------------------------------------------
For the next application we need more notation, especially some reduced
characters of the symmetric groups, enabling one to get rid of the
$n$ dependence of the Kronecker coefficients allowing then an inductive
limit and the definition of universal $S_\infty$ characters. We will
take the liberty of introducing more notation than strictly necessary
so as to allow for vertex operators and some nice graphical representations
of them.
%-----------------------------------------------------------------------
\subsubsection{Stable symmetric group characters}\label{ssubsect:StableSymCharacters}

Studying symmetric groups $S_n$ yields \emph{different} $n$-de\-pen\-dent
products of characters for each $n$. It is highly desirable to remove
this $n$ dependence and this can be achieved using
\emph{reduced characters} introduced by Murnaghan and
Littlewood~\cite{murnaghan:1938a,littlewood:1958a,littlewood:1958b},
see also~\cite{butler:king:1973a,thibon:1991a}. We also will have
occasion to use intermediate \emph{stable permutation characters}, which
we have baptized Thibon characters~\cite{scharf:thibon:wybourne:1993a} above.
To be able to define such reduced characters, we need to allow
nonstandard partitions, that is actually compositions of non negative integers. 
A composition $\Theta$
of an integer $n$ into $p$ parts is a list of non negative integers 
$[\theta_1,\ldots,\theta_p]$ such that $\sum_i \theta_i=n$.
To any composition $\Theta$ one assigns a partition $\vartheta$ by
reordering the parts. 
However, this reordering keeps a sign information
obtainable for example from the definition of the Schur polynomials
in a determinantal (Jacobi-Trudi) form, see~\cite{macdonald:1979a}.
These operators are called raising operators and are defined as
\begin{align}
  R_{i,i+1} [\theta_1,\theta_2,\ldots,\theta_k]
    &= -[\theta_1,\theta_2,\ldots,\theta_{i+1}-1,\theta_i+1,\ldots,\theta_k]
    &&\text{such that}&& \Theta \mapsto \vartheta 
\end{align}
and $\vartheta$ is the unique partition associated to the composition $\Theta$.
Furthermore one demands that a character of the resulting partition is
0 if one gets a trailing negative part, or if one encounters
compositions such that $R_{i,i+1}(\Theta)=-\Theta$, for example
$R_{1,2}[1,2]=-[1,2]$.

Given a tensor irreducible representation $S^\lambda$ of $S_n$ with
$\lambda\vdash n$, it has via the Frobenius characteristic map the
$n$-dependent character $s_\lambda$. The reduced notation removes the
first row and is written as
$\la\lambda\ra=\la\lambda_2,\ldots,\lambda_k\ra$. The character
$s_{\lambda}=\{\lambda\} = \{\lambda_{1},\lambda_{2},\ldots,\lambda_{k}\}$
may be written in reduced notation as
$\la\mu\ra = \{\lambda_{2},\ldots,\lambda_{k}\}$, where the reduced
partition $\mu=(\lambda_{2},\ldots,\lambda_{k})$
has weight $\vert\mu\vert=\vert\lambda\vert-\lambda_{1}$ so that
$\lambda_{1}=n-\vert\mu\vert$. This may be used to recover
$\{\lambda\}=\{n-\vert\mu\vert,\lambda_{2},\ldots,\lambda_{k}\}$ from
$\la\mu\ra$ for any given $n$. The resulting partitions may be
non standard (and may even be negative) when 
$n-\vert\lambda\vert<\lambda_2$ and have to be standardized using the 
raising operators. For example $\la 21\ra$ becomes for $n\geq 5$ the
partition $\{n-3,2,1\}$, and for $n=4:$ $\{1,2,1\}=-\{2,1,1\}$, $n=3:$ 
$\{0,2,1\}=-\{1^3\}$, $n=2:$ $\{-1,2,1\}=0$, $n=1:$ $\{-2,2,1\}=0$.
Treating all $n$ at the same time, getting thus $n$-independence, is
done by using a formal power series in $\Sym\lsqb z\rsqb${}
to get
\begin{align}
\{\lambda\}&\mapsto \la\lambda_2,\ldots,\lambda_k\ra
&&&
\la\lambda\ra_z &= \sum_{n\in\mathbb{Z}} s_{(n-p,\lambda_2,\ldots,\lambda_k)}\,z^n
\end{align}
where $p=\lambda_{2}+\cdots+\lambda_{k}$ and for later use a formal
parameter $z$ has been introduced. We also note that due to standardization
of the partitions one has $s_{(n-p,\lambda_2,\ldots,\lambda_k)}=0$
for $n-p\ll 0$.
%-----------------------------------------------------------------------
\subsubsection{Vertex operators}

Using the series $L(z)$~\eqref{eq-Lser} and $M(z)$~\eqref{eq-Mser} one
can define~\cite{zelevinsky:1981b} a \emph{Bernstein vertex operator}
$V(z)=M(z)L^\perp(\oz)$ in $\End(\Sym)\lsqb z\rsqb$. Recall that the
adjoint operator $L^\perp(\oz)$, with $\oz=1/z$, acts by skewing
\begin{align}
L^\perp(\oz)(s_\mu)
  &= \sum_{n\geq 0} (-1)^n\,e^\perp_n(s_\mu) \oz^{n}
   = \sum_{n\geq 0} \la (-1)^n\,e_n\mid s_{\mu_{(1)}}\ra s_{\mu_{(2)}} \oz^{n}
   =s_{\mu/L(\oz)}
\end{align}
Defining the Thibon characters~\eqref{eq-ThibonCharacter} we used a similar
technique with a multiplicative operator $\lla \lambda\rra=\{\lambda M(1)\}$.
Introducing also a formal parameter and combining these two operators
gives the Bernstein vertex operator in the form
\begin{align}
V(z)
   &=M(z)L^\perp(\oz)
    =\exp\left\{\sum_{n\geq 1} z^k\frac{p_n}{n}\right\}
     \exp\left\{-\sum_{n\geq 1} \oz^k\frac{\partial}{\partial p_n}\right\}
\end{align}
where $n\partial/\partial p_n=p^\perp_n$ was used in the power sum basis
of $\Sym$ (over the rationals)~\cite{macdonald:1979a}, and the expansions
of $M$ and $L$ into this basis. We will have no need to use these explicit
forms, we just remark, that Schur polynomials can be obtained as coefficients
of the action of Bernstein vertex operators on the trivial (vacuum) Schur
function $s_{(0)}=1$. Defining $[z^\lambda]=[z_1^{\lambda_1}\ldots z_l^{\lambda_l}]$
as the operator extracting the coefficient of the monomial $z^{\lambda}$ from an
expression in $\Sym\lsqb z_1,\ldots,z_l\rsqb$ (often written as a contour integration)
we get for $\lambda=(\lambda_1,\ldots,\lambda_l)$
\begin{align}
  s_\lambda
    &= [z^\lambda ]\, V(z_{l})V(z_{l-1})\ldots V(z_1)\, s_{(0)}
\end{align}
In fact each application of a vertex operator adds a part to the
partition. The reduced characters are hence a special case of this process
\begin{align}
\la \mu\ra_{z}
  &= V(z)\,s_\mu
   =\lla L^\perp(\overline{z})\mu\rra
   = \{M(z)\,L^\perp(\overline{z})\,\mu \}
\end{align}
while the Thibon characters employ only the multiplicative operator $M(z)$,
specialized in our previous discussion to $z=1$. A general theory for
vertex operators for reduced groups was developed
in~\cite{fauser:jarvis:king:2010a}.

In graphical terms we get for a vertex operator and Thibon character
(recall $l_{\overline{z}}(x) = \la L(\overline{z}) \mid x\ra $ )
\begin{align}
\begin{pic}
   \node (i1) at (0,1) {};
   \node[rectangle,draw,thick] (v) at (0,0) {$V(z)$ };
   \node (o1) at (0,-1) {};
   \draw[thick] (i1.center) to (v.north);
   \draw[thick] (v.south) to (o1.center);
\end{pic}
&\hskip1ex\cong
\begin{pic}
   \node (i1) at (0.5,1) {};
   \node (u) at (0.5,0.75) {};
   \node[circle,inner sep=2pt,draw,thick,fill=white] (l) at (1,0.25) {};
      \node[xshift=8pt] at (l.east) {$l_{\oz} $};
   \node[circle,inner sep=2pt,draw,thick,fill=white] (M) at (1,-0.25) {};
      \node[xshift=13pt] at (M.east) {$M(z) $};
   \node (m1) at (0,0.25) {};
   \node (m2) at (0,-0.25) {};
   \node (d) at (0.5,-0.75) {};
   \node (o1) at (0.5,-1) {};
\begin{pgfonlayer}{background}
   \draw[thick] (i1.center) to (u.center);
   \draw[thick,out=180,in=90] (u.center) to (m1.center);
   \draw[thick] (m1.center) to (m2.center);
   \draw[thick,out=270,in=180] (m2.center) to (d.center);
   \draw[thick,out=0,in=90] (u.center) to (l.center);
   \draw[thick,out=270,in=0] (M.center) to (d.center);
   \draw[thick] (d.center) to (o1.center);
\end{pgfonlayer}
\end{pic}
;&&&
\begin{pic}
   \node (i1) at (0,1) {};
   \node[rectangle,draw,thick] (v) at (0,0) {$M(z)\cdot$ };
   \node (o1) at (0,-1) {};
   \draw[thick] (i1.center) to (v.north);
   \draw[thick] (v.south) to (o1.center);
\end{pic}
&\hskip1ex\cong
\begin{pic}
   \node (i1) at (0,1) {};
   \node[circle,inner sep=2pt,draw,thick,fill=white] (M) at (1,0.5) {};
      \node[xshift=8pt,yshift=7pt] at (M.east) {$M(z) $};
   \node (m2) at (0,0.5) {};
   \node (d) at (0.5,-0.25) {};
   \node (o1) at (0.5,-1) {};
\begin{pgfonlayer}{background}
   \draw[thick] (i1.center) to (m2.center);
   \draw[thick,out=270,in=180] (m2.center) to (d.center);
   \draw[thick,out=270,in=0] (M.center) to (d.center);
   \draw[thick] (d.center) to (o1.center);
\end{pgfonlayer}
\end{pic}
\hskip-2ex;&&&
\begin{pic}
   \node (i1) at (0,1) {};
   \node[circle,inner sep=2pt,draw,thick,fill=white] (M) at (0.5,0.75) {};
      \node[xshift=8pt,yshift=7pt] at (M.east) {$M(z) $};
   \node[circle,inner sep=2pt,draw,thick,fill=white] (L) at (1.5,0.75) {};
      \node[xshift=8pt,yshift=7pt] at (L.east) {$L(z) $};
   \node (m1) at (1,0.25) {};
   \node (m2) at (0,0.25) {};
   \node (d) at (0.5,-0.25) {};
   \node (o1) at (0.5,-1) {};
\begin{pgfonlayer}{background}
   \draw[thick] (i1.center) to (m2.center);
   \draw[thick,out=270,in=180] (m2.center) to (d.center);
   \draw[thick,out=270,in=180] (M.center) to (m1.center);
   \draw[thick,out=270,in=  0] (L.center) to (m1.center);
   \draw[thick,out=270,in=0] (m1.center) to (d.center);
   \draw[thick] (d.center) to (o1.center);
\end{pgfonlayer}
\end{pic}
&\hskip-2ex\cong
\begin{pic}
   \node (i1) at (0,1) {};
   \node[circle,inner sep=2pt,draw,thick,fill=white] (M) at (1,0.5) {};
      \node[xshift=8pt] at (M.east) {$\eta $};
   \node (m2) at (0,0.5) {};
   \node (d) at (0.5,-0.25) {};
   \node (o1) at (0.5,-1) {};
\begin{pgfonlayer}{background}
   \draw[thick] (i1.center) to (m2.center);
   \draw[thick,out=270,in=180] (m2.center) to (d.center);
   \draw[thick,out=270,in=0] (M.center) to (d.center);
   \draw[thick] (d.center) to (o1.center);
\end{pgfonlayer}
\end{pic}
\cong
\begin{pic}
   \node (i1) at (0,1) {};
   \node (o1) at (0,-1) {};
   \draw[thick] (i1.center) to (o1.center);
\end{pic}
\end{align}
The last tangle equality shows how from the inverse series $M(z)L(z)=1$
one gets an inverse of the branching $\{\lambda\} = \lla\lambda L(1)\rra{}
= \{\lambda L(z)M(z)\}$.
As the skewing operator does not commute with the multiplication, the
situation for vertex operators is more complicated. One obtains the
commutation relation
\begin{align}
L^\perp(z)M(w)
  &= (1-zw)M(w)L^\perp(z){}
;&&&
\begin{pic}
   \node (i1) at (0,1) {};
   \node (u) at (0.5,0.25) {};
   \node[circle,inner sep=2pt,draw,thick,fill=white] (M) at (1,0.75) {};
      \node[xshift=8pt,yshift=7pt] at (M.east) {$M(w) $};
   \node[circle,inner sep=2pt,draw,thick,fill=white] (l) at (1,-0.75) {};
      \node[xshift=8pt,yshift=-7pt] at (l.east) {$l_{z} $};
   \node (m1) at (0,0) {};
   \node (d) at (0.5,-0.25) {};
   \node (o1) at (0,-1) {};
\begin{pgfonlayer}{background}
   \draw[thick,out=270,in=180] (i1.center) to (u.center);
   \draw[thick,out=270,in=0] (M.center) to (u.center);
   \draw[thick] (u.center) to (d.center);
   \draw[thick,out=0,in=90] (d.center) to (l.center);
   \draw[thick,out=180,in=90] (d.center) to (o1.center);
\end{pgfonlayer}
\end{pic}
\cong\hskip1ex
\begin{pic}
   \node (i1) at (0.5,1) {};
   \node (u) at (0.5,0.75) {};
   \node[circle,inner sep=2pt,draw,thick,fill=white] (l) at (1,0.25) {};
      \node[xshift=-14pt] at (l.east) {$l_{z} $};
   \node[circle,inner sep=2pt,draw,thick,fill=white] (M) at (1,-0.25) {};
      \node[xshift=-18pt] at (M.east) {\tiny $M(w) $};
   \node (m1) at (0,0.25) {};
   \node (m2) at (0,-0.25) {};
   \node (d) at (0.5,-0.75) {};
   \node (o1) at (0.5,-1) {};
   \node[circle,inner sep=2pt,draw,thick,fill=white] (M2) at (1.5,0.25) {};
      \node[xshift=8pt,yshift=7pt] at (M2.east) {\tiny $M(w) $};
   \node[circle,inner sep=2pt,draw,thick,fill=white] (l2) at (1.5,-0.25) {};
      \node[xshift=8pt,yshift=-7pt] at (l2.east) {$l_{z} $};
\begin{pgfonlayer}{background}
   \draw[thick] (i1.center) to (u.center);
   \draw[thick,out=180,in=90] (u.center) to (m1.center);
   \draw[thick] (m1.center) to (m2.center);
   \draw[thick,out=270,in=180] (m2.center) to (d.center);
   \draw[thick,out=0,in=90] (u.center) to (l.center);
   \draw[thick,out=270,in=0] (M.center) to (d.center);
   \draw[thick] (d.center) to (o1.center);
   \draw[thick] (M2.center) to (l2.center);
\end{pgfonlayer}
\end{pic}  
\end{align}
which is a direct consequence of the bialgebra
law~\eqref{grph-bialgebra-and-antipode}, $M(w)$ and $L(z)$ being group
like, and the evaluation $\la L(z)\mid M(w)\ra = (1-zw)$ as only $s_{(0)}$
and $s_{(1)}$ terms survive. It is easy to derive many results about
vertex operators using the graphical language. Some care is needed as
evaluations may turn out to be infinite, as for example
$\la M(1)\mid M(1)\ra = \sum_{n\geq0} 1$ and needs regularization,
or may vanish as for $\la L(1) \mid M(1)\ra = 0$.

Vertex operators emerged prominently
in~\cite{sato:miwa:jimbo:1997a,data:kashiwara:miwa:1981a} and following
papers. In the setting for symmetric functions see for
example~\cite{jing:1991a,jing:1991b,jarvis:yung:1992a,jarvis:yung:1993a,%
jarvis:yung:1994a,baker:1996a} and for non-classical subgroups of $\GL(N)${}
see~\cite{fauser:jarvis:king:2010a}.
As our emphasis is on providing examples for hash products we return
to the decomposition of reduced character products.

%-----------------------------------------------------------------------
\subsubsection{Murnaghan-Littlewood inner branching}

Dealing with reduced $S_\infty$ characters, a natural question is to
ask for their now $n$ independent inner multiplication and Kronecker
coefficients $\widetilde{g}^\lambda_{\mu,\nu}$. The classical
result is
\mybenv{Theorem}\textbf{[Murnaghan-Littlewood]}
The inner product of reduced symmetric group characters is given by
the recursive formula
\begin{align}
\la \mu\ra * \la \nu\ra
  &= \sum_{\lambda} \widetilde{g}^\lambda_{\mu,\nu}\la \lambda\ra
   = \sum_{\alpha,\beta,\zeta} \la \mu/(\alpha\zeta)\,\nu/(\beta\zeta)
     (\alpha * \beta )\ra
\end{align}
The $n$-dependent inner product $(\alpha*\beta)$ is computed by
recursion over this formula, as it works on lower weight terms.
\myeenv
Using higher derived hash products we can provide this multiplication
by a threefold hash product.
\mybenv{Theorem:}
The Murnaghan-Littlewood product formula for stable symmetric group
characters is given by the higher derived hash product
\begin{align}
\la x\ra * \la y\ra
  &=\la x \,\#_{m,1,1}\, y\ra
   = \la M\mid x_{(1)} * y_{(1)}\ra
     \,\la(x_{(2)} * y_{(2)})\cdot (x_{(3)}\cdot y_{(3)})\ra
\end{align}
where $\#_{m,1,1}=m\circ*\conv*\conv\mul$ and it will be recalled that
$m(x)=\la M\mid x\ra$. Hence this is a further deformation of the hash
product $\#_{1,1}$ for Thibon characters.
\myeenv

\noindent
\textbf{Proof:}
We can give a direct argument to show that the higher derived hash
evaluates to the Murnaghan Littlewood formula as follows:
\begin{align}
\la \mu\ra \,\#_{m,1,1}\, \la \nu\ra
  &= \la M\mid \mu_{(1)} * \nu_{(1)}\ra
       (\mu_{(2)} * \nu_{(2)})(\mu_{(3)}\nu_{(3)})
     \nn
  &= \sum_{\alpha,\beta,\zeta,\sigma}
     \la M\mid\beta * \sigma\ra
     \,\la
     \left(\alpha/\beta\right)*\left(\zeta/\sigma\right)
     \left(\mu/\alpha\right)\left(\nu/\zeta\right)\ra
    \nn
  &= \sum_{\alpha,\beta,\zeta}
    \la \left(\alpha/\beta\right)*\left(\zeta/\beta\right)
        \left(\mu/\alpha\right)*\left(\nu/\zeta\right) \ra
    \nn
  &= \sum_{\alpha,\beta,\zeta}
      \la \left(\mu/(\alpha\zeta)\right)\left(\nu/(\beta\zeta)\right)
         \left(\alpha * \beta\right) \ra
\end{align}
showing the equivalence.
\qed

It may be instructive to see how this result can be derived from the
definition of reduced and Thibon characters.
\mybenv{Lemma}\label{lem-perpAction}
The outer comultiplication of a skew $L^\perp \{\mu\}$ expands as follows
\begin{align}
\Delta (L^\perp \{\mu\})
   =(L^\perp \mu)_{(1)}\ot(L^\perp \mu)_{(2)}
  &= (L^\perp \{\mu\}_{(1)})\ot (\{\mu\}_{(2)})
   = (\{\mu\}_{(1)}) \ot (L^\perp \{\mu\}_{(2)})
\end{align}
That is the the derivation acts either left or right. Furthermore the
invertibility $LM=1$ induces
\begin{align}
M^\perp L^\perp \{\mu\}
  &= (ML)^\perp \{\mu\} = 1^\perp \{\mu\} = \{\mu\}
\end{align}
This allows us to derive the action of a skew on an inner product
\begin{align}
M^\perp ( A * B )
  &= \la M \mid A_{(1)} * B_{(1)}\ra\, (A_{(2)} * B_{(2)})
  \nonumber \\
  &= \la A_{(1)} \mid B_{(1)} \ra \, (A_{(2)} * B_{(2)})
\end{align}
\myeenv
\noindent
\textbf{Proof:} This involves a short graphical calculation using the
Frobenius Hopf mixed bialgebra law~\eqref{grph-mixedBialgebra} and the
fact that $\la M \mid x*y\ra = \la x\mid y\ra$. Alternatively, to see
the first part, for all $x,y,z$ we have
\begin{align}
   \la\, x\ot y \mid \Delta(L^\perp z) \,\ra 
   &= \la\, x\, y \mid L^\perp z \,\ra 
    = \la\, L\, x\, y \mid z\,\ra 
    = \la\, (Lx)\,y \mid z \,\ra \nn
   &= \la\, (Lx)\ot y \mid \Delta(z)\,\ra 
    = \la\, (Lx)\ot y \mid z_{(1)}\ot z_{(2)} \,\ra \nn
   &= \la\, x\ot y \mid (L^\perp x_{(1)}) \ot z_{(2)} \,\ra
\end{align}
Similarly for the final part
\begin{align}
  \la\, z \mid M^\perp (x*y)\,\ra 
  &= \la\, M\,z \mid x*y \,\ra{}
   = \la\, \delta(M\,z) \mid x\ot y \,\ra \nn
  &= \la\, \delta(M)\,\delta(z) \mid x\ot y \,\ra
   = \sum_\zeta \la\, (\zeta\ot\zeta)\, \delta(z) \mid x\ot y\,\ra \nn
  &= \sum_\zeta \la\, \delta(z) \mid (x/\zeta)\ot(y/\zeta) \,\ra
   = \sum_\zeta \la\, z \mid (x/\zeta)*(y/\zeta) \,\ra \nn
  &= \sum_{\zeta,\eta} \la\, z \mid (\la\,\zeta\mid\eta\,\ra){}
     \ (x/\zeta)*(y/\eta) \,\ra\nn
  &= \la\, z\mid (\la\, x_{(1)}\mid y_{(1)}\,\ra)\ (x_{(2)}*y_{(2)}) \,\ra
\end{align}
Since these are true for all $x,y,z$ the required result follows.
\qed

We can now obtain the higher derived hash product using intermediate
Thibon characters.
\mybenv{Corollary}
The tensor product decomposition of stable symmetric group characters is
obtained from the decomposition of Thibon characters,
Theorem~\ref{thm-innerThibon}, by introducing a further skew by
$L^\perp$ and hence by a vertex operator expansion
\begin{align}
\la \mu\ra\,*\,\la \nu\ra
  &= \lla L^\perp \mu\rra \#_{1,1} \lla L^\perp \nu \rra
   = \lla(L^\perp (\mu\,\#_{\mathsf{m},1,1}\, \nu)\rra
  \nonumber \\
  &= \la \mu\,\#_{\mathsf{m},1,1}\, \nu\ra
\end{align}
defining the higher derived hash product.
\myeenv

\noindent\textbf{Proof:}
We obtain the result either graphically (as the reader may verify) or
algebraically, using the facts from Lemma~\ref{lem-perpAction}, as
follows:
\begin{align}
\la \mu\ra\,*\,\la \nu\ra
  &= \lla L^\perp \mu\rra * \lla L^\perp \nu \rra
  \nonumber \\
  &= \lla (L^\perp \mu)\,\#_{1,1}\, (L^\perp \nu) \rra
   = \lla (L^\perp \mu)_{(1)}*(L^\perp\nu)_{(1)}
         \big(\, (L^\perp \mu)_{(2)} (L^\perp \nu)_{(2)} \,\big)\rra
  \nonumber \\
  &=\lla (\mu_{(1)} * \nu_{(1)})\,
         (L^\perp \mu_{(2)})\, (L^\perp \nu_{(2)}) \rra 
   \nonumber \\      
  &= \lla L^\perp M^\perp \big[\, (\mu_{(1)} * \nu_{(1)})
          (L^\perp \mu_{(2)}) (L^\perp \nu_{(2)}) \,\big] \rra
  \nonumber \\
  &= \lla L^\perp \big[\, M^\perp (\mu_{(1)} * \nu_{(1)})
            (M^\perp L^\perp \mu_{(2)})
            (M^\perp L^\perp \nu_{(2)}) \,\big] \rra
  \nonumber \\
  &= \lla L^\perp\big[\, \la\mu_{(1)}\mid \nu_{(1)}\ra
            (\mu_{(2)}*\nu_{(2)})\, \mu_{(3)}\nu_{(3)} \,\big] \rra
  \nonumber \\
  &= \la \mu\,\#_{m,1,1}\,\nu\ra
\end{align}
We used again the fact that $M$ and $L$ are group like series.
\qed

The proof shows, that an additional application of a Schur function series
$L, M$ or a skew Schur function series $L^\perp, M^\perp$ results in going
up from a hash to a higher hash product by one step.
%-----------------------------------------------------------------------
\section{Closing remarks}\label{sec-Misc}
%-----------------------------------------------------------------------

Having demonstrated the plethora of examples of group character
decompositions of our deformation theory of (character) multiplications
parameterized by Laplace pairings, it remains to give some pointers on
shortcomings and possible further generalizations of the method.
Moreover, the notion of Laplace pairings can be generalized at the cost
of losing several properties of the obtained deformations and it is
clear that not every 2-cocycle needs to be Laplace, but will still
produce associative deformed multiplications.

An obvious comparison has to be made to our
work~\cite{fauser:jarvis:king:wybourne:2005a}, where we used plethystic
series of Schur functions to describe a very general method to obtain
character decompositions for $\GL$-subgroups defined by polynomial
identities. We encountered plethystic pairings $r_{\pi}$, which we have
also employed in~\cite{fauser:jarvis:king:2010c}, which is not only
employing convolution, but also composition (plethysm) of symmetric
functions. The cases where we use the Schur functions $s_{(2)}$
and $s_{(1,1)}$ produce the orthogonal and symplectic branchings,
as only trivial plethysms are involved (the cut comultiplication
produces $s_{(1)}\ot s_{(1)}$ in both cases, and $s_{(1)}(X)=X$ produces
the trivial composition $\Id : X\rightarrow X$). Moving to nonclassical
groups leads to $r_{\pi}$ pairings which can be shown to be not
Laplace by direct computation. The problem is that plethysm (composition)
is a non-linear operation and induces certain scalings related to loop
operators, see Appendix~\ref{FGLloop}. However, the treatment of the
symmetric group reduced characters above, which emerges from a higher
order deformation of symmetry type $(2)$ and $(3)$, shows that not all
cases are beyond the reach of the present theory using the Laplace
subgroup to parameterize deformations.

The deformation theory presented here relies on the fact that the ambient
Hopf algebra is built over free modules. Restricting the number of
indeterminates to a finite $N$ induces syzygies and the corresponding
\emph{modification rules} have to be found and applied. These are
known for classical groups, but for the nonclassical groups obtained
by higher deformations must be developed, due to the present lack of
a general method, case-by-case. A general theory of modification
rules is as of now not available and poses a big challenge.

A way to go further is by studying `loop operators' (related to Adams
operations), that is operators of the form
$[r]=\mul^{(r)}\circ\comul^{(r)}$ or similar operators for the Laplace
pairing or mixed forms. The loop operators are related to Eulerian and Lie
idempotents and in Appendix~\ref{FGLloop} we exemplify that they also inject
the ring of integers into the respective formal group. Moreover these
idempotents are related to Adams operators, and hence to the scaling
properties in use in plethystic branchings. Without entering into details,
it is sufficient to remark that one can use such loop operators to
enforce \emph{modified} Laplace type properties, and even make the 
associative multiplication $\mul$ of $\HB$ `Laplace'. The resulting
deformations are, however, in general no longer associative.

It may be remarked, that the algebras of outer and inner plethysm as
defined by Littlewood, see for
example~\cite{littlewood:1940a,littlewood:1958a,butler:1970a,butler:king:1973a,king:1974a},
has in its non linear part a structure which is very akin to the Laplace
expansion laws~(\ref{eq-LaplaceRight},\ref{eq-LaplaceLeft}). The difference
is, that in the plethystic expansion laws the multiplications and
comultiplications in use are different, hence not from the same self dual
ambient Hopf algebra. This is reminiscent of the fact that bialgebra pairings,
see Definition~\ref{def-BialgebraPairing}, actually mediate between
\emph{different} bialgebras in a matched pair of bialgebras or Hopf algebras.
These expansions are very effective tools to compute such plethysms, as
we showed in~\cite{fauser:jarvis:king:2010a}. An expansion of our
present theory to such a setting is highly desirable.

We mentioned the results of Aguiar et
al.~\cite{aguiar:ferrer:moreira:2004b,aguiar:ferrer:moreira:2004a}
producing a noncommutative version of the hash product. Beside the
problematic notion of plethysm in this setting, our theory of higher
derived hash products should generalize along similar lines. This
is also supported by the original development of Rota and Stein, who
worked in a general supersymmetric, that is graded commutative, setting.

Further opportunities to extend the present theory include the usage of
non polynomial formal group laws. These FGLs would lead to an
\emph{infinite convolution product} of Laplace pairings having its own
intricacies. However, along with the additive versus multiplicative
formal group analogy these deformations should be Hopf, keeping the
original comultiplication of the Rota Stein development.

Rota and Stein showed that the Hopf algebra of symmetric functions is
generated by a single letter in a two stage process~\cite{rota:stein:1994a}
(our alphabet $X$). In the second part of that work~\cite{rota:stein:1994b}
they show how to deal with vector symmetric functions (that is multiple
alphabets) which are related to plane partitions. From a formal point of
view our theory applies to this case as well, but we have not yet
investigated it thoroughly.

We have seen the that the \emph{group like} series $M$ and $L$ play a
central role in the theory of Bernstein (vertex) operators and stable
symmetric group characters. Using power sum plethysms, one may produce
other group like series $G_{n}(t)=M[p_{n}](t)$. The develoment of the
branching alters, but our general theory should cope with this.

Our theory is indeed versatile enough to deal with deformed symmetric
functions, as we have not put any restriction on the involved (Frobenius)
Laplace pairings. Using a $(q,t)$ parameterized $z_{\mu}(q,t)$ factor in
the inner product of power sum symmetric functions, will for example
produce a $(q,t)$-parameterized Schur-Hall scalar product~\cite{macdonald:1979a}.
This is actually just applying a different $(q,t)$ dependent map for
defining the derived pairing. In such a manner one can obtain for
example Jack, zonal, Hall-Littlewood and Macdonald symmetric functions
and our formal theory does not change. This renders the theory presented
here rather flexible and worth studying further.
%-----------------------------------------------------------------------
\medskip
\begin{appendix}
\renewcommand{\theequation}{\Alph{section}-\arabic{equation}}

\medskip
\section{Relation to formal group laws\label{formalgroups}}

In this appendix we provide some hints on the relationship of the theory
of formal group laws to our present development. This sheds some light on the
possibility of using the framework developed here in a more general setting,
especially it may help to produce similar results for fields of finite
characteristics or may provide a guide to the Hopf algebra approach in
the case of non polynomial formal group laws. For a treatment of formal
groups see for
example~\cite{childs:2000a,hazewinkel:1978a,yanagihara:1977a,lazard:1975a,froehlich:1968a}.

\subsection{Definition and basics} We start by giving a definition of
formal group laws.

\mybenv{Definition}\label{fgl}
Let $\mathsf{R}\lsqb X,Y\rsqb$ be the ring of formal power series in a
commutative set of indeterminants $X,Y$ and let $\mathsf{R}\lsqb X\rsqb =
\mathsf{R}+ \mathsf{R}^+\lsqb X\rsqb$ the decomposition into the (valuation)
ring $\mathsf{R}$ and the augmentation ideal $\mathsf{R}^+\lsqb X\rsqb$,
furthermore let $\lambda(X)\in \mathsf{R}^+\lsqb X\rsqb$.
An $n$-dimensional \emph{formal group law} (FGL) $F(X,Y)$ is an $n$-tuple
of formal power series $F_i(X,Y)$, $i\in \{1,\ldots,n\}$ over alphabets
$X,Y,\ldots$ of length (size) $n$
\begin{align}
F(X,Y) &=
  (F_1(X,Y),F_2(X,Y),\ldots,F_n(X,Y))
\end{align}
such that
\begin{align}
i)\,\,\,&~~ F_i(X,Y) = X + Y \mod \deg 2 \nonumber \\
ii)\,&~~ F_i(F(X,Y),Z) = F_i(X,F(Y,Z)) &&\textrm{associativity} \nonumber \\
iii)&~~ F(X,0)=X,~~F(0,Y)=Y  &&\textrm{identity} \nonumber \\
iv)\,&~~ F(X,\lambda(X))= 0 &&\textrm{inverse} \nonumber \\
v)\,\,&~~ F_i(X,Y)=F_i(Y,X) &&\textrm{commutativity}
\end{align}
We write alternatively $X +_F Y = F(X,Y)$ introducing the $+_F$-addition.
\myeenv

\mybenv{Examples}
Generic formal group laws are not polynomial. However for our present
situation the very simple polynomial formal groups are of interest. We
restrict ourselves to 1-dimensional FGLs for notational simplicity.
\begin{itemize}
\item
$\mathbb{G}_a(X,Y) = X +_a Y = X + Y$ the additive formal group law. This
group describes the space of $K$-points $m\in P(\mathbb{G}_a,K)$ where $m$ are
the maximal ideals of $\mathsf{R}$. The antipode $\lambda(X)$ has to fulfil
$X+\lambda(X)=0$, hence we see that $\lambda(X)=-X$.
\item
$\mathbb{G}_m(X,Y) = X +_m Y = X + Y + XY$ the multiplicative formal group
law. This group describes the space
$1+m \in U_1(\mathsf{R}) \cong P(\mathbb{G}_m,K)$ the multiplicative
group of principal units. The antipode $\lambda(X)$ has to fulfil
$X +\lambda(X) +X\lambda(X)=0$, hence one gets $\lambda(X)= -X/(1+X)$.
A slight generalization of the multiplicative group is obtained by
introducing a parameter $b$, an invertible element in $\mathsf{R}$
(group of units). The FGL is then given as
$\mathbb{G}^b_m(X,Y) = X +_m^b Y = X + Y + bXY$, with antipode
$\lambda(X)= -X/(1+bX)$.
\item
Let $S=\mathbb{Z}[\{c_{i,j} \mid i,j\ge 1\}]$ with $\{c_{ij,} \mid i,j\geq 1\}$
an infinite set of formal parameters. Define the universal formal group law
as  $F(X,Y) = X + Y + \sum_{i,j\ge1} c_{i,j}X^iY^j$ such that the relations
in definition \ref{fgl} are fulfilled. Specifying the parameters accordingly
allows the study of formal group laws over finite fields etc.
\end{itemize}
\myeenv

\mybenv{Definition}\label{fglMorphism}
Let $F,G$ be an $n$ and an $m$-dimensional FGL respectively. A homomorphism
of FGLs $\phi$ is an $m$-tuple of formal power series
$\Phi(X) = (\Phi_1(X),\ldots,\Phi_m(X))$ in $\mathsf{R}^+\lsqb X\rsqb^m$
such that
\begin{align}
\Phi_i(F(X,Y)) = G_i(\Phi(X),\Phi(Y))
\end{align}
\myeenv
$\Phi$ is an isomorphism if $n=m$ and if there exists a FGL morphism $\Psi$
such that $\Psi(\Phi(X))=X$ and $\Phi(\Psi(X))=X$. Homomorphisms of formal
group laws map points of $F$ to points of $G$,
$\Phi : P(F) \longrightarrow P(G)$. This leads to the formula
\begin{align}
\Phi(X +_F Y) &= \Phi(F(X,Y)) = G(\Phi(X),\Phi(Y)) \nonumber \\
  &= \Phi(X) +_G \Phi(Y)
\end{align}

\mybenv{Example}
Let $f\in \mathsf{R}^+\lsqb X\rsqb$ with compositional inverse $\overline{f}$,
and $\overline{f}\circ f = \Id$. The FGL $F$ is isomorphic to the additive
formal group if and only if there exists an $f$ as above with
\begin{align}
F(X,Y) &= \overline{f}(\mathbb{G}_a(f(X),f(Y)))
        = \overline{f}(f(X) +_a f(Y))
\end{align}
For example if $\mathbb{Q}\subseteq \mathsf{R}$ then define
\begin{align}
f(x) &= \ln(1+x) = -\sum_{n\ge 1} \frac{(-1)^nx^n}{n} \nonumber \\
\overline{f}(x) &= e^x-1 = \sum_{n\ge 1} \frac{x^n}{n!}
\end{align}
Then follows $f : \mathbb{G}_m \rightarrow \mathbb{G}_a$ via
\begin{align}
\mathbb{G}_m(x,y) &= e^{\ln(1+x)+\ln(1+y)}-1 = (1+x)(1+y)-1
  = x + y + xy = x +_m y
\end{align}
Now~\cite{froehlich:1968a} establishes that in characteristic 0 every
$n$-dimensional formal group law $F$ is isomorphic to the $n$-dimensional
additive formal group $\mathbb{G}_a$. The morphism $\ell_F : F \longrightarrow
\mathbb{G}_a$ is called logarithm. If $\partial_x\ell_F(x)\vert_{x=0}=1$,
it then follows that $\ell_F(x) = x \mod \deg 2$ and $\ell_F$ is uniquely defined.
The multiplicative group is, however, not isomorphic to the additive one
over a valuation ring $\mathsf{R}$ of a local field $K$ since no such logarithm
exists.
\myeenv

\subsection{Relation to Hopf algebras}

\subsubsection{Definitions and identification of coproducts}

Since we are dealing with power series, not with polynomials, we need to
complete the tensor product due to $\mathsf{R}\lsqb X,Y\rsqb \cong
\mathsf{R}\lsqb X\rsqb \hat{\otimes} \mathsf{R}\lsqb Y\rsqb$,
see~\cite{froehlich:1968a}. We set $X\cong X\otimes 1$ and $Y\cong 1\otimes X$
omitting hats on tensors usually. We extend without further notice the
multiplication and especially the comultiplication map
\begin{align}
\Delta :&\mathsf{R}\lsqb X\rsqb \longrightarrow
\mathsf{R}\lsqb X\rsqb \hat{\otimes} \mathsf{R}\lsqb X\rsqb{}
\subset \mathsf{R}\lsqb X\rsqb \hat{\otimes} \mathsf{R}\lsqb Y\rsqb
\end{align}
We can use this fact to show that every FGL defines a coproduct.
\mybenv{Definition}
To every $n$-dimensional formal group law $F$ we define a comultiplication
$\Delta_F$, a counit $\epsilon$ and an antipode $\lambda$ as follows
\begin{align}
\Delta_F(x_i) &= F_i(X\otimes 1, 1\otimes X) \cong F_i(X,Y)
  \nonumber \\
\epsilon(x_i) &= 0,\quad
  \epsilon : \mathsf{R}\lsqb X\rsqb \longrightarrow \mathsf{R}
  \nonumber \\
\lambda(x_i) &= [-1]_i(x) = -x_i \mod \deg 2
\end{align}
\myeenv
For example we have $\Delta_m(x) = \mathbb{G}_m(x,y) = x + y + xy
\cong x\otimes 1 +1\otimes x + x\otimes x$ for the multiplicative group.
It is hence clear, that a morphism of FGLs, Definition~\ref{fglMorphism},
induces a morphism of coproducts induced by the involved FGLs via
\begin{align}
\Delta_F(f)(x) &= f(\Delta_G(x_1),\ldots,\Delta_G(x_n))
\end{align}
In our case, since $\Symm[X]$ is self dual, we induce also multiplications
by this process. We want to stress the following observation:
\begin{align}
s_\lambda(\mathbb{G}_a(X,Y)) =
s_\lambda(X+Y) &= s_{\lambda_{(1)}}(X)s_{\lambda_{(2)}}(Y)
   \cong \Delta s_\lambda \nonumber \\
s_\lambda(\mathbb{G}_m(X,Y)) =
s_\lambda(X+Y+XY)
  &= s_{\lambda_{(1)}}(X+Y)s_{\lambda_{(2)}}(XY) \nonumber \\
  &= s_{\lambda_{(1)}}(X)s_{\lambda_{(2)[1]}}(X)
     s_{\lambda_{(3)}}(Y)s_{\lambda_{(2)[2]}}(Y) \nonumber \\
  &\cong \Delta_{\#_{1,1}} s_\lambda
\end{align}
Here the Schur functions are seen as (polynomial) morphisms of FGLs. We find
that the additive formal group gives rise to the outer coproduct $\Delta$
on $\Symm[X]$, but the (generic) hash product gives rise to the multiplicative
group coproduct $\Delta_m = \Delta_{\#_{1,1}}$. This corresponds to the
branching law which we derived for the Thibon characters $\lla \lambda\rra$
branching $GL(n+m+nm)$ into $GL(n)\times GL(m)$. This relation is rather
fascinating especially for those characters which are no longer polynomial but
depend on Schur function series, as do the Thibon or the Murnaghan-Littlewood
characters, or the dual characters of orthogonal and symplectic
groups~\cite{fauser:jarvis:king:2007c}. Furthermore, note that the inner
comultiplication $\delta$, related to the product $XY$ of alphabets
(and by duality to the inner product $*$ of symmetric functions) is not a
formal group law, and hence \emph{lacks} this nice theory behind it. In
this sense, the work of Aguiar et
al~\cite{aguiar:ferrer:moreira:2004b,aguiar:ferrer:moreira:2004a}
shows that there is an analog of the multiplicative FGL even in the
non-commutative realm.

\subsubsection{The loop operator\label{FGLloop}}

In Section~\ref{sec-Misc} of the paper, we briefly mentioned
`loop operators' $[r]=\mul^{(r-1)}\circ\comul^{(r-1)}$. In this part
of the Appendix we want to show which role these operators play
in the theory of formal groups. This is interesting in terms of finite
fields, where such operators eventually will be nilpotent.

\mybenv{Definition}
There is an injection of the integers $\mathbb{Z}$ into the
endomorphisms ring of the FGL $F$
\begin{align}
 [~] : \mathbb{Z} &\longrightarrow \End(F) \nonumber \\
 [1](X) &= X \nonumber \\
 [2](X) &= F(X,X), \nonumber \\
 [m](X) &= F([m-1](X),X),\quad m>1 \nonumber \\
 [-1](X) &= \lambda(X) \nonumber \\
 [-m](X) &= F([-(m-1)](X),\lambda(X)),\quad m>1
\end{align}
and $[m]\in \End (F)$ for all $m\in\mathbb{Z}$.
\myeenv

\mybenv{Examples}
$\mathbb{G}_a$: For the additive group we find
\begin{align}
 [n](X)=nX=X + \ldots + X,\quad \forall n
\end{align}
is the repetition of alphabets. It is possible to generalize this
map to $\mathbb{Q}$ or even $\mathbb{R}$ or $\mathbb{C}$,
see~\cite{fauser:jarvis:king:2010a}.

$\mathbb{G}_m^b$: In the multiplicative case one obtains easily
\begin{align}
 [n](X)= \frac{1}{b}(1+bX)^n-\frac{1}{b},\quad \forall n>0
\end{align}
showing the dependency of the loop operator $[n]$ on the underlying
FGL (Hopf algebra). For $b=1$ this is the multiplicative group.
\myeenv
\end{appendix}
%-----------------------------------------------------------------------
%%%%%%%%%%%%%%%%%%%%%Bibliography from bibtex database here
%
%{\small
%\bibliography{sql}
%\bibliographystyle{plain}
%\def\topsep{0pt}
%\def\parsep{0pt plus 5pt minus 1pt}
%\def\itemsep{-0.5ex}
%
%\printindex
%}
%-----------------------------------------------------------------------
\end{document}
\eof
%-----------------------------------------------------------------------